\documentclass[12pt,a4paper]{article}
\usepackage{amsmath}
\usepackage{amssymb}
\usepackage{amsxtra}
\usepackage{amscd}
\usepackage{amsthm}
\usepackage{hyperref}
\usepackage[mathscr]{eucal} 
\usepackage{graphicx}
\usepackage{mathtools}
\DeclareMathOperator{\MyProd}{\scalebox{1.4}{$\mathrm{I\kern-0.2ex I}$}}

\theoremstyle{plain}

\newtheorem{sbprop}[subsubsection]{Proposition}

\newtheorem{sblem}[subsubsection]{Lemma}

\theoremstyle{definition}

\newtheorem{para}[subsection]{}

\newtheorem{sbrem}[subsubsection]{Remark}

\newtheorem{sbpara}[subsubsection]{}

\newenvironment{pf}{\proof[\proofname]}{\endproof}
\begin{document}

\title{Height functions for motives}

\author
{Kazuya Kato}

\maketitle

%\numberwithin{equation}{section}
\newcommand{\lr}[1]{\langle#1\rangle}
\newcommand{\ul}[1]{\underline{#1}}
\newcommand{\eq}[2]{\begin{equation}\label{#1}#2 \end{equation}}
\newcommand{\ml}[2]{\begin{multline}\label{#1}#2 \end{multline}}
\newcommand{\ga}[2]{\begin{gather}\label{#1}#2 \end{gather}}
\newcommand{\mc}{\mathcal}
\newcommand{\mb}{\mathbb}
\newcommand{\surj}{\twoheadrightarrow}
\newcommand{\inj}{\hookrightarrow}
\newcommand{\red}{{\rm red}}
\newcommand{\codim}{{\rm codim}}
\newcommand{\rank}{{\rm rank}}
\newcommand{\Pic}{{\rm Pic}}
\newcommand{\Div}{{\rm Div}}
\newcommand{\divi}{{\rm div}}
\newcommand{\Hom}{{\rm Hom}}
\newcommand{\Ext}{{\rm Ext}}
\newcommand{\im}{{\rm im}}
\newcommand{\fil}{{\rm fil}}
\newcommand{\gp}{{\rm gp}}
\newcommand{\Spec}{{\rm Spec}}
\newcommand{\Sing}{{\rm Sing}}
\newcommand{\Char}{{\rm char}}
\newcommand{\Tr}{{\rm Tr}}
\newcommand{\Gal}{{\rm Gal}}
\newcommand{\Min}{{\rm Min}}
\newcommand{\mult}{{\rm mult}}
\newcommand{\Max}{{\rm Max}}
\newcommand{\Alb}{{\rm Alb}}
\newcommand{\gr}{{\rm gr}}
\newcommand{\Ker}{{\rm Ker}}
\newcommand{\Lie}{{\rm Lie}}
\newcommand{\infi}{{\rm inf}}
\newcommand{\et}{{\rm \acute{e}t}}
\newcommand{\pole}{{\rm pole}}
\newcommand{\ti}{\times }
\newcommand{\modu}{{\rm mod}}
\newcommand{\Ab}[1]{{\mathcal A} {\mathit b}/#1}
% Skriptbuchstaben
\newcommand{\sA}{{\mathcal A}}
\newcommand{\sB}{{\mathcal B}}
\newcommand{\sC}{{\mathcal C}}
\newcommand{\sD}{{\mathcal D}}
\newcommand{\sE}{{\mathcal E}}
\newcommand{\sF}{{\mathcal F}}
\newcommand{\sG}{{\mathcal G}}
\newcommand{\sH}{{\mathcal H}}
\newcommand{\sI}{{\mathcal I}}
\newcommand{\sJ}{{\mathcal J}}
\newcommand{\sK}{{\mathcal K}}
\newcommand{\sL}{{\mathcal L}}
\newcommand{\sM}{{\mathcal M}}
\newcommand{\sN}{{\mathcal N}}
\newcommand{\sO}{{\mathcal O}}
\newcommand{\sP}{{\mathcal P}}
\newcommand{\sQ}{{\mathcal Q}}
\newcommand{\sR}{{\mathcal R}}
\newcommand{\sS}{{\mathcal S}}
\newcommand{\sT}{{\mathcal T}}
\newcommand{\sU}{{\mathcal U}}
\newcommand{\sV}{{\mathcal V}}
\newcommand{\sW}{{\mathcal W}}
\newcommand{\sX}{{\mathcal X}}
\newcommand{\sY}{{\mathcal Y}}
\newcommand{\sZ}{{\mathcal Z}}
% Sonderbuchstaben mit Doppellinie
\newcommand{\A}{{\mathbb A}}
\newcommand{\B}{{\mathbb B}}
\newcommand{\C}{{\mathbb C}}
\newcommand{\D}{{\mathbb D}}
\newcommand{\E}{{\mathbb E}}
\newcommand{\F}{{\mathbb F}}
\newcommand{\G}{{\mathbb G}}
\renewcommand{\H}{{\mathbb H}}
\newcommand{\I}{{\mathbb I}}
\newcommand{\J}{{\mathbb J}}
\newcommand{\M}{{\mathbb M}}
\newcommand{\N}{{\mathbb N}}
\renewcommand{\P}{{\mathbb P}}
\newcommand{\Q}{{\mathbb Q}}
\newcommand{\R}{{\mathbb R}}
\newcommand{\T}{{\mathbb T}}
\newcommand{\U}{{\mathbb U}}
\newcommand{\V}{{\mathbb V}}
\newcommand{\W}{{\mathbb W}}
\newcommand{\X}{{\mathbb X}}
\newcommand{\Y}{{\mathbb Y}}
\newcommand{\Z}{{\mathbb Z}}
% Misc
\newcommand{\pic}{{\text{Pic}(C,\sD)[E,\nabla]}}
\newcommand{\ocd}{{\Omega^1_C\{\sD\}}}
\newcommand{\oc}{{\Omega^1_C}}
\newcommand{\al}{{\alpha}}
\newcommand{\be}{{\beta}}
\newcommand{\ta}{{\theta}}
\newcommand{\ve}{{\varepsilon}}
\newcommand{\phe}{{\varphi}}
\newcommand{\om}{{\overline M}}
\newcommand{\sym}{{\text{Sym}(\om)}}
\newcommand{\an}{{\text{an}}}
\newcommand{\bs}{{\backslash}}
\newcommand{\lra}{\longrightarrow}
\newcommand{\Sig}{{\Sigma}}
\newcommand{\sig}{{\sigma}}
\newcommand{\dR}{{{\rm dR}}}
\newcommand{\fg}{{\frak {g}}}
\newcommand{\st}{{{\rm st}}}
\newcommand{\pst}{{{\rm pst}}}
\newcommand{\crys}{{{\rm crys}}}
\newcommand{\pri}{{{\rm prim}}}
\newcommand{\Ad}{{{\rm Ad}}}
\newcommand{\Rep}{{{\rm Rep}}}
\newcommand{\Mh}{{{\rm Mor}_{\rm hor}}}
\newcommand{\hor}{{{\rm hor}}}
\newcommand{\La}{{\Lambda}}
\newcommand{\la}{{\lambda}}

%\vspace*{1.2in}

\begin{center} Dedicated to Professor Alexander Beilinson on his 60th birthday. 
\end{center}

\begin{abstract} We define various height functions for motives over number fields. We compare these height functions with classical height functions on  algebraic varieties, and also with analogous height functions for variations of Hodge structures on curves over $\C$. These comparisons provide 
new questions on motives over number fields.

\end{abstract}

\medskip

{\bf Contents}

\medskip
\S0. Introduction

\S1. Height functions for motives

\S2. Speculations

\setcounter{section}{-1}
\section{Introduction}

\begin{para}\label{0.1}

In this paper, for a number field $F$, we define various height functions 

\medskip

(1) $H_{\star}, H_{\star\diamondsuit}, H_{\spadesuit}, H_{\heartsuit,S}, H_{\diamondsuit}, H_{\clubsuit}, \dots : \quad X(F)\to \R_{>0}$,

\medskip
\noindent
where $X(F)=\{\text{motives over $F$}\}$.
See \S\ref{s:2.2} for this set $X(F)$ of motives, which is defined by fixing the Hodge numbers of the motives. 

 We compare these height functions with height functions

\medskip

(2) $V(F)\to \R_{>0}$

\medskip
\noindent
for an algebraic variety $V$ over $F$. (1) and (2) are similar, but a big difference is that usually, $X(F)$ is not understood as the set of $F$-points of an algebraic variety except the case of a Shimura variety which is regarded as a moduli space of motives. 

The height function $H_{\star}$ is almost the same as the height of pure motives in \cite{Ka2}  and in Koshikawa \cite{Ko1} which generalizes the height of  abelian varieties defined by Faltings (\cite{Fa0}). The height function $H_{\star\diamondsuit}$ is almost the same as the height of mixed motives defined in  \cite{Ka3} and is closely related to 
 the height parings of Beilinson \cite{Be2}, Bloch \cite{Bl} and Gillet-Soul\'e (\cite{GS1}, \cite{GS2}), and to  Beilinson regulators in \cite{Be}.

In the above (2), assume that $V$ is a dense open set of  a projective smooth algebraic variety $\bar V$ and that the complement $D=\bar V\smallsetminus V$ is a divisor on $\bar V$ with normal crossings. We will see that the height function $H_{\spadesuit}$ is  similar to the height function $H_{K+D}$ in (2) where $K$ is a canonical divisor of $\bar V$. On the other hand, the  height functions $H_{\star\diamondsuit}$ in (1) is similar to the height function $H_A$ in (2) associated to an ample divisor $A$ on $\bar V$.

\end{para}

\begin{para}\label{0.2} We compare the height functions in (1) and (2) also with height functions 

\medskip

(3) $h_{\star}, h_{\star\diamondsuit}, h_{\spadesuit}, h_{\heartsuit}, h_{\diamondsuit}, h_{\clubsuit},  \dots:\quad \Mh(C, \bar X(\C))\to \R$

\medskip
\noindent
where $C$ is a projective smooth curve over $\C$ and 

\medskip

$\Mh(C, \bar X(\C))=\{\text{variations of Hodge structure on $C$ with log degeneration}\}.$

\medskip

This notation $\Mh(C,\bar X(\C))$ is explained in \ref{0.3} below. The height function which we denote by $h_{\star}$ above was considered in an old paper \cite{Gri} of Griffiths. 
\end{para}

\begin{para}\label{0.3}
Though we do not have the theory of  moduli spaces of motives, we have the philosophy that (like the complex analytic space $V(\C)$ associated to $V$ in (2))

\medskip

(4) the complex analytic space associated to the moduli space of motives  exists and it is the known period domain $X(\C)$ 
which classifies Hodge structures.

\medskip
 See \S\ref{s:2.2} for this space  $X(\C)$ which is defined by fixing the same Hodge numbers as $X(F)$.

By our previous works \cite{KU} and \cite{KNU} Part III with S. Usui and C. Nakayama, we have a
special toroidal partial compactification $\bar X(\C)$ of the period domain $X(\C)$. (See \S\ref{s:2.2}.)
This space $\bar X(\C)$ may be regarded as the analogue for (1) of $\bar V(\C)$ in (2). 

A variation of Hodge structure on $C$ with log degeneration at a finite subset $R$ of $C$ and with the given Hodge numbers determines (and is determined by) a horizontal morphism $C\smallsetminus R \to X(\C)$ (the period map) which extends to a morphism $C\to \bar X(\C)$. $\Mh(C, \bar X(\C))$ denotes the set of morphisms $C\to \bar X(\C)$ such that for some finite subset $R$ of $C$, the image of $C\smallsetminus R$ is contained in $X(\C)$ and the morphism $C\smallsetminus R\to X(\C)$ is horizontal. 

As is explained in \S\ref{s:2.2b}, the comparison of $$V(F)\subset V(\C) \subset \bar V(\C)\quad \text{and}\quad X(F)\to X(\C)\subset \bar X(\C)$$ suggests us that 
$H_{\spadesuit}$ in (1) and $h_{\spadesuit}$ in (3) are similar to  $H_{K+D}$ in (2) and $h_{\heartsuit}$ in (3) is similar to $H_D$ in (2). 

We expect that what happens in the world of motives is suggested by what happens on period domains.
\end{para}

\begin{para}

This (4) 
 and the comparisons of the above (1), (2), (3) provide 
many new questions on motives over number fields which we discuss in Section 2. We will present Questions 1--9. 
For example, we ask

\medskip

(a) how to formulate the motive versions for (1) of  Vojta conjectures \cite{Vo} for (2) (Questions 4,5,6),

(b)  how to formulate the motive versions for (1)  of the theories of the number of points of $V(F)$ of bounded height for (2) (Questions 7,8,9). 

\medskip

Vojta conjectures are stated by using the height functions  $H_{K+D}$ and $H_A$ in (2), and the conjectures of Batyrev and Manin  \cite{BM} on the number of points of bounded height are stated by using relations between the divisors $A$ and $K$.
We will think what are analogues of these for motives in \S\ref{2.4b} and \S\ref{2.5b}, respectively. 

In \S\ref{s:2.1}, we relate the conjectural  Mordell-Weil theorems for the $K$-groups of algebraic varieties over number fields to the finiteness of the number of mixed motives of bounded height with fixed pure graded quotients. 
In \S\ref{2.3b}, basing on the philosophy (4) and the analogy between (1) and (3), we speculate that the analytic theory of curvature forms on Griffiths period domains may be strongly related to the theory of  motives over a number field (Questions 1, 2, 3). 

Though the author can not present enough evidences, he expects that the answers to all these Questions are Yes. He hopes that the new subjects and new problems  in this paper can stimulate new studies in arithmetic geometry.

 \end{para}

 \begin{para}

 We dedicate this paper to Professor Alexander Beilinson on his 60th birthday, who played leading roles in the arithmetic studies of motives. \S2 of this paper contains dreams, may be, too much. But we think he likes dreams in mathematics and hope that he thinks such paper is fine.

 Concerning various subjects in this paper, the author learned many things  from the papers \cite{Ko1} and \cite{Ko2} of Teruhisa Koshikawa and from discussions with 
 him, 
  and  is very thankful to him.

  The author thanks Sampei Usui and Chikara Nakayama that the joint works \cite{KU} and \cite{KNU} with them on extended period domains are useful in this work. 
 The author thanks Spencer Bloch for stimulating discussions on heights of motives.

\end{para}

\section{Height functions for motives}

In this \S1, we define the height functions $H_r(M)$ ($r\in \Z$),  $H_{\clubsuit}(M)$ and $H_{\spadesuit}(M)$  for motives $M$ over a number field. 

These are different from the height of a motive considered in \cite{Ka2}, \cite{Ka3}, \cite{Ko1}, \cite{Ko2}. In this \S1, we slightly change the definition of this height. We denote this modified height by $H_{\star}(M)$ in the case $M$ is pure, and by $H_{\star\diamondsuit}(M)$ in the case $M$ is mixed.

The definitions of $H_r(M)$, $H_{\clubsuit}(M)$, $H_{\spadesuit}(M)$ and $H_{\star}(M)$  are given in \S1.4 after the preliminaries in \S1.1--\S1.3. The definition of $H_{\star\diamondsuit}(M)$ is given in \S1.7 after preliminaries in \S1.5. In \S1.6,  we consider the Hodge theoretic analogues of these height functions.

Heights of motives were studied also by  A. Venkatesh independently (unpublished).  
 
\subsection{Motives}
Let $F$ be a number field, that is, a finite extension of the rational number field $\Q$. For a place $v$ of $F$, let $F_v$ be the local field of $F$. If $v$ is non-archimedean,   let $\F_v$ be the residue field of $v$.

\begin{sbpara}\label{motQ}

We consider mixed motive with $\Q$-coefficients over $F$. How to define it may not be fixed yet in the mathematics of today, but we think that the definition of Jannsen \cite{Ja} is the best for the present paper and we use it. 

Let $M$ be a mixed motive with $\Q$-coefficients over $F$. Then we have the \'etale realization $M_{et}$ of $M$ and the de Rham realization $M_{\dR}$ of $M$. For each embedding $a: F\to \C$, we have the Betti realization $M_{a, B}$ of $M$. These have the following properties.

Let $\hat \Z= \varprojlim_n \Z/n\Z$ and let ${\bf A}^f_\Q=\hat \Z \otimes_\Z \Q$  be the non-archimedean component of the adele ring of $\Q$. Fix an algebraic closure $\bar F$ of $F$. 
Then $M_{et}$ is a free 
 ${\bf A}^f_\Q$-module of finite rank with a continuous action of $\Gal(\bar F/F)$. For each prime number $p$, we denote by $M_{et,\Q_p}$ the $\Q_p$-component of $M_{et}$.

 $M_{\dR}$ is a finite-dimensional $F$-vector space with a decreasing filtration (Hodge filtration). 

$M_{a, B}$ is a finite-dimensional $\Q$-vector space.

These realizations are related to each other (see \ref{comp} below). In particular, $\dim_{\Q_p}(M_{et,\Q_p})=\dim_F(M_{dR})=\dim_\Q(M_{a,B})$, and this dimension is called the rank of $M$.

$M$ has an increasing filtration $W$, called the weight filtration, by sub mixed motives $W_w=W_wM$ ($w\in \Z$) with $\Q$-coefficients over $F$. This weight filtration induces weight filtrations on the above realizations. 
A mixed motive $M$ with $\Q$-coefficients over $F$  is called a pure motive of weight $w$ if the weight filtration $W$ satisfies $W_w=M$ and $W_{w-1}=0$. We have a pure motive $\gr^W_wM=W_wM/W_{w-1}M$ of weight $w$ for each $w\in \Z$. 
\end{sbpara}

\begin{sbpara}\label{comp} We give only rough descriptions of the relations between realizations. See \cite{Ja} for details.

(1) We have comparison isomorphisms  $\C\otimes_{\Q} M_{a,B} \cong \C \otimes_F M_{\dR}$  and ${\bf A}^f_\Q \otimes_{\Q} M_{a,B} \cong M_{et}$, which are compatible with the weight filtrations, and the former isomorphism produces a mixed Hodge structure $M_{a,H}$. 

(2) Let $v$ be a non-archimedean place of $F$ and let $p=\text{char}(\F_v)$ (the characteristic of $\F_v$). Then $M_{et,\Q_p}$ is a de Rham representation of $\Gal(\bar F_v/F_v)$, $D_{\dR}(M_{et, \Q_p})=D_{dR}(F_v, M_{et, \Q_p})$ is identified with $F_v \otimes_F M_{\dR}$, and the filtration $(D_{dR}^r(M_{p,et}))_r$ and the weight filtration of the former are  identified with the Hodge filtration and the weight filtration of the latter, respectively.

\end{sbpara}

\begin{sbpara} Roughly speaking, a mixed motive with $\Q$-coefficients over $F$ in the sense of \cite{Ja} is a system of realizations which comes from geometry. Morphisms are homomorphisms of systems of realizations. 

The category of mixed motives with $\Q$-coefficients over $F$ is an abelian category. It is a Tannakian category, and has direct sums, tensor products, symmetric powers $\text{Sym}^r$, exterior powers $\wedge^r$ ($r\in \Z$),  the dual $M \mapsto M^*$, and Tate twists $M\mapsto M(r)$ ($r\in \Z$). 

\end{sbpara}

\begin{sbpara}

For a pure motive $M$ with $\Q$-coefficients over $F$ of weight $w$, a polarization of $M$ is a morphism $M\otimes M \to \Q(-w)$ which gives a polarization of the Hodge structure $M_{a,H}$  for any $a:F\to \C$. 

\end{sbpara}

\begin{sbpara}\label{semi-st}
 
Let $M$ be a mixed motive with $\Q$-coefficients over a number field $F$. 

For a non-archimedean place $v$ of $F$, we say $M$ is of good reduction (resp. semi-stable reduction) at $v$ if it satisfies the following two conditions (i) and (ii). Let $p=\text{char}(\F_v)$.

\medskip

(i) The representation of $\Gal(\bar F_v/F_v)$ on $M_{et,\Q_p}$ is crystalline (resp. semi-stable).

(ii) For any prime number $\ell\neq p$, the action of the inertia subgroup of $\Gal(\bar F_v/F_v)$ on $M_{et,\Q_\ell}$ is trivial (resp. unipotent). 

\medskip
We say $M$ is of bad reduction at $v$ if it is not of good reduction at $v$.

We say $M$ is of semi-stable reduction if it is of semi-stable reduction at any non-archimedean place of $F$. 
\end{sbpara}

\begin{sbpara}\label{motZ}

By a mixed motive with $\Z$-coefficient over $F$, we mean a 
 mixed motive $M$  with $\Q$-coefficients over $F$  endowed with a $\Gal(\bar F/F)$-stable $\hat \Z$-lattice $T$ in  $M_{et}$. 
 
 For a mixed motive $M$ with $\Z$-coefficients over $F$, we denote this $T$ by $M_{et,\hat \Z}$, and denote its $p$-adic part for a prime number $p$ by $M_{et,\Z_p}$. 
 
 If $M$ is a mixed motive with $\Z$-coefficients over $F$, we have a $\Z$-lattice of the Betti-realization  $M_{a,B}$ for $a:F\to \C$ as the inverse image of $M_{et, \hat \Z}$ in $M_{a,B}$ via the comparison isomorphism in \ref{comp}. We will denote 
 this $\Z$-lattice as $M_{a,B, \Z}$. 
 
\end{sbpara}

\begin{sbpara} When we consider a polarization of a pure motive with $\Z$-coefficients, we do not assume that it is integral (that is, we do not assume that the image of the induced pairing $M_{et,\hat \Z} \otimes M_{et, \hat \Z} \to {\bf A}^f_\Q(-w)$ is contained in $\hat \Z(-w)$). This convention should be good when we consider the duality, for the integrality of the polarization is not preserved in the duality. 
\end{sbpara}

\subsection{Ideas of the definitions of the height functions for motives}

\begin{sbpara}\label{Cases12}
In this \S1, we will consider three cases.

\medskip

Case (pure): $M$ is a pure motive with $\Z$-coefficients over a number field $F$ of weight $w$.

Case (mixed): $M$ is a mixed motive  with $\Z$-coefficients over a number field $F$.

Case (mixed-pol): $M$ is a mixed motive  with $\Z$-coefficients over a number field $F$ whose pure graded quotients $\gr^W_wM$ ($w\in \Z$) are polarized. 

\medskip

In this \S1, we will define $H_r(M)\in \R_{>0}$ ($r\in \Z$), $H_{\star}(M)$, $H_{\clubsuit}(M) \in \R_{>0}$ in Case (pure),    $H_{\spadesuit}(M)\in \R_{>0}$ in Case (mixed),   and $H_{\diamondsuit}(M)$ and $H_{\star\diamondsuit}(M)$ in Case (mixed-pol). 
The height function $H_{\star}(M)$ (resp. $H_{\star\diamondsuit}(M)$) is a modified version of the height defined in 
 \cite{Ka2}, \cite{Ko1}, \cite{Ko2} (resp. \cite{Ka3}). The height function  $H_{\diamondsuit}(M)$ is the ratio of $H_{\star\diamondsuit}(M)$ and $\prod_w H_{\star}(\gr^W_wM)$.

The height function $H_{\spadesuit}(M)$ has the importance that it is similar to the height function $H_{K+D}$ in the comparison of (1) and (2) in \ref{0.1}, and is related to differential forms and a canonical measure on a period domain. $H_{\clubsuit}(M)$ is a simpler version of $H_{\spadesuit}(M)$ and is also related to differential forms and a canonical measure on some period domain. 
$H_{\star}(M)$ and $H_{\star\diamondsuit}(M)$ have  the importance that in the comparison with (2) in \ref{0.1}, it is (conjecturally) like the height on $V(F)$ associated to an ample line bundle on $\bar V$. $H_r(M)$ for $r\in \Z$ have the importance that these $H_{\clubsuit}(M)$ and $H_{\spadesuit}(M)$, and   $H_{\star}(M)$  are expressed by using $H_r(M)$. 

\end{sbpara}

\begin{sbpara}\label{idea1}  We describe the ideas of the definitions of 

\medskip

(i) $H_r(M)$ ($r\in \Z$), (ii) $H_{\star}(M)$, (iii)  $H_{\clubsuit}(M)$, (iv) $H_{\spadesuit}(M)$.
\medskip

Here in (i), (ii), (iii), we assume $M$ is in Case (pure). In (iv), we assume $M$ is in Case (mixed-pol). For each (i)--(iv), 
consider the one-dimensional $F$-vector  space $L$ defined as

\medskip

(i)  $L= (\text{det}_F\gr^r M_{\dR})\otimes_F (\text{det}_F\gr^{w-r}M_{\dR})^{-1}$, 

(ii) $L= \otimes_{r\in \Z}  ((\text{det}_F\gr^r M_{\dR})\otimes_F (\text{det}_F\gr^{w-r}M_{\dR})^{-1})^{\otimes r}$.

(iii) $L=\text{det}_F(\text{End}_F(M_{dR})/F^0)^{-1}$,

(iv) $L= \text{det}_F(\text{End}_{F,W, \langle\; \rangle}(M_{dR})/F^0)^{-1}$,

\medskip

Here in (iv), 
$\text{End}_{F,W, \langle \; \rangle}(M_{dR})$ denotes the $F$-vector space of all $F$-linear maps $h: M_{dR}\to M_{dR}$ which respect the weight filtration $W$ such that $\langle h(x), y\rangle+ \langle x, h(y)\rangle_w=0$ for any $w\in \Z$ and any $x,y \in \gr^W_wM_{\dR}$ where  $\langle\;,\rangle: \gr^W_wM_{\dR}\times \gr^W_wM_{\dR} \to F$ is the pairing defined by the polarization of $\gr^W_wM$.  Both in (iii) and (iv), $F^0$ denotes the $0$-th Hodge filtration.

Here and in the following, the notation $F$ for a number field and that for Hodge filtration may be confusing, but the author hopes that the reader can distinguish them. To avoid the confusion, we will denote $\gr^r_F$ for the Hodge filtration $F$ just as $\gr^r$.

As is explained in \S1.4, for each place $v$ of $F$, we have a canonical metric $|\;\;|_v$ on $L$. If $v$ is archimedean, this comes from Hodge theory. If $v$ is non-archimedean, this comes from $p$-adic Hodge theory. 
In the definition of $|\;\;|_v$ for $v$ non-archimedean, we use the method of Koshikawa \cite{Ko1}, \cite{Ko2} which he used for the height of a pure motive (and which was different from the original method in \cite{Ka2} by the author).  For a non-zero element $x$ of $L$, $|x|_v=1$ for almost all places $v$ of $F$.

We define 

\medskip

$H_r(M):= \prod_v |x|_v^{-1/2}$ in (i), \quad $H_{\star}(M):= \prod_v |x|_v^{-1/2}$ in  (ii),

$H_{\clubsuit}(M):= \prod_v |x|_v^{-1}$ in (iii), \quad $H_{\spadesuit}(M):= \prod_v |x|_v^{-1}$ in (iv).

\medskip

By the product formula, this does not depend on the choice of $x\in L\smallsetminus \{0\}$. 

In (iv), $H_{\spadesuit}(M)$ is independent of the choice of the polarizations of $\gr^W_wM$ ($w\in \Z$), and hence $H_{\spadesuit}(M)$ is defined for $M$ in Case (mixed). 

\end{sbpara}

\begin{sbpara}\label{idea2}  We describe the idea of the definition of $H_{\star\diamondsuit}(M)$ with $M$ in Case (mixed-pol). It is defined in the form 
$$H_{\star\diamondsuit}(M):= (\prod_{w\in \Z} H_{\star}(\gr^W_wM)) \cdot H_{\diamondsuit}(M)$$
where $H_{\diamondsuit}(M) \in \R_{>0}$ describes, roughly speaking, how $M$ is far from being the direct sum of $\gr^W_wM$ ($w\in \Z$). This $H_{\diamondsuit}(M)$ is defined in \cite{Ka3} and is reviewed in \S1.7.

\end{sbpara}

\subsection{Review on $p$-adic Hodge theory}

 We review some parts of the work of Bhatt, Morrow, Scholze  \cite{BMS} on integral $p$-adic Hodge theory.

\begin{sbpara} Let $K$ be a complete discrete valuation field of mixed characteristic $(0,p)$ with perfect residue field $k$.  
Let $\bar K$ be an algebraic closure of $K$ and let $C$ be the completion of $\bar K$. Let $\bar k$ be the residue field of $\bar K$, which is an algebraic closure of $k$. Let $O_K$, $O_{\bar K}$ and $O_C$ be the valuation rings of $K$, $\bar K$ and $C$, respectively. 

The main subject we review in this \S1.3 from \cite{BMS} is that, for a de Rham representation $V$ of $\Gal(\bar K/K)$, a $\Z_p$-lattice $T\subset V$ gives an $O_C$-lattice $D_{\dR,O_C}(T)$ in  $C \otimes_K D_{\dR}(V)$.

\end{sbpara}

\begin{sbpara}\label{Forings}
We review the rings $W(O_{C^{\flat}})$, $B_{\dR}^+$ and $B_{\dR}$ of Fontaine.  

Let $O_{C^{\flat}}:= \varprojlim_{n\geq 0} O_{\bar K}/pO_{\bar K}=\varprojlim_{n\geq 0} O_C/pO_C$ where the inverse limits are taken with respect to the homomorphisms $x\mapsto x^p$.  Then $O_{C^{\flat}}$ is a perfect ring and it is an integral domain. The field of fractions of $O_{C^{\flat}}$ is denoted by $C^{\flat}$. (Actually,  $O_{C^{\flat}}$ is a complete valuation ring of rank $1$, and $C^{\flat}$ is an algebraically closed field.)

Consider the ring $W(O_{C^{\flat}})$ of Witt vectors with entries in $O_{C^{\flat}}$. The $p$-th power map of $O_{C^{\flat}}$ induces a ring automorphism $\varphi: W(O_{C^{\flat}}) \to W(O_{C^{\flat}})$  called Frobenius. 

We have a unique ring homomorphism $\theta: W(O_{C^{\flat}})\to O_C$ which lifts the composition $W(O_{C^{\flat}})\to O_{C^{\flat}}\to O_C/pO_C$ (the last arrow is $(a_n)_{n\geq 0}\mapsto a_0$). This homomorphism $\theta$ is surjective and its kernel $\frak p$ is a principal prime ideal. 
The local ring of $W(O_{C^{\flat}})$ at $\frak p$ is a discrete valuation ring whose completion is denoted by $B_{\dR}^+$. $B_{\dR}$ is defined to be the field of fractions of $B_{\dR}^+$. 

\end{sbpara}

\begin{sbpara}\label{ThFar} Let $\sC_1$ be the category of free $W(O_{C^{\flat}})$-modules $M$ of finite rank endowed with a Frobenius linear bijection $\varphi_M: M[1/\xi] \to M[1/\varphi(\xi)]$ where $\xi$ is a generator of the ideal $\frak p=\text{Ker}(\theta: W(O_{C^{\flat}})\to O_C)$. (Such $M$ with $\varphi_M$ is called a (finite free) Breuil-Kisin-Fargues module.) On the other hand, 
let $\sC_2$ be the category of  pairs $(T, \Xi)$, where $T$ is a free $\Z_p$-module of finite rank and $\Xi$ is a $B_{\dR}^+$-lattice in $B_{\dR}\otimes_{\Z_p} T$. Then we have an equivalence of categories $$\sC_1\simeq \sC_2.$$
This is a theorem of Fargues. A proof is given in \cite{SP}. 

When an object $M$ of $\sC_1$ and an object $(T, \Xi)$ of $\sC_2$ correspond, then
$$\Xi= B_{\dR}^+\otimes_{W(O_{C^{\flat}})} M.$$
\end{sbpara}

\begin{sbpara}\label{dRint} Let $V$ be a de Rham representation of $\Gal(\bar K/K)$ and let $T$ be a $\Z_p$-lattice in $V$. Then we have an $O_C$-lattice in $C\otimes_K D_{\dR}(V)$, which we denote by $D_{\dR,O_C}(T)$, defined as follows. 
We have
$B_{\dR}\otimes_{\Z_p} T= B_{\dR}\otimes_{\Q_p} V= B_{\dR}\otimes_K D_{\dR}(V)$. 
Take $\Xi:= B_{\dR}^+\otimes_K D_{\dR}(V)$. Let $M$ be the object of $\sC_1$ corresponding to the object $(T, \Xi)$ of $\sC_2$. 
Then we have $B_{\dR}^+\otimes_{W(O_{C^{\flat}})}  M= B_{\dR}^+\otimes_K D_{\dR}(V)$. By using the canonical homomorphism $B_{\dR}^+\to C$, we have $C\otimes_{W(O_{C^{\flat}})} M= C \otimes_K D_{\dR}(V)$. 
Let $$D_{\dR, O_C}(T):= O_C\otimes_{W(O_{C^{\flat}})} M \subset C \otimes_K D_{\dR}(V).$$ 

This construction of the $O_C$-lattice $D_{\dR,O_C}(T)$ in $C\otimes_K D_{\dR}(V)$ is compatible with $\oplus$, $\otimes$, the dual, $\text{Sym}^r$, $\wedge^r$ ($r\geq 0$), and Tate twists. 

Here the compatibility with Tate twists mean the following. For a de Rham representation $V$ and for  $r\in \Z$, $D_{\dR}(V)$ and $D_{\dR}(V(r))$ are identified as $K$-vector spaces. We have $D_{\dR,O_C}(T)= D_{\dR,O_C}(T(r))$ in $C \otimes_K D_{\dR}(V)= C\otimes_K D_{\dR}(V(r))$.

\end{sbpara}

\begin{sbpara}\label{dRint3} For $V$ and $T$ as in \ref{dRint}, we have also an $O_C$-lattice $D^r_{\dR,O_C}(T)= D_{\dR,O_C}(V)\cap (C \otimes_K D^r_{\dR}(V))$ of $C \otimes_K D^r_{\dR}(V)$ and an $O_C$-lattice $\gr^rD_{\dR,O_C}(T)= D^r_{\dR,O_C}(T)/D^{r+1}_{\dR,O_C}(T)$ of $C\otimes_K \gr^r D_{\dR}(V)$ for $r\in \Z$. 

\end{sbpara}

\begin{sbpara}\label{geom} Geometric meanings of $D_{\dR, O_C}(T)$.

Let  $X$ be a  proper smooth scheme over $K$, let $T=H^m_{et}(X_{\bar K}, \Z_p)$, $V= H^m_{et}(X_{\bar K}, \Q_p)$, and identify $D_{dR}(V)$ with $H^m_{dR}(X/K)$ and $\gr^rD_{\dR}(V)$ with $H^{m-r}(X, \Omega^r_{X/K})$. 

(1) Assume $X=\frak X\otimes_{O_K} K$ for a 
 proper smooth scheme $\frak X$ over $O_K$. Assume that $H^i(\frak X, \Omega^j_{\frak X/O_K})$ are torsion free for  all $i,j$. Then $T$ is torsion free, $D_{\dR,O_C}(T)=O_C\otimes_{O_K} H^m_{\dR}(\frak X/O_K)$, and $\gr^rD_{\dR,O_C}(T)=O_C\otimes_{O_K} H^{m-r}(\frak X, \Omega^r_{\frak X/O_K})$. This follows from \cite{BMS}. 

(2) Assume $X= \frak X\otimes_{O_K} K$ for a proper regular flat scheme $\frak X$ over $O_K$ of semi-stable reduction. Assume that $H^i(\frak X, \Omega^j_{\frak X/O_K}(\log))$ are torsion free for all $i,j$. Then $T$ is torsion free, $D_{\dR,O_C}(T)=O_C\otimes_{O_K} H^m_{\log \dR}(\frak X/O_K)$, and $\gr^rD_{\dR,O_C}(T)=O_C\otimes_{O_K} H^{m-r}(\frak X, \Omega^r_{\frak X/O_K}(\log))$. 
This follows from the works of Faltings \cite{Fa1}, \cite{Fa2} under a mild assumption (see \S3.3 of \cite{Ko1}), and is proved in general by $\check{\rm C}$esnavi$\check{\rm c}$ius \cite{Ce}.

\end{sbpara}

\begin{sbpara}\label{st}  Let $K_0$ be the field of fractions of $W(k)$, let $K_{0,\text{ur}}\subset \bar K$ be the maximal unramified extension of $K_0$, and let $\hat K_{0,\text{ur}}=\text{frac}(W(\bar k))\subset C$ be the completion of $K_{0,\text{ur}}$.

Let $V$ be a de Rham representation of $\Gal(\bar K/K)$. Then $V$ is potentially  semi-stable and we have a $K_{0,\text{ur}}$-vector space $D_{\text{pst}}(V)$ of dimension $\dim_{\Q_p}(V)$.
Let $T$ be a $\Z_p$-lattice in $V$,  let $\Xi= B_{\dR}^+\otimes_K 
D_{\dR}(V) $, and $M$ be the object of $\sC_1$ corresponding to the object $(T, \Xi)$ of $\sC_2$. 
Let
$$D_{\pst,W(\bar k)}(T):= W(\bar k) \otimes_{W(O_{C^{\flat}})} M$$
where the homomorphism $W(O_{C^{\flat}})\to W(\bar k)$ is induced by $O_{C^{\flat}} \to O_C/pO_C\to \bar k$.

We have a canonical isomorphism
$$ {\hat K}_{0,\text{ur}}
 \otimes_{K_{0,\text{ur}}}D_{pst}(V)\cong {\hat K}_{0,\text{ur}} \otimes_{W(\bar k)} D_{\pst,W(\bar k)}(T)$$
of vector spaces over $\hat K_{0,\text{ur}}$. 
This isomorphism follows from the theory of Kisin on Beuil-Kisin modules (\cite{Ki}) and the relation of Breuil--Kisin modules and Breuil--Kisin--Fargues modules explained in \cite{BMS}. 
Hence $D_{\pst,W(\bar k)}(T)$ is regarded as a $W(\bar k)$-lattice in $\hat K_{0,\text{ur}}\otimes_{K_{0,\text{ur}}} D_{\pst}(V)$.

This construction of the $W(\bar k)$-lattice $D_{\pst,W(\bar k)}(T)$ in $\hat K_{0,\text{ur}} \otimes_{K_{0,\text{ur}}} D_{\pst}(V)$  is compatible with $\oplus$, $\otimes$, the dual, $\text{Sym}^r$, $\wedge^r$ ($r\geq 0$), and Tate twists. 
Here the compatibility with Tate twists mean the following. For a de Rham representation $V$ and for  $r\in \Z$, $D_{\pst}(V)$ and $D_{\pst}(V(r))$ are identified as $K_{0,\text{ur}}$-vector spaces. We have $D_{\pst,W(\bar k)}(T)= D_{\pst, W(\bar k)}(T(r))$ in $\hat K_{0,\text{ur}} \otimes_{K_{0,\text{ur}}} D_{\pst}(V)= \hat K_{0,\text{ur}} \otimes_{K_{\text{ur}}} D_{\pst}(V(r))$.

\end{sbpara}

\begin{sbpara}\label{st2}
Geometric meanings of $D_{\pst,W(\bar k)}(T)$. 

Let  $X$ be a  proper smooth scheme over $K$, let $T=H^m_{et}(X_{\bar K}, \Z_p)$ and let $V= H^m_{et}(X_{\bar K}, \Q_p)$.

\medskip

(1) Assume $X=\frak X\otimes_{O_K} K$ for a 
 proper smooth scheme $\frak X$ over $O_K$ and let $Y= X \otimes_{O_K} k$. Assume that $H^i_{\crys}(Y/W(k))$ are torsion free for all $i$.  Then $T$ is torsion free. Via the classical isomorphism $D_{\pst}(V)\cong K_{0,\text{ur}}\otimes_{W(k)} H^m_{\crys}(Y/W(k))$, we have  $D_{\pst,W(\bar k)}(T)=W(\bar k)\otimes_{W(k)} H^m_{\crys}(Y/W(k))$.  This follows from \cite{BMS}. 

(2) The generalization of (1) to the semi-stable reduction case (using log crystalline cohomology) follows from 
the works of Faltings \cite{Fa1}, \cite{Fa2} under a mild assumption (see \S3.3 of \cite{Ko1}). The author learned that the proof  in general was given by a recent joint work of B. Bhatt, M. Morrow, P. Scholze, and also by a recent joint work of K. $\check{\rm C}$esnavi$\check{\rm c}$ius
and T. Koshikawa. 
 \medskip
 
 We will not use this \ref{st2} in this paper, but we will use 
the integral structure of $\hat K_{0,\text{ur}}\otimes_{K_{0,\text{ur}}} D_{\pst}(V)$ in \ref{st} in the definition of the height function $H_{\heartsuit,S}$ for motives in \S2.3. 

\end{sbpara}

\subsection{Height functions $H_r(M)$, $H_{\clubsuit}$, $H_{\spadesuit}$, and $H_{\star}$}

We define the height functions for motives, $H_r$ ($r\in\Z$), $H_{\star}$, $H_{\clubsuit}$ and $H_{\spadesuit}$. 

We first define $H_r$ ($r\in \Z$) in \ref{defht1} after preparations in \ref{Hr1}, \ref{Hr2}, and then define $H_{\star}$, $H_{\clubsuit}$ and $H_{\spadesuit}$ in \ref{Hscs}.

\begin{sbpara}\label{Hr1} 
We give a preparation on Hodge theory. 

Let $H=(H_\R, F)$ be a pure $\R$-Hodge structure of weight $w$ ($H_\R$ is the $\R$-structure and $F$ is the Hodge filtration on $H_\C=\C\otimes_\R H_\R$). Consider the one-dimensional $\C$-vector space  
$I=(\text{det}_\C \gr^r H_\C) \otimes_\C (\text{det}_C \gr^{w-r}H_\C)^{-1}$. We define a metric $|\;\;|$ on $I$ by
$$|a\otimes b^{-1}|:= |a/\bar b|=|\bar a/b|$$
where $a$ is a  basis of $\text{det}_\C \gr^r H_\C$, $b$ is a basis of 
$\text{det}_C \gr^{w-r}H_\C$, $\bar a$ is a basis of $\text{det}_\C\gr^{w-r}H_\C$ obtained from $a$ as below,    $\bar b$ is a basis of $\text{det}_\C\gr^rH_\C$ obtained from $b$ as below, and we take the ratios $a/\bar b, \bar a/b\in \C^\times$ of two bases of one-dimensional $\C$-vector spaces and take their absolute values $|a/\bar b|,\; |\bar a/b|\in \R_{>0}$. 

The map $a\mapsto \bar a$ and $b\mapsto \bar b$ are given by
the semi-linear isomorphism
$\text{det}_\C\gr^{w-r}H_\C\cong \text{det}_\C\gr^rH_\C$ 
(semi-linearity is for the complex conjugation) induced by the composition of the bijections

\medskip

(*) \;  $\gr^{w-r} H_\C \cong H^{w-r,r}\to  H^{r,w-r}\cong \gr^rH_\C$.

\medskip

Here $H^{p,q}_\C$ ($p+q=w$) denotes the component of $H_\C$ of Hodge type $(p,q)$, and the central arrow in (*)  is the semi-linear isomorphism induced by $\C\otimes_\R H_\R\to \C\otimes_\R H_\R\;;\;z\otimes v \mapsto \bar z \otimes v$ ($z\in \C$, $v\in H_\R$).

\end{sbpara}

\begin{sbpara}\label{Hr2} In \ref{Hr1}, if a polarization $p$ on $H$ is given, we have the Hodge metric on $H_\C$ defined by $p$, hence $\gr^rH_\C$ and $\gr^{w-r}H_\C$ have the induced Hermitian metrics, and hence $I=(\text{det}_\C\gr^rH_\C) \otimes_\C (\text{det}_C\gr^{w-r}H_\C)^{-1}$ has the induced  metric. This metric coincides with the metric defined in \ref{Hr1} which was defined not using a polarization. 

This can be seen as follows. Let $a$ be a basis of $\text{det}_\C\gr^rH_\C$ and $b$ be a basis of $\text{det}_\C \gr^{w-r}H_\C$, and let $|a|_p$ and $|b|_p$ be the metic defined by the Hodge metric on $H_\C$ associated to $p$, respectively. Write $a=z\bar b$ in $\text{det}_\C\gr^rH_\C$ with $z\in \C^\times$. Let $\langle\;,\;\rangle:\gr^rH_\C\times \gr^{w-r}H_\C\to \C$ be the non-degenerate bilinear pairing induced by $p$, and denote the induced bilinear pairing  $\text{det}_C\gr^rH_\C\times \text{det}_\C\gr^{w-r}H_\C\to \C$
also by $\langle\;,\;\rangle$. Then we have
$|a|_p=|\langle a, \bar a\rangle|^{1/2}= |z||\langle b, \bar b\rangle|^{1/2}=|z||b|_p$ and hence 
$|a|_p|b|_p^{-1}=|z|$.

\end{sbpara}

\begin{sbpara}\label{defht1}

Let $M$ be as in Case (pure) (\ref{Cases12}). Let $r\in \Z$ and let 
$$L_r(M):= (\text{det}_F\gr^rM_{\dR})\otimes_F (\text{det}_F \gr^{w-r} M_{\dR})^{-1}.$$ 
For each place $v$ of $F$, we define a metic $|\;\;|_v$ on $L_r(M)$
as follows.

 If $v$ is archimedean, take a homomorphism $a:F\to \C$ which induces the place $v$ and let $H$ be the $\R$-Hodge structure $\R\otimes_\Q M_{a,H}$ (\ref{comp} (1)).  We have $\C \otimes_F L_r(M) \cong I$ with $I$ of $H$ as above. For $x\in L_r(M)$, by using the above metric $|\;\;|$ on $I$, we define 
$|x|_v=|x|$ if $v$ is a real place, and $|x|_v= |x|^2$ if $v$ is a complex place. 

Let $v$ be a non-archimedean place and  let $C$ be the completion of an algebraic closure of $F_v$.

In general, for a one-dimensional $C$-vector space $J$  endowed with an $O_C$-lattice $J_{O_C}$, we have the associated  metric $|\;\;|$ on $J$ defined as
$$|x|= \inf\{ |a|\;|\;a\in C,  x \in aJ_{O_C}\}$$
where $a \mapsto |a|$ is the absolute value of $C$ whose restriction to $F_v$ coincides with the standard absolute value of $F_v$.

Let $p=\text{char}(\F_v)$ and consider $D_{\dR}(M_{et,\Q_p})$ of the representation $M_{et,\Q_p}$ of $\Gal(\bar F_v/F_v)$. 
We have the metric 
 on 
 \medskip
 
 $C \otimes_F L_r(M)= C\otimes_{F_v}
(\text{det}_{F_v} \gr^r D_{\dR}(M_{et, \Q_p})) \otimes_{F_v} (\text{det}_{F_v} \gr^{w-r} D_{dR}(M_{et, \Q_p}))^{-1}$  

\medskip
\noindent
 associated to the $O_C$-lattice 
 
 \medskip
 $(\text{det}_{O_C}\gr^rD_{\dR,O_C}(M_{et, \Z_p})) \otimes_{O_C}(\text{det}_{O_C}\gr^{w-r}D_{\dR,O_C}(M_{et, \Z_p}))^{-1}$ 
 
 \medskip
 \noindent
(\ref{dRint3}). This metric induces a metric $|\;\;|_v$ on $L_r(M)$.

For $x\in L_r(M)\smallsetminus \{0\}$, we have $|x|_v=1$ for almost all places $v$ of $F$. This follows from (1) of \ref{geom}. 

We define $$H_r(M):= \prod_v |x|_v^{-1/2}.$$
By the product formula, this is independent of $x$.

\end{sbpara}

\begin{sbprop}\label{dtetc} Let $M$ be  in Case (pure).

(1)  $H_r(M^*)=H_{-r}(M)^{-1}$. Here $M^*$ denotes the dual of $M$. 

(2) $H_r(M(i))= H_{r+i}(M)$ for the $i$-th Tate twist.

(3) $H_r(M)H_{w-r}(M)=1$ and 
$\prod_{r \in \Z} \; H_r(M)=1$.

(4) For a finite extension $F'$ of $F$, we have 
$H_r(M')=H_r(M)^{[F':F]}$ where $M'$ denotes the pure motive with $\Z$-coefficients over $F'$ induced by $M$.

(5)  Let $n\in \Q_{>0}$.
 Let $M'$ be the pure motive with $\Z$-coefficients over $F$ whose underlying pure motive with $\Q$-coefficients is the same as that of $M$ such that $M'_{et,\hat \Z}=n\cdot M_{et,\hat \Z}$.Then $$H_r(M')=H_r(M).$$.

(6) Let $M'$ be  in Case (pure). Then $H_r(M\oplus M')=H_r(M)H_r(M')$.

\end{sbprop}

These are evident.

\begin{sbpara}\label{Hscs} We define $H_{\star}$, $H_{\clubsuit}$, $H_{\spadesuit}$.

For $M$ in Case (pure), let 
$$H_{\star}(M)= \prod_{r\in \Z} \; H_r(M)^r.$$
This $H_{\star}(M)$ is a modified version of the height of $M$ in Case (pure) given in \cite{Ka2} and in Koshikawa  \cite{Ko1}, \cite{Ko2}.  The precise relation to the height by Koshikawa is given in  \ref{met3} below.

For $M$ in Case (pure), let
$$H_{\clubsuit}(M) = \prod_{r<0} H_r(M^*\otimes M)^{-1}$$
where $M^*$ denotes the dual of $M$. 

For $M$ in Case (mixed), let
$$H_{\spadesuit}(M)=$$ $$
 \prod_{w<0, r<0} H_r(\gr^W_w(M^*\otimes M))^{-1}\cdot \prod_{w:\text{odd}, r<w} H_r(\text{Sym}^2 \gr^W_wM)^{-1}\cdot  \prod_{w:\text{even}, r<w} H_r(\wedge^2 \gr^W_wM)^{-1}.$$

\end{sbpara}

$H_{\clubsuit}(M)$ (resp. $H_{\spadesuit}(M)$) is expressed by using $H_r(M)$ (resp. $H_r(\gr^W_wM)$) as in \ref{rsc} below.
 
\begin{sbpara}\label{ab}
(1) For a map $h:\Z\to \Z_{\geq 0}$ with finite support and for $r\in \Z$, let 
$$a(h,r)= \sum_{i<r} h(i)\in \Z_{\geq 0}.$$

 (2) For a map $h:\Z^2\to \Z_{\geq 0}$ with finite support and for $w,r\in \Z$, define $b(h,w,r)\in \Z_{\geq 0}$ as follows. 
Let 
$A(w,r)= \{(w',r')\in \Z^2\;|\; w'<w, r'<r\}$, $B(w,r)=\{(w',r')\in \Z^2\;|\; w'>w, r'>r\}$, $C(w,r)= \{(w',r')\in \Z^2\;|\; w'=w, r'<r, r+r'\neq w\}$. 
Define
$$b(h,w,r):= (\sum_{(w',r')\in A(w,r)\cup C(w,r)} h(w',r'))-(\sum_{(w',r')\in B(w,r)} h(w',r'))+e(w,r)$$
where $e(w,r)=h(w,r)-(-1)^w$ if $r>w/2$ and $e(w,r)=0$ if $r\leq w/2$.
\end{sbpara}
\begin{sbprop}\label{rsc}  (1) Let $M$ be in Case (pure). Then
$$H_{\clubsuit}(M)= \prod_{r\in \Z}\; H_r(M)^{a(h,r)}$$
where $a(h,r)$ is as in \ref{ab} (1) with $h(r)= \dim_F \gr^rM_{\dR}$ (\ref{ab} (1)). 

(2) Let $M$ be in Case (mixed). Then 
$$H_{\spadesuit}(M)= \prod_{w,r\in \Z}\; H_r(\gr^W_wM)^{b(h,w,r)}$$
where $b(h, w,r)$ is as in \ref{ab} (2) with $h(w,r)= \dim_F \gr^r\gr^W_w(M_{\dR})$ (\ref{ab} (2)). 

\end{sbprop}

 This  Proposition \ref{rsc} will be proved  in this \S1.4 later. 
 
 In the rest of this \S1.4, we give comments  on $H_{\star}$ first, then on $H_{\clubsuit}$ and then on $H_{\spadesuit}$, and 
then give a Proposition \ref{Htprop2} on the properties of these height functions.

\begin{sblem}\label{alg22}  Let $k$ be a field, let $w\in \Z$, and let $H$ be a finite dimensional $k$-vector space endowed with a finite decreasing filtration such that $h(r)=h(w-r)$, where $h(r):=\dim_k \gr^r H$, for any $r\in \Z$. For $r\in \Z$, let $L_r= (\text{det}_k \;\gr^r H)\otimes_k (\text{det}_k \;\gr^{w-r} H)^{-1}$. 
Then we have  a canonical isomorphism
$$(\otimes_{r\in \Z} (\text{det}_k\; \gr^r H)^{\otimes r})^{\otimes 2} \cong (\otimes_{r\in \Z} L_r)^{\otimes r} \otimes_k (\text{det}_k H)^{\otimes w}.$$

\end{sblem}

\begin{pf} We have canonical isomorphisms
$$(\otimes_{r\in \Z} (\text{det}_k \gr^r H)^{\otimes r})^{\otimes 2} \cong (\otimes_{r\in \Z} (\text{det}_k \gr^r H)^{\otimes r}) \otimes_k (\otimes_{r\in \Z} (\text{det}_k \gr^{w-r} H)^{\otimes w-r})$$ $$\cong (\otimes_{r\in \Z} L_r^{\otimes r})\otimes_k (\otimes_{r\in \Z} (\text{det}_k \gr^{w-r} H)^{\otimes w}) \cong (\otimes_{r\in \Z} L_r^{\otimes r})\otimes_k (\text{det}_k H)^{\otimes w}.$$ 
\end{pf}

\begin{sbpara}\label{met23} As a preparation for  \ref{met3}, we consider motives of rank one. Let $M$ be as in Case (pure) and assume that $M$ is of rank $1$. For each place $v$ of $F$, we define a metric $|\;\;|_v$ on the one-dimensional $F$-vector space $M_{\dR}$ as follows. 

Assume first that $v$ is archimedean and let $a: F\to \C$ be a homomorphism which induces $v$. For $x\in M_{\dR}$, denote the image of $x$ under the isomorphism $\C\otimes_F M_{\dR} \cong  \C\otimes_\Z M_{a,B,\Z}$
as $z \otimes \gamma$ where $z\in \C$ and $\gamma$ is a $\Z$-basis of  $M_{a,B,\Z}$, Define  $|x|_v= (2\pi)^{-w}z\bar z$ if $v$ is a complex place and $|x|_v= (2\pi)^{-w/2}(z\bar z)^{1/2}$ if $v$ is a real place.
. 

Assume next $v$ is non-archiemdean. Let $|\;\;|_v$ be the metric on $M_{\dR}$ induced by the metric on $C\otimes_F M_{\dR}$ associated to the $O_C$-lattice $D_{\dR,O_C}(M_{et,\Z_p})$ in $C\otimes_F M_{\dR}$. 

Define $|M|:= \prod_v |x|_v^{-1}\in \R_{>0}$ for a non-zero element $x$ of $M_{\dR}$. By the product formula, this is independent of the choice of $x$. 

We have 

\medskip

(1) $|M(r)|=|M|$ ($r\in \Z$) for Tate twists. 

(2) If $M$ is the pure motive $\Z$, $|M|=1$. 

(3) For a finite extension $F'$ of $F$, $|M'|=|M|^{[F':F]}$ where $M'$ denotes the pure motive with $\Z$-coefficients over $F'$ induced by $M$. 

\medskip

A philosophy on motives is that there is a finite extension $F'$ of $F$ such that $M$ becomes isomorphic to $\Z(r)$ over $F'$ for some $r\in \Z$. By (1), (2), (3), this philosophy tells that $|M|$ should be always equal to $1$. 

\end{sbpara}

\begin{sbpara}\label{met3}

In \cite{Ka2}, the author defined the height of $M$  in Case (pure). The definition was improved by Koshikawa in \cite{Ko1}. These definitions use metrics $|\;\;|_v$ on the one-dimensional $F$-vector space 
$$L:=\otimes_{r\in \Z} (\text{det}_F \gr^rM_{dR})^{\otimes r}$$
for places $v$ of $F$. The height is defined as $\prod_v |x|_v^{-1}$ for a non-zero element $x$ of $L$. 
\end{sbpara}

\begin{sbprop}\label{met3}

Assume that $M$ is in Case (pure). Let $H'_{\star}(M)$ be the height of $M$ defined in \cite{Ko1}. Then 
$$H'_{\star}(M) = H_{\star}(M) |\text{det}(M)|^{w/2}.$$

\end{sbprop}

\begin{pf}
Via the isomorphism \ref{alg22}  for $k=F$ and $H=M_{\dR}$, let $x\in L\smallsetminus \{0\}$ and denote the image of $x^{\otimes 2}$ in $\otimes_{r\in \Z} L_r(M)^{\otimes r} \otimes_F (\text{det}_F(M_{\dR}))^{\otimes w}$ as $y\otimes z^{\otimes w}$ where 
$y$ is a non-zero element of $\otimes_{r \in \Z} L_r(M)^{\otimes r}$ and $z$ is a non-zero element of $\text{det}_F (M_{\dR})= (\text{det}(M))_{\dR}$. The definition of $|x|_v$ in \cite{Ko1} shows that 
$$|x|_v = |y|^{1/2}_v |z|_v^{w/2}$$
where $|y|_v$ is as in \ref{defht1} and $|z|_v$ is as in \ref{met23}.  
By taking the product for all $v$, we have the formula in \ref{met3}.
\end{pf}

\begin{sbpara}
The philosophy on motives in \ref{met23} tells that $|\text{det}(M)|$ should be $1$ and hence that $H'_{\star}(M)$ should coincide with $H_{\star}(M)$. 

\end{sbpara}

\begin{sbpara}\label{alg1} The following facts concerning determinant modules are well known. Let $S=(S, \sO_S)$ be a locally ringed space.

(1) For $i=1,2$, let $V_i$ be a vector bundle over $S$ of rank $n(i)$. Fix a bijection $(s, t): \{1, 2, \dots, n(1)n(2)\}\to \{1, \dots, n(1)\}\times \{1, \dots, n(2)\}$. Then concerning $\text{det}_{\sO_S}(V_1\otimes_{\sO_S} V_2)$, we have a unique isomorphism
$$(\text{det}_{\sO_S}(V_1))^{\otimes n(2)} \otimes_{\sO_S} (\text{det}_{\sO_S}(V_2))^{\otimes n(1)}\overset{\cong}\to \text{det}_{\sO_S}(V_1\otimes_{\sO_S} V_2)$$
which sends $(\wedge_{i=1}^{n(1)} a_i)^{\otimes n(2)}\otimes (\wedge_{i=1}^{n(2)} b_i)^{\otimes n(1)}$ to $\wedge_{i=1}^{n(1)n(2)} (a_{s(i)}\otimes b_{t(i)})$ for any local sections $a_i$ of $V_1$ and $b_i$ of $V_2$. 

(2) Let $V$ be a vector bundle over $S$ of rank $n$. Let $m=n(n+1)/2$ and fix a bijection  $(s,t): \{1,\dots, m\}\to \{(i,j)\in \{1,\dots, n\}\times \{1,\dots, n\}\;|\; i\leq j\}$. Then concerning $\text{det}_{\sO_S}(\text{Sym}^2(V))$, we have a unique isomorphism 
$$(\text{det}_{\sO_S}(V))^{\otimes (n+1)} \overset{\cong}\to  \text{det}_{\sO_S}(\text{Sym}^2_{\sO_S}(V))$$
which sends $(\wedge_{i=1}^n a_i)^{\otimes (n+1)}$ to $\wedge_{i=1}^m (a_{s(i)}a_{t(i)})$ for any local sections $a_i$ of $V$. 

(3) Let $V$ be a vector bundle over $S$ of rank $n$. Let $m=n(n-1)/2$ and fix a bijection  $(s,t): \{1,\dots, m\}\to \{(i,j)\in \{1,\dots, n\}\times \{1,\dots, n\}\;|\; i< j\}$. Then concerning $\text{det}_{\sO_S}(\wedge^2_{\sO_S}(V))$, we have a unique isomorphism 
$$(\text{det}_{\sO_S}(V))^{\otimes (n-1)} \overset{\cong}\to  \text{det}_{\sO_S}(\wedge^2_{\sO_S}(V))$$
which sends $(\wedge_{i=1}^n a_i)^{\otimes (n-1)}$ to $\wedge_{i=1}^m (a_{s(i)}\wedge a_{t(i)})$ for any local sections $a_i$ of $V$. 

In the following \ref{pfcl} and \ref{pfsp}, we use only the case $S=\Spec(k)$ for a field $k$ (hence $V_1, V_2, V$ are merely finite-dimensional vector spaces over $k$). But we will use the general case of \ref{alg1} in \ref{Hanalo}. 
\end{sbpara}

\begin{sbpara}\label{pfcl} We prove Proposition \ref{rsc} (1).

Define
$$L^{(1)}:= \otimes_{r<0}  L_r(M^*\otimes M)^{-1}, \quad L^{(2)} := \otimes_{r\in \Z} L_r(M)^{\otimes a(h,r)}$$
where the notation $L_r$ is as in 1.4.3. 
We give an isomorphism $L^{(1)}\cong L^{(2)}$ which preserves the canonical metrics at all places of $F$. \ref{rsc} (1) will follow from it. 

Since $\gr^r(M^*\otimes M)_{\dR}= \oplus_i (\gr^{i-r}M_{\dR})^* \otimes \gr^i M_{\dR}$, we have by \ref{alg1} (1) 
$$L_r(M^*\otimes M) 
\cong \otimes_{i\in \Z} L_i(M) ^{\otimes n(r,i)}$$
where $n(r,i)=\dim_F\gr^{i-r}M_{\dR}- \dim_F\gr^{i+r}M_{dR}$. This isomorphism preserves the metrics at all places of $F$. These isomorphisms for $r<0$  give the isomorphism $L^{(1)}\cong L^{(2)}$ which preserves the metrics  at all places of $F$.

\end{sbpara}

\begin{sbpara}\label{metclub} We give comments on $H_{\clubsuit}$. 
Let $M$ be in Case (pure) and let $L= \text{det}_F((\text{End}_F M_{\dR})/F^0)^{-1}$.

(1) For each place $v$ of $F$, we have a canonical metric $|\;\;|_v$ on $L$ as follows, and we have 
$$H_{\clubsuit}(M)= \prod_v |x|^{-1}_v\quad \text{for } \quad x\in L\smallsetminus \{0\}.$$

Since  $\gr^r(M^*\otimes M)_{\dR}= \oplus_{i\in \Z}  \Hom_F(\gr^iM_{dR}, \gr^{i+r}M_{\dR})$ and $\gr^{-r}(M^*\otimes M)_{dR}= \oplus_{i\in \Z} \Hom_F(\gr^{i+r}M_{\dR}, \gr^iM_{\dR})$ are dual $F$-vector spaces of each other, we have a canonical isomorphism $(\text{det}_F(\gr^r(M^*\otimes M)_{\dR}))^{\otimes 2}\cong L_r(M^*\otimes M)$. Hence we have 
$$L^{\otimes 2} \cong L^{(1)}.$$
Hence the metric on $L^{(1)}$ at $v$ induces the metric $|\;\;|_v$ on $L$.

(2) For a non-archimdean place $v$ of $F$, the above metric $|\;\;|_v$ on $L$ coincides with the metric induced from the $O_C$-lattice in $C\otimes_F L$, where $C$ is the completion of an algebraic closure of $F_v$,  which is induced from the $O_C$-lattice in $C\otimes_F \text{End}_F(M_{\dR})/F^0$ defined as the image of the $O_C$-lattice $D_{\dR,O_C}\text{End}_{\Z_p}(M_{et,\Z_p})$ in $D_{\dR}(\text{End}_{\Q_p}(M_{et,\Q_p}))=C\otimes_F \text{End}_F(M_{\dR})$. Here $p=text{char}(\F_v)$ and $D_{\dR}$ is defined by the local field $F_v$.

\end{sbpara}

\begin{sbpara}\label{pfsp}  We prove Proposition \ref{rsc} (2). 
Let $$L^{(1)}:=P\otimes_F Q \otimes_F R \quad \text{where}\quad 
P=\otimes_{w<0, r<0} L_r(\gr^W_w(M^*\otimes M))^{-1}, $$ $$Q=\otimes_{w:\text{odd}, r<w} L_r(\text{Sym}^2 \gr^W_wM)^{-1}, \quad R= \otimes_{w:\text{even}, r<w} L_r(\wedge^2 \gr^W_wM)^{-1},$$
$$L^{(2)}:= \otimes_{w,r} L_r(\gr^W_wM)^{\otimes b(h, w,r)}.$$
We give an isomorphism $L^{(1)}\cong L^{(2)}$ which preserves the canonical metrics at all places of $F$. \ref{rsc} (2) will follow from it.

This isomorphism comes from the following isomorphisms (1) and (2) which preserve canonical metrics at all places of $F$. 

(1) $\otimes_{w<0,r<0} L_r(\gr^r \gr^W_w(M^*\otimes M))^{-1}\cong  \otimes_{w,r\in \Z} L_r(\gr^W_wM)^{\otimes s(w,r)}$ where $s(w,r)= (\sum_{(w',r')\in A(w,r)} h(w',r'))- (\sum_{(w',r')\in B(w,r)} h(w',r'))$.

(2) $\otimes_{r<w} L_r(\text{Sym}^2 \gr^W_w M)^{-1}$ (resp.  
$\otimes_{r<w} H_r(\wedge^2 \gr^W_w M)^{-1}$) $\cong  \otimes_{r\in \Z} L_r(\gr^W_wM)^{\otimes t(w,r)}$
in the case $w$ is odd (resp. even) where 
 $t(w,r)= (\sum_{(w', r')\in C(w,r)} h(w', r')) +e(w,r)$.

Here $A(w,r),B(w,r), C(w,r)$ are as in \ref{ab}.

We prove (1). By \ref{alg1} (1), we have an isomorphism  
$$L_r((\gr^W_{w_1}M)^* \otimes \gr^W_{w_2}M)  \cong  (\otimes_{i\in \Z} L_i(\gr^W_{w_2}M)^{\otimes m(i)}) \otimes_F (\otimes_{i\in \Z} L_i(\gr^W_{w_1}M)^{\otimes n(i)})$$
where $m(i)= \dim_F \gr^{i-r}\gr^W_{w_1}M$ and $n(i)=- \dim_F \gr^{i+r} \gr^W_{w_2}M_{\dR}$.
This isomorphism  preserves the metrics $|\;\;|_v$ for any $v$.
Hence 
$$H_r((\gr^W_{w_1}M)^* \otimes \gr^W_{w_2}M)= \prod_i H_i(\gr^W_{w_2}M)^{m(i)} \cdot \prod_i H_r(\gr^W_{w_1})^{n(i)}.$$
Since $gr^W_w(M^*\otimes M)= \oplus_{w_2-w_1=w} (\gr^W_{w_1}M)^*\otimes \gr^W_{w_2}M)$, this gives (1).

(2) for $w$ odd (resp. even) is obtained in the similar way by using 
 \ref{alg1} (2) (resp. \ref{alg1} (3)) instead of \ref{alg1} (1).

\end{sbpara}

\begin{sbpara}\label{met1} We give comments on $H_{\spadesuit}$. 

Let $M$ be in Case (mixed-pol) and let  $L=\text{det}_F((\text{End}_{F,W, \langle\;,\;\rangle} M_{\dR})/F^0)^{-1}$.

(1) For each place $v$ of $F$, we have a canonical metric $|\;\;|_v$ on $L$ as follows, and we have 
$$H_{\spadesuit}(M)= \prod_v |x|^{-1}_v\quad \text{for } \quad x\in L\smallsetminus \{0\}.$$ This metric $|\;\;|_v$ is induced by the following isomorphism $L^{\otimes 2}\cong L^{(1)}$ with $L^{(1)}$ as in \ref{pfsp}.

We have
$$L^{-1}\cong (\otimes_{w<0,r<0} R_{w,r})\otimes (\otimes_{w\in \Z, r<0} S_{w,r})\quad \text{where}$$
$$R_{w,r}=\text{det}_F (\gr^r\gr^W_w((M^*\otimes M)_{\dR})), \quad S_{w,r}= \text{det}_F(\gr^r\text{End}_{F, \langle;,\;\rangle}(\gr^W_wM_{\dR})).$$
Since the polarization gives the duality between $\gr^r\gr^W_w((M^*\otimes M)_{\dR})$ and $\gr^{w-r}gr^W_w((M^*\otimes M)_{\dR})$, we have an isomorphism

\smallskip

(1.1) $R_{w,r}^{\otimes 2}\cong  L_r(\gr^W_w(M^*\otimes M))$.

\smallskip

Let $w$ be an odd integer. 
We have an isomorphism
$\text{Sym}^2_F\gr^W_w(M_{\dR}) \overset{\cong}\to \text{End}_{F, \langle,\rangle_w} (\gr^W_wM_{\dR})$
which sends 
$xy\in \text{Sym}^2_F\gr^W_w(M_{\dR})$ ($x,y\in \gr^W_wM_{\dR}$) to  $f\in \text{End}_{F, \langle,\rangle_w} (\gr^W_wM_{\dR})$ defined by 
$f(z)=\langle z, x\rangle_w  y + \langle z, y\rangle_w  x$. Here $\langle\;,\;\rangle_w$ is the pairing defined by the polarization of $\gr^W_wM$. It induces an isomorphism 
$\gr^{w+r}\text{Sym}^2_F\gr^W_w(M_{\dR}) \overset{\cong}\to \gr^r\text{End}_{F, \langle,\rangle_w} (\gr^W_wM_{\dR})$.
Since  $\gr^r\text{End}_{F, \langle,\rangle_w} (\gr^W_wM_{\dR})$ and $\gr^{-r}\text{End}_{F, \langle,\rangle_w} (\gr^W_wM_{\dR})$ are dual of each other, we obtain a canonical isomorphism

\smallskip

(1.2) $S_{w,r}^{\otimes 2} \cong L_{r+w}(\text{Sym}^2\gr^W_wM)$ for $w$ odd. 

\smallskip

Similarly for $w$ even, we have an isomorphism 
$\wedge^2_F\gr^W_w(M_{\dR}) \overset{\cong}\to \text{End}_{F, \langle,\rangle_w} (\gr^W_wM_{\dR})$
which sends 
$x\wedge y\in \wedge^2_F\gr^W_w(M_{\dR})$ ($x,y\in \gr^W_wM_{\dR}$) to  $f\in \text{End}_{F, \langle,\rangle_w} (\gr^W_wM_{\dR})$ defined by 
$f(z)=\langle  z, x\rangle_w y - \langle z, y\rangle_w  x$. 
 It induces an isomorphism 
$\gr^{w+r}\wedge^2_F\gr^W_w(M_{\dR}) \overset{\cong}\to \gr^r\text{End}_{F, \langle,\rangle_w} (\gr^W_wM_{\dR})$.
Since  $\gr^r\text{End}_{F, \langle,\rangle_w} (\gr^W_wM_{\dR})$ and $\gr^{-r}\text{End}_{F, \langle,\rangle_w} (\gr^W_wM_{\dR})$ are dual of each other, we obtain a canonical isomorphism

\smallskip

(1.3) $S_{w,r}^{\otimes 2} \cong L_{r+w}(\wedge^2\gr^W_wM)$ for $w$ even.

\smallskip
The canonical metrics on the right hand sides of (1.1), (1.2), (1.3) induce a canonical metric on $L$. 

\medskip

(2) For an archimedean place $v$ of $F$, the metric $|\;\;|_v$ on $L$ defined in (1) coincides with the Hodge metric defined by the polarizations of $\gr^W_wM$. This fact follows from 1.4.2 and from the fact that the isomorphism $L^{\otimes 2}\cong  L^{(1)}$ preserves the Hodge metrics at $v$ defined by the polarizations as is easily seen.

\medskip
(3) However, for a non-archimedean place $v$ of $F$, the metric $|\;\;|_v$ on $L$ defined in (1) need not coincide with 
the metric $|\;\;|'_v$ induced from the $O_C$-lattice in $C\otimes_F L$, where $C$ is the completion of an algebraic closure of $F_v$,  which is induced from the $O_C$-lattice in $C\otimes_F \text{End}_{F,W,\langle\;,\;\rangle}(M_{\dR})/F^0$ defined as the image of the $O_C$-lattice $D_{\dR,O_C}(\text{End}_{\Z_p,W,\langle\;,\;\rangle}(M_{et,\Z_p}))$ in $C\otimes_{F_v} D_{\dR}(\text{End}_{\Q_p,W, \langle\;,\;\rangle}(M_{et,\Q_p}))=C\otimes_F \text{End}_{F,W,\langle\;,\;\rangle}(M_{\dR})$. Here $p=\text{char}(\F_v)$ and $D_{\dR}$ is defined by the local field $F_v$. See the example below.

\medskip

(4) We can define another  height function $H'_{\spadesuit}(M)$ for $M$ in Case (mixed-pol) by using the above metric $|\;\;|_v'$ for non-archimedean places $v$ and the metric $|\;\;|_v$ for archimedean places of $v$. $H'_{\spadesuit}$ has good properties like $H_{\spadesuit}$, for example, the results on $H_{\spadesuit}$ in \ref{Htprop2} below hold also for $H'_{\spadesuit}$. In this paper, we use $H_{\spadesuit}(M)$, not $H'_{\spadesuit}(M)$, because the former is independent of the choices of polarizations on $\gr^W_wM$ and expressed by using $H_r(\gr^W_wM)$ and hence the arguments become simpler.

{\bf Example.} Let $F=\Q$, $M=\Z\oplus \Z(1)$, take $a, b\in \Q_{>0}$, and consider the polarization on $\gr^W_wM$ which is the direct sum of $a$ times the evident morphism $\Q\otimes \Q\to \Q$ for $w=0$ and $b$ times the evident morphism $\Q(1)\otimes \Q(1)\to \Q(2)$ for $w=-2$. Then $L$ is canonically isomorphic to $\Q$ with a canonical base $e$.  We have $|e|_v=|a/b|^{1/2}_v$ for any place $v$ of $\Q$, $|e|'_v=|e|_v$ for the archimedean place of $\Q$, but $|e|'_v=1$ for any non-archimedean place $v$ of $\Q$. Hence $H_{\spadesuit}(M)=1$ (by the product formula), but $H'_{\spadesuit}(M)= (b/a)^{1/2}$. 

\end{sbpara}

\begin{sbprop}\label{Htprop2}

(1) Concerning the dual, we have $$H_{\star}(M^*)=H_{\star}(M), \quad H_{\clubsuit}(M^*)= H_{\clubsuit}(M), 
\quad H_{\spadesuit}(M^*)=H_{\spadesuit}(M)$$ for $M$ in Case (pure), Case (pure), Case (mixed), respectively.

(2) Concerning Tate twists, we have $$H_{\star}(M(r))=H_{\star}(M), \quad H_{\clubsuit}(M(r))= H_{\clubsuit}(M), \quad H_{\spadesuit}(M(r))=H_{\spadesuit}(M)$$ for $M$ in Case (pure), Case (pure), and Case (mixed), respectively.

(3) For a finite extension $F'$ of $F$, we have 
$$H_{\star}(M')=H_{\star}(M)^{[F':F]}, \quad H_{\clubsuit}(M')=H_{\clubsuit}(M)^{[F':F]}, \quad H_{\spadesuit}(M')=H_{\star}(M)^{[F':F]}$$
for $M$ in Case (pure), Case (pure), Case (mixed), respectively, where $M'$ denotes the motive over $F'$ induced by $M$.

(4) Let $n\in \Q_{>0}$. For $M$ in Case (pure) (resp. Case (mixed)), if $M'$ is obtained from $M$ by changing $M_{et, \hat \Z}$ as $M'_{et,\hat \Z}= n \cdot M_{et, \hat \Z}$ (but $M'$ with $\Q$-coefficients is the same as $M$ with $\Q$-coefficients), then $H_{\star}(M')=H_{\star}(M)$ and $H_{\clubsuit}(M')=H_{\clubsuit}(M)$ (resp. $H_{\spadesuit}(M')=H_{\spadesuit}(M)$). 

(5) For direct sums, for $M$ and $M'$  in Case (pure), $H_{\star}(M\oplus M')= H_{\star}(M)H_{\star}(M')$.

\end{sbprop}

\begin{pf} This follows from the result \ref{dtetc} on the height function $H_r$. 
\end{pf}

\begin{sbrem}
This is a remark related to the above (4) in \ref{Htprop2}.

Koshikawa proved that the set $$\{H_{\star}(M')\;|\;\text{$M'$ is isogenous to $M$}\}$$ is finite. The author expects that his method works to have the same result for $H_r$ and hence for $H_{\clubsuit}$ and $H_{\spadesuit}$. 
\end{sbrem}

\begin{sbpara}

{\bf Example.} Let $A$ be an abelian variety of dimension $g$ over a number field $F$ and let $M$ be the pure motive $H^1(A)$ with $\Z$-coefficients over $F$. Then $H_{\star}(M)$ is the usual height of $A$ of Faltings if $A$ is of semi-stable reduction (this is proved in \cite{Ko1} modulo some power of $2$, but this power of $2$ disappears by the work \cite{Ce}). We have
  $$H_1(M)=H_{\star}(M), \quad H_0(M)=H_{\star}(M)^{-1}, \quad H_r(M)=1\;\;\text{for}\;\; r\neq 0,1,$$ $$H_{\clubsuit}(M)=H_{\star}(M)^{2g}, \quad 
    H_{\spadesuit}(M)= H_{\star}(M)^{g+1}.$$
     \end{sbpara}

\subsection{Reviews on relative monodromy filtrations and the splittings of Deligne}\label{prim}
We review relative monodromy filtrations and related subjects.

\begin{sbpara}\label{2.4.b}

Let $A$ be a ring and let $\Upsilon$ be an $A$-module  endowed with a finite increasing filtration $W$ by $A$-submodules, and let $N: \Upsilon \to \Upsilon$ be a nilpotent $A$-homomorphism such that $NW_w\subset W_w$ for any $w\in \Z$. A finite  increasing filtration $\sW$ on $\Upsilon$ by $A$-submodules is called {\it a relative monodromy filtration of $N$ with respect to $W$} if it satisfies the following two conditions (i) and (ii).

(i) $N\sW_w\subset \sW_{w-2}$ for any $w\in \Z$.

(ii) For any $w\in \Z$ and any integer $m\geq 0$, we have an isomorphism  $N^m: \gr^{\sW}_{w+m}\gr^W_w \overset{\cong}\to \gr^{\sW}_{w-m}\gr^W_w$.

A relative monodromy filtration $\sW$ of $N$ with respect to $W$ is unique if it exists (\cite{De} 1.6.13).

\end{sbpara}

\begin{sbpara}
If the relative monodromy filtration $\sW$ exists, the filtration on $\Hom_A(\Upsilon, \Upsilon)$ induced by $\sW$ is a relative monodromy filtration of $\text{Ad}(N): h \mapsto Nh-hN$ with respect to the filtration on $\Hom_A(\Upsilon, \Upsilon)$ induced by $W$. 

\end{sbpara}

\begin{sbpara}\label{pi} The relative monodromy filtration exists if the filtration $W$ is pure, that is, if $W_w=\Upsilon$ and $W_{w-1}=0$ for some $w\in\Z$. 

Hence in general, the relative monodromy filtration $\sW$ on $\gr^W_w$ exists and we can discuss $\gr^{\sW}_m\gr^W_w$.

For $w\in \Z$ and for any integer $m\geq 0$, we have a direct sum decomposition
$$\gr^{\sW}_{w+m}\gr^W_w= (\gr^{\sW}_{w+m}\gr^W_w)_{\pri}\oplus N(\gr^{\sW}_{w+m+2}\gr^{\sW}_w)$$
where $(\gr^{\sW}_{w+m}\gr^W_w)_{\pri}$ is defined as 
the kernel of $N^{m+1}: \gr^{\sW}_{w+m}\gr^W_w\to \gr^{\sW}_{w-m-2}\gr^W_w$ and is called the primitive part of $\gr^{\sW}_{w+m}\gr^W_w$ (\cite{De}).

\end{sbpara}

\begin{sbpara} We now use the theory of the splitting of Deligne (\cite{SC}) which we review in \ref{NdSC} below. This theory
tells, roughly speaking, that a splitting of the relative monodromy filtration $\sW$ (satisfying some conditions) extends canonically to a double splitting of the pair $(W, \sW)$. 

\end{sbpara}

\begin{sbpara}\label{Assume12}
We now assume the following (1) and (2). 

\medskip

(1) The relative monodromy  filtration $\sW$ exists.

\medskip

(2)  There is a decomposition $$\Upsilon=\oplus_{m\in \Z} \Upsilon_m$$
as an $A$-module satisfying the following conditions (i)--(iii).

(i) $\sW_m= \oplus_{i\leq m} \Upsilon_i$ for any $m\in \Z$,

(ii) $N(\Upsilon_m)\subset \Upsilon_{m-2}$ for any $m\in \Z$.

(iii) $W_w=\oplus_m W_w\cap \Upsilon_m$ for any $w\in \Z$.
\medskip

Then for  $d\geq 2$, as is shown in \ref{Nindep} below, we have a canonically defined element $$N_d\in (\gr^{\sW}_{-2}\gr^W_{-d}\Hom_A(\Upsilon, \Upsilon))_{\pri}$$ 
which should be called the component of $N$ of weight $-d$. This $N_d$ will play an important role in \S1.7.

\end{sbpara}

\begin{sbpara}\label{NdSC} Assume (1) (2) in \ref{Assume12}. For a decomposition of $\Upsilon$ as in (2) in \ref{Assume12}, by the theory of Deligne splitting in 
\cite{SC}, there is a unique decomposition $\Upsilon=\oplus_{i,j} \Upsilon_{i,j}$ as an $A$-module having the following properties (i)--(iii).

(i) $\Upsilon_m=\oplus_i \Upsilon_{i,m}$ for any $m$.

(ii)  $W_w=\oplus_{i,m,i\leq w} \Upsilon_{i,m}$ for any $w$.

(iii)  In the isomorphism $$\Hom_A(\Upsilon,\Upsilon)\cong \oplus_{w,m} \;\gr^{\sW}_m\gr^W_w\Hom_A(\Upsilon,\Upsilon)$$
given by this decomposition $(\Upsilon_{i,j})_{i,j}$ of $\Upsilon$, write the image of $N\in \Hom_A(\Upsilon,\Upsilon)$ in $\oplus_{w,m} \;\gr^{\sW}_m\gr^W_w\Hom_A(\Upsilon,\Upsilon)$ as $\sum_{d\geq 0} N_d$ with $N_d\in \gr^{\sW}_{-2}\gr^W_{-d}\Hom_A(\Upsilon, \Upsilon)$. Then $N_1=0$ and 
 $N_d\in (\gr^{\sW}_{-2}\gr^W_{-d}\Hom_A(\Upsilon, \Upsilon))_{\pri}$ for $d\geq 2$.

\end{sbpara}

\begin{sbprop}\label{Nindep}

This element  $N_d$ of $(\gr^{\sW}_{-2}\gr^W_{-d}\Hom_A(\Upsilon,\Upsilon))_{\pri}$ is independent of the choice of the decomposition of $\Upsilon$ as in (2) of \ref{Assume12}. 

\end{sbprop}

\begin{pf}  Assume we are given another decomposition $\Upsilon=\oplus_i \Upsilon_i'$ as in (2) of \ref{Assume12}. 
Let $g$ be the unique automorphism of $(\Upsilon,\sW)$ which induces the identity map on $\gr^{\sW}$ such that $g\Upsilon_i=\Upsilon_i'$ for all $i$. Then $gW_w=W_w$ for all $w$ and $gN=Ng$. For the Deligne's splitting $\Upsilon=\oplus_{i,j} \Upsilon_{i,j}'$ associated to $(\Upsilon_i')_i$, 
we have $\Upsilon_{i,j}'=g\Upsilon_{i,j}$. If $N_d'$ denotes the element of $(\gr^{\sW}_{-2}\gr^W_{-d} \Hom_A(\Upsilon,\Upsilon))_{\pri}$ defined by this  $(\Upsilon'_i)_i$, we have $N_d'=\Ad(g)N_d$, but $\Ad(g)$ acts on $\gr^{\sW}_{-2}\gr^W_{-d}\Hom_A(\Upsilon, \Upsilon)$ trivially. Hence $N_d'=N_d$. 
\end{pf}

This $(N_d)_{d\geq 2}$ tells how the pair $(W, N)$ is far from being a direct sum of pure objects: 
\begin{sblem} Assume (1) and (2) in \ref{Assume12} and assume that $A$ is a field. 
Then $N_d=0$ for all $d\geq 2$ if and only if there is an $A$-linear splitting $\Upsilon\cong \gr^W$ of $W$ via which $N$ corresponds to $(\gr^W_w(N)_w$ on $\gr^W$.  

\end{sblem}
This is seen easily

\subsection{Hodge analogues}

In this \S1.6, we consider the analogy between the following two subjects. 

\medskip

(i) A mixed motive $M$ over a number field. 

(ii) A  variation $\sH$ of mixed Hodge structure on a smooth curve over the complex number field $\C$. 

\medskip

The Hodge analogue $h_{\star}(\sH)$ of the height $H_{\star}(M)$   was considered in an old paper Griffiths \cite{Gri}.  In this \S1.6, we explain it and a  generalization $h_{\star\diamondsuit}(\sH)$ of it  to $\sH$ mixed. 
We also discuss Hodge analogues $h_r(\sH)$ ($r\in \Z$), $h_{\clubsuit}(\sH)$ and $h_{\spadesuit}(\sH)$ of the height functions  $H_r(M)$ ($r\in \Z$), $H_{\clubsuit}(M)$ and $H_{\spadesuit}(M)$, respectively.

\begin{sbpara}\label{Hodgeside} Let $C$ be a projective smooth connected curve over $\C$ considered analytically and let $\sH=(\sH_\Q, W, F)$ be a  variation of mixed $\Q$-Hodge structure on  $U= C\smallsetminus R$ where $R$ is a finite subset of $C$. Here $\sH_\Q$ is the local system on $U$ with $\Q$-coefficients, $W$ is the weight filtration, and $F$ is the Hodge filtration.

We assume the following (i)--(iii).

\medskip

(i) The pure graded quotients $\gr^W_w\sH$ are polarizable for any $w\in \Z$.

(ii) Let $x\in R$. Then
the local monodromy of $\sH_\Q$ at $x$ is unipotent. Furthermore,  the local monodromy logarithm at $x\in R$ has a relative monodromy filtration with respect to $W$. 

More precisely, let $x\in C$ and take an open neighborhood $\Delta$ of $x$ in $C$ such that $\Delta\cong \{z\in \C\;|\;|z|<1\}$ and such that $\Delta^*= \Delta\smallsetminus \{x\}\subset U$, and denote by 
$\sH_{\Q,x}$the stalk of $\sH_\Q$ at a point $y$ of $\Delta^*$.  Then the local monodromy group $\pi_1(\Delta^*)$ at $x$, which is canonically isomorphic to $\Z(1)=\Z\cdot 2\pi i $, acts on $\sH_{\Q,x}$. $\sH_{\Q,x}$  depends on the choice of $y$ but the stalk at $y$ and that at $y'\in \Delta^*$ are  isomorphic via an isomorphism which is canonical modulo the action of the local monodromy group at $x$. The first part of (ii) says that this action of the local monodromy group at $x$ on $\sH_{\Q,x}$ is unipotent. Let $N'_x: H_{\Q,x}\to H_{\Q,x}$ be the logarithm of the action of the canonical generator of the local monodromy group at $x$. Then the second part of (ii) says that the relative monodromy filtration of $N'_x$ with respect to the weight filtration $W$ on $\sH_{\Q,x}$ (\ref{2.4.b}) exists.

Note that by the unipotence of the local monodromy in the condition (i), by the theory of canonical extension of Deligne, we have a vector bundle $\sH_{\sO}$ on $C$ which extends the vector bundle $\sO_U \otimes_\Q \sH_\Q$ on $U$ canonically.

(iii) The Hodge filtration on $\sH_{\sO}|_U$ extends to subbundles $F^p=F^p\sH_{\sO}$ of $H_{\sO}$.

\end{sbpara}

\begin{sbrem} (1) It is known that variations of mixed Hodge structures on $U$ of geometric origin satisfy these conditions (i)--(iii) if the unipotence in (ii) is replaced by quasi-unipotence. The unipotence at any $x$ is similar to the semi-stability \ref{semi-st} of a motive over a number field. 
The reason why we put the condition of the unipotence, not quasi-unipotence, will be explained in \ref{reasons} (2). 

(2) Everything in this \S1.6 works for variations of $\R$-Hodge structure. But we consider $\Q$-Hodge structure,  for it is nearer to motives. 

\end{sbrem}

\begin{sbpara}\label{Hanalo} 
If $\sH$ is pure, 
$$h_r(\sH):= \text{deg}(\gr^r \sH_{\sO})= \frac{1}{2}\cdot (\text{deg}(\gr^r\sH_{\sO})- \text{deg}(\gr^{w-r}\sH_{\sO}))\quad 
(r\in \Z), $$
$$h_{\star}(\sH):=\sum_r  \; r\cdot h_r(\sH)=\text{deg}(\otimes_{r\in \Z} (\text{det}_{\sO_C} \gr^r \sH_{\sO})^{\otimes r})= \sum_{r \in \Z} \text{deg} F^r\sH_{\sO},$$
$$h_{\clubsuit}(\sH):=-\text{deg}(\sE nd_{\sO_C}(\sH_{\sO})/F^0)$$ 
 are similar to $H_r(M)$, $H_{\star}(M)$ and $H_{\clubsuit}(M)$  for $M$ in Case (pure), respectively. 
 Here the second $=$ in the definition of $h_r(\sH)$ is due to the perfect duality $\gr^r \sH_{\sO} \times \gr^{w-r}\sH_{\sO} \to \sO_C$ given by a polarization of $\sH$. 
    
 If $\gr^W_w\sH$ ($w\in \Z$) are polarized, 
 $$h_{\spadesuit}(\sH):= -\text{deg} (\sE nd_{\sO_C,W, \langle\;,\;\rangle}(\sH_{\sO})/F^0).$$ 
is similar to 
$H_{\spadesuit}(M)$ of Case (mixed).

The height functions $h_{\clubsuit}$, and $h_{\spadesuit}$ are expressed by using $h_r$ as 
$$h_{\clubsuit}(\sH)=\sum_r a(h, r)h_r(\sH),$$ 
$$h_{\spadesuit}(\sH)=\sum_{w,r} b(h, w,r) h_r(\gr^W_w\sH),$$ 
where $a(h, r)$ is defined as in \ref{ab} (1) with $h(r)= \text{rank}_{\sO_C}\gr^r H_{\sO}$ and $b(h, w,r)$ is defined as in \ref{ab} (2) with $h(w, r)=\text{rank}_{\sO_C} \gr^r\gr^W_w \sH_{\sO}$. These formulas are shown using the isomorphisms (1), (2), (3) in \ref{alg1}.

In particular, $h_{\spadesuit}(\sH)$ is independent of the choices of the polarizations of $\gr^W_w\sH$ ($w\in \Z$). 

\end{sbpara}

\begin{sbrem}  In Griffiths \cite{Gri}, the line bundle  $\otimes_{r\in \Z} \text{det}_{\sO_C} (\gr^r \sH_{\sO})^{\otimes r}$ is called the canonical bundle. When the author wrote the paper \cite{Ka2}, he did not know that his definition of the height of a pure motive in \cite{Ka2} was similar to the above $h_{\star}(\sH)$ and regrets that he did not refer to \cite{Gri} in \cite{Ka2}. 
 The author learned this analogy from Koshikawa as is written in \cite{Ka3}.

 \end{sbrem}

   \begin{sbpara}\label{hmix} If the pure graded quotients of $\sH$ are polarized, the Hodge analogue $h_{\star\diamondsuit}(\sH)$ of $H_{\star\diamondsuit}(M)$ with $M$ mixed (\S1.7) is defined as 
 $$h_{\star\diamondsuit}(\sH)=(\sum_{w\in \Z} h_{\star}(\gr^W_w\sH))+ h_{\diamondsuit}(\sH)$$
 where $h_{\diamondsuit}(\sH)$ tells, roughly speaking, how $\sH$ is far from being the direct sum of $\gr^W_w\sH$ ($w\in \Z$). This $h_{\diamondsuit}(\sH)$ has the shape
 $$h_{\diamondsuit}(\sH)= \sum_{w\in \Z, d\geq 1}  h_{\diamondsuit, w, d}(\sH),$$
$$h_{\diamondsuit,w,1}(\sH) =  \langle \alpha_w,\alpha_w\rangle,$$
and for $d\geq 2$, $$h_{\diamondsuit, w, d}(\sH)= \sum_{x\in C}  h_{\diamondsuit, w, d,x}(\sH),\quad 
h_{\diamondsuit, w, d,x}(\sH)=  (\langle N_{x,w, d}, N_{x,w, d}\rangle_{N_{x,0}})^{1/d}.$$

 We will define each term of $h_{\diamondsuit}(\sH)$. A rough explanation is as follows:
 $\alpha_w$ is the extension class of the exact sequence $0\to \gr^W_{w-1}\sH\to W_w\sH/W_{w-2}\sH\to \gr^W_w\sH\to 0$,
 $\langle\;,\;\rangle$ is the height pairing in Beilinson \cite{Be2}, and we have  $\langle \alpha_w, \alpha_w\rangle\in \Q_{\geq 0}$. On the other hand, for $d\geq 2$, $h_{\diamondsuit, w, d,x}(\sH)$ is, roughly speaking,  the ''length'' of the $\Hom(\gr^W_w, \gr^W_{w-d})$-component of the monodromy logarithm at $x$. 
  
 We give the definition of $h_{\diamondsuit, w, d,x}(\sH)$ for $d\geq 2$ in  \ref{Na} after preparations \ref{Nb}--\ref{Nd}. 
  We give the definition of 
 $h_{\diamondsuit,w, 1}(\sH)$ in \ref{d=1a} after preparations in \ref{d=1b}--\ref{d=1c}.

    \end{sbpara}

  \begin{sbpara}\label{Nb} We recall the theory of degeneration of Hodge structure. (See \cite{CKS},  \cite{KM}, \cite{SW}, \cite{SZ}.)
  
  Assume $\sH$ is pure of weight $w$, and let $x\in C$. Let 
 $\sW$ be the relative monodromy filtration of $N'_x: \sH_{\Q,x}\to \sH_{\Q,x}$ (\ref{Hodgeside})  with respect to $W$.  Let $\sH_{\sO}(x)$ be the fiber of the vector bundle $\sH_{\sO}$ at $x$ and let $F(x)$ be the Hodge filtration on $\sH_{\sO}(x)$. 
   
  The theory of nilpotent orbit in degeneration (\cite{SW}) shows that we have an isomorphism
 $\C\otimes_\Q \sH_{\Q,x}\cong \sH_{\sO}(x)$ which is canonical modulo the actions of $\exp(zN'_x)$ ($z\in \C$). When we identify these vector spaces by this isomorphism, we have
 
 \medskip
 
 (i) $(\sH_{\Q,x}, \sW, F(x))$ is a mixed $\Q$-Hodge structure. 
 
 \medskip
 
 Furthermore, $N'_xF^p(x)\subset F^{p-1}(x)$ for any $p\in \Z$. Hence we have a homomorphism
 $$N_x:=(2\pi i)^{-1}N'_x\;:\; (\sH_{\Q,x}, \sW, F(x)) \to (\sH_{\Q,x}, \sW, F(x))(-1)$$ of mixed $\Q$-Hodge structures. 
 
 The induced map $N_x: \sH_{\sO}(x)\to \sH_{\sO}(x)$ coincides with the composition $\sH_{\sO}(x) \overset{\nabla}\to \sH_{\sO}(x) \otimes_{\C} \Omega^1_C(\log R)(x) \to \sH_{\sO}(x)$
 where the first arrow is induced by the connection $\nabla: \sH_{\sO}\to \sH_{\sO}\otimes_{\sO_C}\Omega^1_C(\log R)$ and the second arrow is the residue map.

\end{sbpara}

\begin{sbpara}\label{Nc} In \ref{Nb}, assume further that $\sH$ is polarized. 
Let $m\geq 0$. Then the polarization 
 $\sH_\Q\otimes \sH_\Q \to \Q(-w)$ induces a non-degenerate $\Q$-bilinear form $\langle \; ,\;\rangle: \gr^{\sW}_{w-m}\sH_{\Q,x} \otimes \gr^{\sW}_{w+m}\sH_{\Q,x}\to \Q(-w)$. It induces a non-degenerate symmetric bilinear form 
$$ \langle \;,\;\rangle_{N_x}: \gr^{\sW}_{w+m}\sH_{\Q,x}\times \gr^{\sW}_{w+m}\sH_{\Q,x}\to \Q(-w-m),$$
 $$\langle a, b\rangle_{N_x}:= \langle N_x^ma, b\rangle_w$$
By \cite{SW} Theorem 6.16,  the restriction of this to
   $$\langle \;,\;\rangle_{N_x}: (\gr^{\sW}_{w+m}\sH_{\Q,x})_{\pri}\times (\gr^{\sW}_{w+m}\sH_{\Q,x})_{\pri}\to \Q(-w-m),$$
 gives a polarization of the Hodge structure $(\gr^{\sW}_{w+m}\sH_{\Q,x})_{\pri}$ of weight $w+m$.

\end{sbpara}

\begin{sbpara}\label{Nd} Let $\sH$ be as in \ref{Hodgeside} and assume that the pure graded quotients of $\sH$ are polarized. Let $w\in \Z, d\geq 2$.

Let $$\sH'=(\gr^W_w\sH)^* \otimes \gr^W_{w-d}\sH\subset  \gr^W_{-d}(\sH^*\otimes \sH),$$ whose underlying $\Q$-local system is $\sH om_{\Q}(\gr^W_w\sH_\Q, \gr^W_{w-d}\sH_\Q)$. Then $\sH'$ is pure of weight $-d$ and is polarized.  Let 
$$P= \Hom_{\text{HS}}(\Q, (\gr^{\sW}_{-2}\sH')_{\pri}(-1)) =  (\gr^{\sW}_{-2}\sH'_{\Q,x})_{\pri}(-1)\cap F^0((\gr^{\sW}_{-2}\sH'_{\C,x})_{\pri}(-1))$$ $$\subset \gr^{\sW}_{-2}\Hom_{\Q}(\gr^W_w\sH_{\Q,x}, \gr^W_{w-d}\sH_{\Q,x})_{\pri}(-1)$$
where $\Hom_{\text{HS}}$ denotes the set of homomorphisms of Hodge structures. Then the restriction
$$\langle \;,\;\rangle_{N_{x,0}}: P \times P \to \Q$$
of $\langle \;,\;\rangle_{N_{x,0}}$ in \ref{Nc} (we take $\sH'$ as the $\sH$ in \ref{Nc} and denote $\gr^W(N_x)$ by $N_{x,0}$)  is a positive definite symmetric bilinear form because this restriction coincides with the restriction of the positive definite Hodge metric in \ref{Nc}.

\end{sbpara}

\begin{sbpara}\label{Na} By \cite{CKS} and \cite{KM}, the condition (2) in \ref{Assume12} is satisfied by $\Upsilon= \sH_{\Q,x}$ with the weight filtration $W$ and  $N=N'_x$ (\ref{Hodgeside}).  
Hence 
$$N_{x,d}\in 
 (\gr^{\sW}_{-2}\gr^W_{-d}\Hom_{\Q}(\sH_{\Q,x}, \sH_{\Q,x}))_{\text{prim}}(-1)$$ is defined.
 Write 
 $$N_{x,d}= \sum_{w\in \Z} N_{x,w,d}$$
 where $N_{x,w,d}$ is the $(\gr^{\sW}_{-2}\Hom(\gr^W_w\sH_\Q, \gr^W_{w-d}\sH_\Q))_{\pri}(-1)$-component of $N_{x,d}$. Then
 $$ N_{x,w,d}\in P\subset   (\gr^{\sW}_{-2}\sH'_{\Q,s})_{\pri}(-1).$$

Hence $\langle N_{x,w,d}, N_{x,w,d}\rangle_{N_{x,0}}\in\Q_{\geq 0}$ is defined. 
It is $0$ if and only if $N_{x,w,d}=0$ (since the form $\langle\;,\;\rangle_{N_{x,0}}$  on $P$ is positive definite).

Define $$h_{\diamondsuit,w,d,x}(\sH):= (\langle N_{x,w,d}, N_{x,w,d}\rangle_{N_{x,0}})^{1/d}.$$
(It may seem strange to take the $d$-th root here. The reason of this will be explained in \ref{reasons} (1).) We have

\medskip

(i) $h_{\diamondsuit,w,d,x}(\sH)=0$ if and only if $N_{x,w,d}=0$.

\medskip

For $x\in U$, we have $N_x=0$ and we have
$h_{\diamondsuit,w,d,x}(\sH)=0$. 
\end{sbpara}

 \begin{sbpara}\label{d=1b} We give a preparation for the definition of $h_{\diamondsuit, w,1}(\sH)$.

 Let $(C, \sH)$ be 
 as  in \ref{Hodgeside} and assume that $\sH$ is pure of weight $-1$. Let $j:U\overset{\subset}\to C$ be the inclusion map. By Zucker \cite{Zu}, $H^1(C,j_*\sH_\Q)$ is endowed with a structure of a pure Hodge structure of weight $0$.
   
 The Hodge filtration on $\C\otimes_\Q H^1(C,j_*\sH_\Q)$ is given as follows (\cite{Zu}). Let $\sG= \sH_{\sO} \otimes_{\sO_C} \Omega^1_C$ and let $\sF\subset \sH_{\sO}$ be the inverse image of $\sG$ under the connection $\nabla: \sH_{\sO}\to \sH_{\sO}\otimes_{\sO_C} \Omega^1_C(\log R)$. Then $\sF$ is a vector bundle. We have an exact sequence $0\to j_*\sH_{\C} \to \sF \overset{\nabla}\to \sG\to 0$ of sheaves on $C$ and thus $\C\otimes_\Q H^1(C, j_*\sH_\Q)= H^1(C, [\sF\overset{\nabla}\to \sG])$. (Here $H^1(C, [\dots])$ denotes the hyper-cohomology of the complex $[\dots]$.) Let $F^p\sF:=\sF\cap F^p\sH_{\sO}$ and $F^p\sG:= F^{p-1}\sH_{\sO}\otimes_{\sO_C} \Omega^1_C$. Then the Hodge filtration in question is given by $F^p(\C\otimes_\Q H^1(C, j_*\sH_\Q))= H^1(C, [F^p\sF \overset{\nabla}\to F^p\sG])$.

  \end{sbpara}
  
  \begin{sbpara}\label{d=1bc} 
  
Let the notation be as in \ref{d=1b}. 
 
 Let  $\Ext^1(\Q, \sH)$ be the 
 abelian group of extension classes in the category of all objects in \ref{Hodgeside} ($C$ is fixed but $U$ moves).  We have
 $$\Ext^1(\Q, \sH)\cong 
  \text{Ker}(H^0(C, j_*\sH_\Q\bs \sF/F^0\sF) \overset{\nabla}\to H^0(C, \sG/F^0\sG))$$
 where  $\nabla$ is the connection. For an element $a$ of $\Ext^1(\Q, \sH)$, we have an element of the right hand side of this isomorphism as follows. Let $0\to \sH \to \tilde \sH\to \Q\to 0$ be the exact sequence corresponding to 
 $a$ on some $U'=U\smallsetminus R'$ with $R'$ a finite subset of $U$.  
 Locally on $U'$, lift $1\in \Q$ to a local section $a_1$ of $F^0\tilde \sH_{\sO}$ and also to a local section $a_2$ of $\tilde \sH_\Q$. Then $a_1-a_2$ gives a  well defined section on $U'$ of the sheaf $\sH_\Q\bs \sH_{\sO}/F^0$. By using the Griffiths transversality  and the condition  (ii) in \ref{Hodgeside} for $\tilde \sH$, we see that this section belongs to the right hand side of the above isomorphism. 
 
 Consider the long exact sequence of cohomology groups associated to the exact sequence of complexes of sheaves on $C$
 $$0\to j_*\sH_\Q \to [\sF/F^0\sF\to \sG/F^0\sG]  \to [j_*\sH_\Q\bs \sF/F^0\sF \to \sG/F^0\sG]\to 0.$$
  Since $H^1(C, [\sF/F^0\sF \to \sG/F^0\sG])= (\C\otimes_\Q H^1(C, j_*\sH_\Q))/F^0$, the kernel of $H^1(C, j_*\sH_\Q) \to H^1(C, [\sF/F^0\sF \to \sG/F^0\sG])$ is identified with 
  $$P:= \Hom_{\text{HS}}(\Q, H^1(C, j_*\sH_\Q))= H^1(C, j_*\sH_\Q) \cap F^0(\C\otimes_\Q H^1(C, j_*\sH_\Q))$$ where $\Hom_{\text{HS}}$ denotes the set of homomorphisms of Hodge structures.  
 Hence the long exact sequence of cohomology groups gives an exact sequence
 
 \medskip
 
 (1) $0 \to H^0(C, j_*\sH_\Q) \to H^0(C, [\sF/F^0\sF \to \sG/F^0\sG]) \to \text{Ext}^1(\Q, \sH) \to P \to 0$.

 \end{sbpara}
  
 \begin{sbpara}\label{d=1c} Let the notation be as in \ref{d=1b}. Assume furthermore that $\sH$ is polarized. 

The polarization $\sH_\Q\otimes \sH_\Q\to \Q(1)$ and the cup product induce a symmetric $\Q$-bilinear form
$$\langle\;,\;\rangle: H^1(C, j_*\sH_\Q) \times H^1(C, j_*\sH_\Q) \to H^2(C, \Q(1))=\Q.$$
 On $P\subset H^1(C, j_*\sH_\Q)$ (\ref{d=1bc}), this pairing 
 coincides with the positive definite Hermitian Hodge metric given by the polarization of the Hodge structure of $H^1(C, j_*\sH_\Q)$. Hence this cup product $$\langle\; ,\;\rangle : P\times P\to \Q$$ is a positive definite symmetric $\Q$-bilinear form.

\end{sbpara}

  \begin{sbpara}\label{d=1a} Let $\sH$ be as in \ref{Hodgeside} and assume that the pure graded quotients of $\sH$ are polarized. 
   
 Let $\alpha_w\in \Ext^1(\gr^W_w\sH, \gr^W_{w-1}\sH)$ be the  class of the exact sequence $0 \to \gr^W_{w-1}\sH\to W_w\sH/W_{w-2}\sH\to \gr^W_w\sH\to 0$ where $\Ext^1$ is the abelian group of extension classes in the category of all objects in \ref{Hodgeside} ($C$ is fixed and $U$ varies). Let $\sH'=  (\gr^W_w\sH)^* \otimes \gr^W_{w-1}\sH\subset \gr^W_{-1}(\sH^*\otimes \sH)$. Then $\sH'$  is pure of weight $-1$ and is polarized. We have 
  $\Ext^1(\gr^W_w, \gr^W_{w-1})= \Ext^1(\Q, \sH')$.  Let $\beta_w\in P$ be the image of $\alpha_w$ under the map $\text{Ext}^1(\Q,\sH') \to P$ (\ref{d=1bc}). Here we take $\sH'$ as $\sH$ in \ref{d=1bc}. 
  We define $$h_{\diamondsuit, w,1}(\sH)=\langle \alpha_w, \alpha_w\rangle:= \langle \beta_w, \beta_w\rangle\in \Q_{\geq 0}.$$ This is zero if and only if $\beta_w$ is zero. 
 \end{sbpara}

\begin{sbprop}\label{CC'}
Let $C'\to C$ be the integral closure of $C$ in a finite extension of the function field of $C$ and let $\sH'$ be the pullback of $\sH$ to $C'\times_C U$. Then
$$h_r(\sH')=[C':C]h_r(\sH), \quad h_{\star}(\sH')=[C':C]h_{\star}(\sH), \quad h_{\clubsuit}(\sH')=[C':C]h_{\clubsuit}(\sH),$$ $$h_{\spadesuit}(\sH')=[C':C]h_{\spadesuit}(\sH),$$ $$h_{\star\diamondsuit}(\sH')=[C':C]h_{\star\diamondsuit}(\sH), \quad h_{\diamondsuit}(\sH')=[C':C]h_{\diamondsuit}(\sH)$$ 
where $\sH$ is pure in the first line, mixed in the second line, and mixed with polarized pure graded quotients in the third line. We have similar formulas for $h_{\diamondsuit, w,1}$ and $h_{\diamondsuit, w,d,x}$ for $d\geq 2$. 
\end{sbprop} 

\begin{pf}
All formulas except the last two are easily seen. It is sufficient to prove the last formula. The case $d=1$ is easy. We assume $d\geq 2$ and we consider $h_{w,d,x}$.  If $x'\in C'$ lies over $x\in C$ and if the ramification index of $C'\to C$ at $x'$ is $e(x',x)$, the standard generator of the local monodromy group at $x'$ is sent to the $e(x',x)$-th power of that at $x$, and hence 
their logarithms satisfy $N_{x'}= e(x',x)N_x$. Since $\langle N_{x,w,d}, N_{x,w,d}\rangle_{N_{x,0}}= \langle \Ad(N_{x,0})^{d-2}(N_{x,w,d}), N_{x,w,d}\rangle$, we have $\langle N_{x',w,d}, N_{x',w,d}\rangle_{N_{x',0}}= e(x',x)^d\cdot \langle N_{x,w,d}, N_{x,w,d}\rangle_{N_{x,0}}$, and hence $\langle N_{x',w,d}, N_{x',w,d}\rangle_{N_{x',0}}^{1/d}= e(x',x)\cdot \langle N_{x,w,d}, N_{x,w,d}\rangle_{N_{x,0}}^{1/d}$. Hence $[C': C]h_{\diamondsuit, d,x}(\sH)= \sum_{x'} e(x',x)h_{\diamondsuit, d,x}(\sH)= \sum_{x'}h_{\diamondsuit,d,x'}(\sH')$ where $x'$ ranges over all points of $C'$ lying over $x$. Hence $[C':C]h_{\diamondsuit,w,d}(\sH)=h_{\diamondsuit,w,d}(\sH')$. \end{pf}

\begin{sbrem}\label{reasons}  (1) The above proof of \ref{CC'} shows that to have the  formulas in the third line, it is important to take the $d$-th root in  the definition of $h_{\diamondsuit,w,d,x}$ for $d\geq 2$ as in \ref{Na}.

(2) If we do not assume the unipotence of the local monodromy in \ref{Hodgeside} and assume only the quasi-unipotence, we still have the definition of these height functions $h_r(\sH)= \text{deg}(\gr^r\sH_{\sO})$ etc., by using the canonical extension $\sH_{\sO}$ of Deligne, but the above formulas in \ref{CC'} do not hold. 

\end{sbrem}

The following \ref{Hprop} was proved by Griffiths (\cite{Gri}) and Peters (\cite{Pe}) in the case $\sH$ is pure (the method of their proof will be sketched in \S\ref{2.3b}). 

\begin{sbprop}\label{Hprop}
Assume $\sH$ is pure (resp. the pure graded quotients of $\sH$ are polarized). 

\medskip

(1) $h_{\star}(\sH)\geq 0$ (resp. $h_{\star\diamondsuit}(\sH)\geq 0$). 

\medskip

(2) Assume $h_{\star}(\sH)=0$ (resp. $h_{\star\diamondsuit}(\sH)=0$) and assume that $\sH_\Q$ has a $\Z$-structure. Then over some finite \'etale covering of $C$, the pullback of $\sH$ to $C'$ becomes a constant pure (resp. mixed) Hodge structure. 
\end{sbprop}

In the mixed case, 
(1) follows from the pure case. To reduce (2) to the pure case, it is sufficient to prove the following proposition.

\begin{sbprop}\label{Hprop2}
  If $\gr^W_w\sH$ ($w\in \Z$) are constant and $h_{\diamondsuit}(\sH)=0$, then $\sH$ is constant. 

\end{sbprop}

\begin{pf}  By the induction on the length of the weight filtration and by $\Ext^1(\gr^W_w,  W_{w-1})= \Ext^1(\Q, (\gr^W_w)^*\otimes W_{w-1})$,  we are reduced to the following situation. $\sH=W_0\sH$, $\gr^W_0\sH=\Q$,   and $W_{-1}\sH$ is a constant variation associated to a mixed Hodge structure $H'$. By $H_{\diamondsuit, 0,d}(\sH)=0$ for all $d\geq 2$, $N_x=0$ for all $x$ and hence $\sH$ has no degeneration (we can take $U=C$). Let $\text{Ext}^1_C(\Q, H')$ be the abelian group of extensions in the category of all $\sH$ as in \ref{Hodgeside} with $U=C$. Then our $\sH$ gives an element $\alpha$ of $\text{Ext}^1_C(\Q, H')$. By similar arguments as in \ref{d=1bc}, we have 
$$\text{Ext}^1_C(\Q, H')\subset  H^0(C, H'_\Q\bs (\sO_C\otimes_\C (H'_\C/F^0))).$$
Consider the long exact sequence of cohomology groups associated to the exact sequence 
$$0 \to H'_\Q \to \sO_C \otimes_\C (H'_\C/F^0)\to H'_\Q\bs (\sO_C \otimes_\C (H'_\C/F^0))\to 0$$
of sheaves on $C$. We have $$H^0(C, \sO_C \otimes_\C (H'_\C/F^0))=H'_\C/F^0.$$
Let $\beta\in H^1(C, H'_\Q)= \Hom(H_1(C, \Q), H'_\Q)$ be the image of $\alpha$ under the connecting homomorphism. 

\medskip
{\bf Claim.} $\beta$ is a homomorphism $H_1(C) \to H'$ of mixed Hodge structures.

\medskip
To prove this, it is sufficient to show that the image of $H^0(C, \Omega^1_C) = F^0(\C\otimes H_1(C, \Q))$ in $H'_\C$ is contained in $F^0H'_\C$. But this follows from the fact that the image of $\beta$ in $H^1(C, \sO_C\otimes H'_\C/F^0)= H^1(C,\sO_C)\otimes H'_\C/F^0= \Hom_\C(H^0(C, \Omega^1_C), H'_\C/F^0)$ is zero because $\beta$ is in the image of the connecting homomorphism.

Since $H_1(C)$ is pure of weight $-1$, we have $$\Hom_{\text{MHS}}(H_1(C), H')\overset{\subset}\to \Hom_{\text{MHS}}(H_1(C), \gr^W_{-1}H')$$ $$\overset{\subset}\to \Hom_{\Q}(H_1(C,\Q), \gr^W_{-1}H'_\Q)= H^1(C, \gr^W_{-1}H'_\Q).$$
Since $H_{\diamondsuit,0,1}(\sH)=0$, the image of $\beta$ in $H^1(C, \gr^W_{-1}H'_\Q)$ is $0$ by the last remark in \ref{d=1a}. Hence $\beta=0$ and hence $\alpha$ comes from $H'_\C/F^0$. Hence the extension class of $\sH$ comes from $H'_\Q\bs H'_\C/F^0$, that is, $\sH$ is constant. 
\end{pf}

\subsection{The heights $H_{\star\diamondsuit}(M)$ of mixed motives}

 \begin{sbpara}

 $H_{\star\diamondsuit}(M)$ for $M$ in Case (mixed-pol) is defined assuming some conjectures which we can check to be true in many cases. It is defined as 
 $$H_{\star\diamondsuit}(M)=(\prod_{w\in \Z} H_{\star}(\gr^W_wH))\cdot H_{\diamondsuit}(M)$$
 where $H_{\diamondsuit}(M)$ tells, roughly speaking, how $M$ is far from being the direct sum of $\gr^W_wM$ ($w\in \Z$). The definition of $H_{\diamondsuit}(M)$ was given in the short summary \cite{Ka3}. We review it in this \S1.7. It has the shape 
$$H_{\diamondsuit}(M)= \prod_{w\in \Z, d\geq 1} H_{\diamondsuit, w,d}(\sH). $$
$H_{\diamondsuit, w, 1}(M)$ is defined by the theory of the height pairing of Beilinson \cite{Be2} , Bloch \cite{Bl} and Gillet-Soul\'e \cite{GS2} using the extension class of the exact sequence 
$0\to  \gr^W_{w-1}M\to W_wM/W_{w-2}M\to \gr^W_wM\to 0$.
For $d\geq 2$, $$H_{\diamondsuit,w,d}(M)= \prod_v  H_{\diamondsuit,w,d,v}(M),$$
where $v$ ranges over all places of $F$. For $v$ archimedean, $H_{\diamondsuit,w,d,v}(M)$ is closely related to Beilinson regulator in \cite{Be} as is seen in \S2.1. For $v$ non-archimedean, $H_{\diamondsuit,w,d,v}(M)$ is defined by using the monodromy operator  at $v$ (roughly speaking, it presents the ''archimedean size'' of the monodromy operator).
   
   For the definitions of $H_{\diamondsuit}$(M), $H_{\diamondsuit, w,d}(M)$ ($d\geq 1$), $H_{\diamondsuit, w,d,v}(M)$ ($d\geq 2$), we do not need that $M$ is with $\Z$-coefficients.

\end{sbpara}

\begin{sbpara}\label{Mdgeq2} Let $M$ be a motive with $\Q$-coefficients over a number field $F$ whose pure graded quotients $\gr^W_wM$ ($w\in \Z$) are polarized. 
Let $w\in \Z, d\geq 2$ and let $v$ be a non-archimedean place of $F$. We define $H_{\diamondsuit,w, d,v}(M)$ under the assumptions (i) and (ii) below. The Hodge analogues of these (i) and (ii) are \ref{Nd} (i) and \ref{Na}, respectively, which are true. We expect that the present (i) and (ii) are also true. 
Assuming (i), (ii), we will define $$H_{\diamondsuit,w, d, v}(M):= \sharp(\F_v)^{l(v)}\quad\text{with}\;\;  l(v)=\langle N_{v,w,d}, N_{v,w,d}\rangle^{1/d}_{N_{v,0}}\in \R_{\geq 0}$$ where $\sharp(\F_v)$ denotes the order of the residue field $\F_v$ of $v$ and $\langle N_{v,w,d}, N_{v,w,d}\rangle_{N_{v,0}}\in \Q_{\geq 0}$ is defined in the same way as $\langle N_{x,w,d}, N_{x,w,d}\rangle_{N_{x,0}}\in \Q_{\geq 0}$ in \ref{Na} but this time assuming (i) and (ii). 

Let $p=\text{char}(\F_v)$. 

Let $\ell$ be a prime number. We define $\langle N_{v,w,d}, N_{v,w,d}\rangle_{N_{v,0}}\in \Q_{\ell}$ for each $\ell$ by an $\ell$-adic method assuming (i) below.

Let $A=\Q_{\ell}$ if $\ell\neq p$, and let $A= F_{v,0,\text{ur}}$ if $\ell=p$. Here $F_{v,0}$ is (as in \ref{st}) the field of fractions of $W(\F_v)$.  Define a finite dimensional $A$-vector space $\Upsilon$ as follows. 
If $\ell\neq p$, let $\Upsilon=M_{et, \Q_\ell}$. 
If $\ell=p$, let $\Upsilon=D_{\pst}(F_v, M_{et,\Q_p})$.

We have the monodromy operators $N'_v: \Upsilon\to \Upsilon$ and 
$N_v: \Upsilon \to \Upsilon(-1)$ defined as follows. 
Here in the case $\ell=p$, $\Upsilon(-1)= \Upsilon$ as an $A$-vector space, but the Frobenius operator $\varphi$ on $\Upsilon(-1)$ is $p$ times $\varphi$ on $\Upsilon$. 

Assume first $\ell\neq p$. let $I_v$ be the inertia subgroup of $\Gal(\bar F_v/F_v)$ and let $a: I_v\to\Z_{\ell}(1)$ be the surjective homomorphism defined as $\sig\mapsto (\nu_n)_n$ where $\nu_n$ is an $\ell^n$-th root of $1$ defined by $\sig(\varpi^{1/\ell^n})= \nu_n  \varpi^{1/\ell^n}$ with $\varpi$ a prime element of $F_v$. Let $\sig$ be an element of $I_v$ such that $a(\sig)$ is a generator of $\Z_{\ell}(1)$, and let $N_v'=\log(\sig):= \frac{1}{m}\log(\sig^m): \Upsilon\to \Upsilon$ where $m\geq  1$ is an integer such that the action of $\sig^m$ on $\Upsilon$ is unipotent. This $N'_v$ depends on $\sig$ but is independent of the choice of $m$. Let $$N_v:= a(\sig)^{-1}N'_v= a(\sig)^{-1}\log(\sig): \Upsilon\to \Upsilon(-1).$$ Then $N_v$ is independent of the choice of $\sig$. Next in the case $\ell=p$, let $N_v': \Upsilon\to \Upsilon$ be  the monodromy operator and let $N_v=N'_v : \Upsilon\to \Upsilon(-1)$. 

We have $N'_vW_w\Upsilon\subset W_w\Upsilon$ for any $\ell$. 

We assume for any $\ell$
\medskip

(i) The weight monodromy conjecture is true for $\gr^W_w\Upsilon$ for any $w\in \Z$.

\medskip

 This means the following (i') is true for each $w\in \Z$. 
 
 \medskip
  Let $\sW$ be the relative monodromy filtration on $\gr^W_w\Upsilon$ for $N'_v:\gr^W_w\Upsilon\to \gr^W_w\Upsilon$.   Let $\bar \F_v$ be the residue field of $F_{v, \text{ur}}$, which is an algebraic closure of $\F_v$. 
 
Assume $\ell\neq p$. Let $\sig$ be an element of $\Gal(\bar F_v/F_v)$ whose image in $\Gal(\bar \F_v/\F_v)$ is $\bar \F_v\to \bar \F_v\;;\;x\mapsto x^{1/\sharp(\F_v)}$. 
 
 \medskip
 
 (i') for the case $\ell\neq p$: For any $m\in\Z$, the eigenvalues of the action of $\sig$ on $\gr^{\sW}_m\gr^W_w\Upsilon$ are algebraic numbers whose all conjugates over $\Q$ have complex absolute value  $\sharp(\F_v)^{m/2}$.
 
 \medskip
 
 Assume $\ell=p$. Take a finite extension ${}'F_v\subset \bar F_v$ of $F_v$ such that the action of the inertia subgroup of $\Gal({}'\bar F_v/{}'F_v)$ on $\Upsilon$ is trivial.  Let ${}'\Upsilon=D_{\st}({}'F_v, M_{et,\Q_p})$ be the fixed part of $\Gal({}'\bar F_v/{}'F_v)$ in $\Upsilon$. Let ${}'\F_v$ be the residue field of ${}'F_v$ and let ${}'A:={}'F_{v,0}={}'F_v\cap F_{v,0,\text{ur}}$ be the field of fractions of $W({}'\F_v)$. Then $\Upsilon= A\otimes_{{}'A} {}'\Upsilon$. Let $d=[{}'\F_v:\F_p]$. 
 
 \medskip
 
 (i') for the case $\ell=p$: For any $m\in \Z$, the eigenvalues of the ${}'A$-linear map  $\varphi^d$ on $\gr^{\sW}_m\gr^W_w({}'\Upsilon)$ are algebraic numbers whose all conjugates over $\Q$ have complex absolute value  $p^{dm/2}$.
 
 \medskip

The weight-monodromy conjecture was proved under certain (not so strong) assumption by Scholze in \cite{S3}.

Assume the weight-monodromy conjecture (i).
Then (1) and (2) in \ref{Assume12} are satisfied by $(\Upsilon, W, N'_v)$. In fact, if $\ell\neq p$ (resp. $\ell=p$), we have the relative monodromy filtration $\sW$ of $N'_v$ with respect to $W$ and  $\Upsilon_m$ in (2) of \ref{Assume12} as follows: $\Upsilon_m$ is the part of $\Upsilon$ on which all eigenvalues of $\sig$ (resp. $\varphi^d$) in (i') have complex absolute values $\sharp(\F_v)^{m/2}$ (resp. $p^{dm/2}$), and $\sW_m:=\oplus_{i\leq m} \Upsilon_i$. 
Hence $N_{v,d}\in (\gr^{\sW}_{-2}\gr^W_{-d}\Hom_A(\Upsilon, \Upsilon))_{\pri}(-1)$ is defined. Write
$$N_{v,d}=\sum_{w\in\Z} N_{v,w,d}$$
where $N_{v,w,d}$ is the $(\gr^{\sW}_{-2}\Hom(\gr^W_w\Upsilon, \gr^W_{w-d}\Upsilon)_{\pri})(-1)$-component of $N_{v,d}$. 

Let $M':= (\gr^W_wM)^*\otimes \gr^W_{w-d}M$. Then $M'$ is a pure motive of weight $-d$, and the polarizations of $\gr^W_wM$ and $\gr^W_{w-d}M$ induce a polarization  $M'\otimes M' \to \Q(d)$. It induces a non-degenerate $A$-bilinear form
$$\langle\;,\;\rangle: (\gr^{\sW}_{2-2d}\Hom(\gr^W_w\Upsilon, \gr^W_{w-d}\Upsilon))(1-d)\times (\gr^{\sW}_{-2}\Hom(\gr^W_w\Upsilon, \gr^W_{w-d}\Upsilon))(-1)\to A.$$ 
This induces a non-degenerate symmetric $A$-bilinear form
$$\langle\; , \;\rangle_{N_{v,0}}: (\gr^{\sW}_{-2}\Hom(\gr^W_w\Upsilon, \gr^W_{w-d}\Upsilon)_{\pri})(-1)\times (\gr^{\sW}_{-2}\Hom(\gr^W_w\Upsilon, \gr^W_{w-d}\Upsilon)_{\pri})(-1)\to A$$ 
$$(a, b)\mapsto \langle \Ad(N_{v,0})^{d-2}(a), b\rangle$$
where $N_{v,0}=\gr^W(N_v)$. 
Thus we obtain  $\langle N_{v,w,d}, N_{v,w,d}\rangle_{N_{v,0}}\in A$. If $\ell=p$,  $\langle N_{v,w,d}, N_{v,w,d}\rangle_{N_{v,0}}$ belongs to the part in $A$ fixed by $\varphi$, which is $\Q_p$.

Following the analogy with \ref{Na}, we conjecture that the following (ii) is true. 
\medskip

(ii)  $\langle N_{v,w,d}, N_{v,w,d}\rangle_{N_{v,0}}\in \Q_{\geq 0}$ and is independent of $\ell$. 
It is zero if and only if $N_{v,w,d}=0$.

\end{sbpara}

\begin{sbpara} For $w\in \Z, d \geq 2$ and for $v$ non-archimedean, assuming (i) and (ii), we define
$$H_{\diamondsuit, w,d, v}(M)= \sharp(\F_v)^{l(v)}\geq 1\quad \text{where}\;\; l(v)= (\langle N_{v,w,d}, N_{v,w, d}\rangle_{N_{v,0}})^{1/d}.$$

The reason why we take the $d$-th root here is to have the result for a base change $F'/F$ in \ref{F':F} (2) below, the same as the case of Hodge analogue  \ref{reasons} (2). 
\end{sbpara}

\begin{sbpara}\label{2.4.c} The independence of $\ell$ and the positivity  in (ii) in \ref{Mdgeq2} is explained by some philosophy on motives as in \ref{2.4.c2} below. We give a preparation for motives over an algebraic closure $\bar k$ of a finite field $k$ of characteristic $p$. 

We define the category of motives with $\Q$-coefficients over $\bar k$ to be the category of Grothendieck motives with $\Q$-coefficients over $\bar k$ modulo homological equivalence, which he defined using projective smooth schemes over $k$. Here, modulo homological equivalence means that a morphism is regarded as $0$ if it induces zero homomorphisms on the $\ell$-adic \'etale realizations for all $\ell\neq p$ and on the crystalline realization. (It is conjectured that this homological equivalence is the same as the rational equivalence, and also the same as the numerical equivalence.)  

Note that by Weil conjecture proved by Deligne, for any motive $M$ with $\Q$-coefficients over $\bar k$, we have a unique direct sum decomposition $M=\oplus_w M_w$ where $M_w$ is of weight $w$.  

There is a notion of polarization, but we do not review it. 

\end{sbpara}

\begin{sbpara}\label{2.4.c2} Let the notation be as in \ref{Mdgeq2}. A philosophy on motives is that the following holds:

Assume that $M$ is pure of weight $w$. Then there is a family of pure motives $(\frak P_m)_m$ of weight $m$ 
with $\Q$-coefficients over the algebraic closure $\bar \F_v$ of the residue field $\F_v$ of $v$, and a direct summand $\frak P_{m,\pri}$ of $\frak P_m$ for each $m\geq w$ having the following properties (1)--(3):

(1) For a prime number $\ell$, if $\ell\neq p$ (resp. $\ell=p$), $\gr^{\sW}_m \Upsilon$ is identified with the $\ell$-adic \'etale realization $\frak P_{m,et,\Q_\ell}$ (resp.  the crystalline realization $\frak P_{m,\crys}$),  $(\gr^{\sW}_m \Upsilon)_{\pri}$ is identified with $\frak P_{m,\pri,et,\Q_\ell}$ (resp.  $\frak P_{m,\pri, \crys}$), and $N_v: \gr^{\sW}_m \Upsilon\to (\gr^{\sW}_{m-2}\Upsilon)(-1)$ comes from a morphism $N_v: \frak P_m \to  \frak P_{m-2}(-1)$ of motives which is independent of $\ell$. 
This is similar to  \ref{Nb} in Hodge theory. 

(2) For $m\in \Z$, the homomorphism $\gr^{\sW}_{w-m}\Upsilon \otimes \gr^{\sW}_{w+m} \Upsilon \to A(-w)$ induced by the polarization $M \otimes M \to \Q(-w)$ of $M$ comes from a morphism  $\frak P_{w-m}\otimes \frak P_{w+m}\to \Q(-w)$ of motives which is independent of $\ell$. For $m\geq 0$, this morphism and the morphism $N_v^m: \frak P_{w+m} \to \frak P_{w-m}(-m)$  induce a polarization $\frak P_{w+m, \pri} \otimes \frak P_{w+m,\pri}
\to \Q(-w-m)$ of $\frak P_{w+m,\pri}$. This is similar to \ref{Nc} in Hodge theory.

(3) In particular, if $w=-d$ for an integer $d\geq 2$ and if we take $m=1$, we have the polarization
$\langle \;,\;\rangle_{N_v} : \frak P_{-2,\pri}(-1)\otimes \frak P_{-2,\pri}(-1) \to \Q$. If we put $P:=\Hom(\Q, \frak P_{-2,\pri}(-1))$ (the set of morphisms of motives), the  pairing $P\times P \to \Q$ induced by the polarization is positive definite.

Apply this by taking the pure motive $M'=(\gr^W_wM)^* \otimes \gr^W_{w-d}M$ of weight $-d$ in \ref{Mdgeq2} as $M$ above. We expect that
$N_{v, w,d}\in P$. This explains that $\langle N_{v,w,d}, N_{v, w,d}\rangle_{N_{v,0}}$ should belong to $\Q_{\geq 0}$ and be independent of $\ell$, and that it 
should be $>0$ if $N_{v,w,d}\neq 0$.

\end{sbpara}

\begin{sbpara}\label{deltareview} Now we consider an Archimedean place $v$. We define $H_{\diamondsuit, w,d,v}(M)$ for $d\geq 2$ (not using any conjecture). 
 We give first Hodge theoretic preparations.

Let $H$ be a mixed $\R$-Hodge structure. We review a homomorphism $\delta:\gr^WH_\R \to \gr^WH_\R$. 

As is explained in \cite{CKS} 2.20,  there is a unique pair $(s, \delta)$ where $s$ is a splitting $\gr^WH_\R\cong H_\R$ of $W$ and $\delta$  is a nilpotent $\R$-linear map $\gr^WH_\R\to \gr^WH_\R$ such that the original Hodge filtration $F$ of $H$ is expressed as $F=s(\exp(i\delta)\gr^W(F))$ and such that the Hodge $(p,q)$-component $\delta_{p,q}$ of $\delta$ with respect to $\gr^W(F)$ is zero unless $p<0$ and $q<0$.

\end{sbpara}

\begin{sbpara} Let $v$ be an archimedean place of $F$. Let $d\geq 2$.

Take a homomorphism $a: F\to \C$ which induces $v$. 
Let $\delta_v: \gr^WM_{a,B,\R}\to \gr^WM_{a,B,\R}$ (here $M_{a,B,\R}=\R\otimes_{\Q} M_{a,B}$)  be the $\delta$ (\ref{deltareview}) of the mixed $\R$-Hodge structure $\R\otimes_\Q M_{a,H}$. 
 Let $H'= (\gr^W_wM_{a,H})^*\otimes \gr^W_{w-d}M_{a,H}\subset \gr^W_{-d}(M_{a,H}^*\otimes M_{a,H})$ which is a polarized Hodge structure of 
weight $-d$. Let $\delta_{v,w,d}\in H'_\R$ 
be the $\Hom(\gr^W_wM_{a,B,\R}, \gr^W_{w-d}M_{a,B,\R})$-component of $\delta_v$. The Hodge metric $(\;,\;)$ of $H'$ defines $(\delta_{v,w,d}, \delta_{v,w,d})\in \R_{\geq 0}$.
Let 
$$H_{\diamondsuit, w,d,v}(M)= \exp(2(\delta_{v,w,d}, \delta_{v,w,d})^{1/d})\quad \text{if $v$ is a complex place},$$ 
$$H_{\diamondsuit, w,d,v}(M)=\exp((\delta_{v,w,d},\delta_{v,w,d})^{1/d})\quad \text{if $v$ is a real place}. $$
 The reason why we take here the $d$-th root is that we are following the analogy between $N_v$ for $v$ non-archimedean and $\delta_v$ in Hodge theory. 
 
  \end{sbpara}
 
\begin{sbpara}\label{1.7ex1} {\bf Example.} 
Let $a\in F^\times$. We have the mixed motive $M$ with $\Z$-coefficients over $F$ associated to $a$ having an exact sequence $0\to \Z(1)\to M \to \Z\to 0$. In this case, (i) and (ii) in \ref{Mdgeq2} are satisfied at all non-archimedean places of $F$. We have 
$$H_{\diamondsuit, 0, 2,v}(M)= \max(|a|_v, |a|_v^{-1})$$ for all places $v$ of $F$, and $H_{\diamondsuit}(M)=H_{\diamondsuit, 0, 2}(M)= \prod_v H_{\diamondsuit,0,2,v}(M)$. 
If $F=\Q$ and if $a=b/c$ with $b,c\in \Z$ such that $(b,c)=1$, we have $H_{\diamondsuit}(M)= \max(|b|^2, |c|^2)$. 

\medskip

In fact, $N_{v,0,2}= \text{ord}_v(a)$ for $v$ non-archimedean, 
$\delta_{v,0,2} = -i \log(|a|_v)\in \R(1)= \R \cdot 2\pi i$ if $v$ is real, and 
$\delta_{v,0,2} = -2^{-1}i \log(|a|_v)\in \R(1)= \R \cdot 2\pi i$ if $v$ is complex.
The above result on $H_{\diamondsuit,0,2,v}(M)$ follows from this. 

\end{sbpara}

\begin{sbpara}\label{BeBl} The definition of  $H_{\diamondsuit,w,1}$ in \ref{d=1} below is a simple application of the theory of  height parings of Beilinson, Bloch and Gillet-Soul\'e. 

 Let $M$ be a pure motive of weight $-1$ with $\Q$-coefficients over a number field $F$. 
 Let $$E(M):=\text{Ext}^1(\Q, M)$$ be the extension group in the category of mixed motives with $\Q$-coefficients over $F$. If $M=H^{2r-1}(X)(r)$ for a projective smooth scheme $X$ over $F$ with $r\in \Z$, we expect that 
  $$E(M)\cong (CH^r(X)\otimes\Q)_0$$
  where $(\;)_0$ denotes the part  homologically equivalent to $0$.

   \medskip

In  Beilinson \cite{Be2}, Bloch \cite{Bl} and Gillet-Soul\'e \cite{GS2},  the height paring 
$$\langle\;,\;\rangle:  E(M) \times E(M^*(1))\to \R_{>0}.$$
is defined assuming some conjectures. 
We describe one method following  Scholl \cite{Sl}. 

Let $a\in E(M)$, $b\in E(M^*(1))$. 
Then $a$ corresponds to a mixed motive $S_a$ with $\Q$-coefficients over 
$F$ such that $W_0S_a=S_a$, $W_{-2}S_a=0$, $\gr^W_0S_a=\Q$, $\gr^W_{-1}S_a=M$, 
and since $\text{Ext}^1(\Q, M^*(1))\cong \text{Ext}^1(M, \Q(1))$ (by duality), 
$b$ corresponds to a mixed motive $S_b$ with $\Q$-coefficients over $F$ such that 
$W_{-1}S_b=S_b$, $W_{-3}S_b=0$, $\gr^W_{-1}S_b=M$, $\gr^W_{-2}S_b=\Q(1)$.

We assume the following (i)$=$ (i-1)$+$(i-2)  (a philosophy of motives explained below tells that this (i) should be always true.)
 
 \medskip
 
 (i-1) There is a mixed motive $S$ with $\Q$-coefficients over $F$ such that 
 $$W_0S=S, \quad W_{-3}S=0, \quad S/W_{-2}S=S_a, \quad W_{-1}S=S_b.$$
 
 \medskip
 
 The group $F^\times\otimes \Q$ acts on the set of isomorphism classes of $S$ as in (i-1) as follows. 
 For $S$ as in (i-1) and for $c\in F^\times\otimes \Q$, $c$ sends the class of $S$ to the class of $S'$ where $S'$ is obtained from $S$ and $c$ as follows:  Let $E_c$ be the mixed motive with $\Q$-coefficients over $F$ corresponding to $c$ with an exact sequence $0\to \Q(1)\to E_c\to \Q\to 0$.  Take the fiber product $J$ of $S\oplus E_c\to \Q^2\leftarrow \Q$ and then  take the pushout $S'$ of $\Q(1)\leftarrow \Q(1)^2\to J$. Here the morphism $\Q\to \Q^2$ is $x \mapsto (x,x)$ and the morphism $\Q(1)^2\to \Q(1)$ is $(x,y)\mapsto x-y$. 
 \medskip

 (i-2) The action of $F^\times\otimes \Q$ on the set of  isomorphism classes of $S$ as in (i-1) is transitive. 
   
\medskip

We also assume

\medskip
(ii) $M$ satisfies the weight-monodromy conjecture at each non-archimedean place of $F$, 

\medskip

By (ii), for each non-archimedean place $v$ of $F$ and for each prime number $\ell$, the representation $M_{et,\Q_{\ell}}$ of $\Gal(\bar F_v/F_v)$ gives 
$$N_{v,0,2}\in \Q_{\ell}.$$
In fact, for a prime number 
$\ell \neq p=\text{char}(\F_p)$, the $\ell$-adic method gives $N_{v,0,2}\in \Hom_{\Q_{\ell}}(\gr^W_0M_{et, \Q_{\ell}}, \gr^W_{-2}M_{et, \Q_{\ell}})(-1)= \Hom_{\Q_{\ell}}(\Q_{\ell}, \Q_{\ell}(1))(-1)=\Q_{\ell}$, and the $p$-adic method similarly gives $N_{v,0,2}$ in the fixed part of $D_{\pst}(\Q_p(1))(-1)$ by Frobenius which is $\Q_p$.
We assume that

\medskip

(iii) This $N_{v,0,2}$ belongs to $\Q$ and is independent of $\ell$. 

\medskip

On the other hand, for each archimedean place $v$ of $F$, we have $$\delta_{v,0,2}\; \text{of}\; S\; \in \R(1)=\R\cdot 2\pi i.$$

Define $t(v)= \sharp(\F_v)^{N_{v,0,2}}$ if $v$ is non-archimedean, $t(v)=\exp(i \delta_{v,0,2})$ if $v$ is a real place, and $t(v)=\exp(2i \delta_{v,0,2})$ if $v$ is a complex place. Define
$$\langle a, b\rangle:= \prod_v t(v)$$
where $v$ ranges over all places of $F$. If $S$ is replaced by its modification $S'$ obtained by using the action of $c\in F^\times$  in  (i-2),  
$t(v)$  is replaced by  $t(v)|c|_v$ for any place $v$ of $F$. By product formula, $\langle a, b\rangle$ is independent of the choice of $S$ as in (i-1).

If we have a polarization $p:M\to M^*(1)$ of $M$, we have a pairing 

\medskip

(*) $\quad \langle\;,\;\rangle: E(M) \times  E(M)\to \R_{>0}\quad (a, b) \mapsto \langle a, p(b)\rangle.$

\medskip
\noindent
It is a symmetric pairing.
It is expected that the following (iv) is true.

\medskip

(iv) The logarithm of the above paring  $\log\langle\;,\;\rangle: E(M)\times E(M) \to \R$ is  positive definite.

\medskip

 The philosophical reasons of (i-1) and (i-2) are as follows. 
We  have an exact sequence $$0\to \text{Ext}^1(\Q, \Q(1)) \to \text{Ext}^1(\Q,S_b)\to \text{Ext}^1(\Q, M)\to \text{Ext}^2(\Q, \Q(1))$$ associated to the exact sequence $0\to \Q(1) \to S_b\to M \to 0$ of mixed motives. A philosophy of motives tells $\text{Ext}^i(\Q, \Q(1))= H^{i-1}(\Spec(F), {\bf G}_m)$. Hence  we should have  $\text{Ext}^2(\Q, \Q(1))=\Q\otimes \text{Pic}(F)=0$. Hence the above exact sequence tells that there should be an element $s$ of $\text{Ext}^1(\Q,S_2)$ whose image in $\text{Ext}^1(\Q, M)$ coincides with the class of $S_a$ and this $s$, which gives $S$ in (i-1),  is determined modulo the image of $\text{Ext}^1(\Q, \Q(1))$ which should be identified with $F^\times \otimes \Q$.

\end{sbpara}

\begin{sbpara}\label{d=1} 
We define
 $H_{\diamondsuit,w,1}$. 

Let $M$ be a mixed motive with $\Q$-coefficients over a number field $F$ whose pure graded quotients $\gr^W_wM$ ($w\in \Z$) are polarized. 
Let $M'= (\gr^W_wM)^*\otimes \gr^W_{w-1}M\subset \gr^W_{-1}(M^*\otimes M)$. This is a polarized pure motive of weight $-1$. We assume 
that (i)--(iii) in \ref{BeBl} are true if we take this $M'$ as $M$ in \ref{BeBl}.

 \medskip
 Let $\alpha_w\in \text{Ext}^1(\gr^W_wM, \gr^W_{w-1}M)=\text{Ext}^1(\Q, M')$ be the extension class of
 $0\to \gr^W_{w-1}M\to W_wM/W_{w-2}M\to \gr^W_wM\to 0$.  
  We define $$H_{\diamondsuit,w,1}(M)= \langle \alpha_w, \alpha_w\rangle$$
 by using the pairing $\langle\;,\;\rangle$ in (*) of \ref{BeBl}. If we assume (iv) of \ref{BeBl}, $H_{\diamondsuit, w,1}(M)\geq 1$, and we have $H_{\diamondsuit, w,1}(M)=1$ if and only of $\alpha=0$.

\end{sbpara}

\begin{sbrem} 
In the case $M=H^{2r-1}(X)(r)$, under certain assumptions, the rationality (iii) in \ref{BeBl}  is proved by understanding  $N_{v,0,2}$ in \ref{BeBl}  as a certain intersection number on the special fiber at $v$ of a proper regular integral models over $O_F$. See \cite{Be2}, \cite{Bl}, \cite{GS2}.

\end{sbrem}

 \begin{sbprop}\label{F':F}  Let $M$ be in Case (mixed-pol). Assume (i) and (ii) in \ref{Mdgeq2} and (i)--(iii)  in \ref{d=1}.
 
 (1) For Tate twists and for the dual, we have $H_{\diamondsuit}(M(r))= H_{\diamondsuit}(M)$ and $H_{\diamondsuit}(M^*)=H_{\diamondsuit}(M)$, and similar formulas for $H_{\diamondsuit,w,1}$ and $H_{\diamondsuit, w,d,v}$ for $d\geq 2$ and for each place $v$ of $F$. 
 
 (2) Let $F'$ be a finite extension of $F$, and let $M'$ be the motive over $F$ induced by $M$.
 We have  $H_{\diamondsuit}(M')= H_{\diamondsuit}(M)^{[F':F]}$ and similar formulas for $H_{\diamondsuit,w,1}$ and 
 $H_{\diamondsuit, w,d,v}$ for $d\geq 2$ and for each place $v$ of $F$.
 \end{sbprop}

\begin{sbrem} (1) If we like to define $H_{\diamondsuit,w,d,v}$ for $d\geq 2$ and for $v$ non-archimedean  only under the assumption (i) in \ref{Mdgeq2} not assuming the rationality (ii) in \ref{Mdgeq2}, we can do as follows.  Choose a prime number $\ell$ and a field homomorphism $\Q_{\ell}\to \C$ and define $l(v)$ as the absolute value of the the image of $\langle N_{v,w,d}, N_{v,w,d}\rangle_{N_{v,0}}\in \Q_{\ell}$ (which is defined $\ell$-adically) in $\C$. Then define $H_{\diamondsuit,w,d,v}(M)$ to be $\sharp(\F_v)^{l(v)}$.

(2) If we like to define $H_{\diamondsuit, w,1}$ assuming only (i) and (ii) in \ref{BeBl} not assuming the rationality (iii), we can do as follows. For each non-archimedean place $v$ of $F$, choose a prime number $\ell$ and a field homomorphism $\Q_{\ell}\to \C$ and let $\text{Re}(N_{v,w,2})$ be the real part of the image of $N_{v,w,2}\in \Q_{\ell}$ (which is defined $\ell$-adically) in $\C$. Then define $t(v)$ in \ref{BeBl} to be $\sharp(\F_v)^{\text{Re}(N_{v,w,2})}$.

With the definitions in these (1) and (2),  we still have the result in the above \ref{F':F}.
 
\end{sbrem}

\section{Speculations}

In this \S2, we  present our speculations on motives. 

 In \S2.1, we consider how the height function  $H_{\diamondsuit}$ of mixed motives with fixed pure graded quotients are related to the papers \cite{Be} of Beilinson and \cite{BK} of Bloch and the author on motives.  The main things of \S2.1 are presented in \ref{2.1thm}  and \ref{fin}.
 
 In \S2.2, we review period domains and their toroidal partial compactifications, and define the set $X(F)$ of motives. In \S\ref{s:2.2b}, we define  height functions $h_{\heartsuit}$ and $H_{\heartsuit,S}$, and show that the height functions $h_{\spadesuit}$ and $h_{\heartsuit}$ are equal to the degrees of the pullbacks of $''K+D''$ and $''D''$ of  a toroidal partial compactification  $\bar X(\C)$ 
 of a period domain $X(\C)$, respectively,  via the period map. This \S\ref{s:2.2b} gives us philosophies how to compare (1)--(4) in Introduction. 
 
 Basing on these philosophies,  we present  our speculations on the relation between  curvature forms on Griffiths period domains and motives in \S\ref{2.3b}, on the motive version of Vojta conjectures in \S\ref{2.4b}, and on the motive version of Manin-Batyrev conjectures in \S\ref{2.5b}. 
 
\subsection{Mixed motives with fixed pure graded quotients}\label{s:2.1}

\begin{sbpara}\label{S2.1}
In this \S2.1, we discuss the relations between the following three subjects.

\medskip

(a) The finiteness of the number of mixed motives of bounded height with fixed pure graded quotients.

(b) The conjectural Mordell-Weil theorems for $K$-groups of algebraic varieties over number fields.

(c) Some conjectures in the works \cite{Be} of Beilinson and in \cite{BK} of Bloch and the author.

\medskip

Assuming conjectures which appeared in the definition of $H_{\diamondsuit}$ in \S1.7, we will have (c) $\Rightarrow$ (a) $\Rightarrow$ (b) (\ref{2.1thm}, \ref{fin}). 

The more precise meanings of (a), (b), (c) are explained in \ref{(a)}, \ref{(b)}, \ref{(c)} below, respectively. 

\end{sbpara} 

\begin{sbpara}\label{(a)}
We explain (a) in \ref{S2.1}.

\medskip

Let $F$ be a number field. Assume that we are given a polarized pure motive $M_w$ of weight $w$ with $\Z$-coefficients over $F$ for each $w\in \Z$ satisfying $M_w=0$ for almost all $w$.

We expect that the following (1) is true. 
\medskip

(1) If $C\in \R_{>0}$, there are only finitely many isomorphism classes of mixed motives $M$ with $\Z$-coefficients over $F$ such that $\gr^W_wM=M_w$ for all $w\in \Z$ and such that $H_{\diamondsuit}(M)\leq C$.

\medskip
This is regarded as a motivic version of \ref{Hprop2}  in Hodge theory.

\end{sbpara}

\begin{sbpara}\label{(b)} We explain (b) in \ref{S2.1}.

Let $d\geq 1$ be an integer. Let $M$ be a pure motive with $\Z$-coefficients over $F$ of weight $-d$, and let $\text{Ext}^1(\Z, M)$ be the group of extension classes of $0\to M\to \tilde M \to \Z\to 0$ where $\tilde M$ is a mixed motive
with $\Z$-coefficients over $F$. The following is expected in the case $M=H^m(X)(r)$ for a smooth projective scheme $X$ over $F$ and  for $m, r\in \Z$ such that $d=2r-m$. 

\medskip
(1) If $d\geq 2$, $$\Q\otimes \text{Ext}^1(\Z, M)\cong \gr^r_{\gamma} (\Q\otimes K_{d-1}(X))$$
where $\gamma$ is the Gamma-filtration. 

(2) If $d=1$, as in \ref{BeBl}, $\Q\otimes \text{Ext}^1(\Z, M) \cong (\Q\otimes_\Z CH^r(X))_0$. 

\medskip

We expect that the following (1) and (2) are true. 

\medskip

(1) If $d\neq 2$, $\text{Ext}^1(\Z, M)$ is finitely generated as an abelian group.

(2) If $d=2$, if we fix a finite set $S$ of places of $F$ which contains all archimedean places of $F$, the subgroup of $\text{Ext}^1(\Z, M)$ consisting of classes of mixed motives which are of good reduction outside $S$ is finitely generated as an abelian group. 

\medskip

This is related to the subject (a) because $\text{Ext}^1(\Z, M)$ is classifying mixed motives $\tilde M$ with $\Z$-coefficients over $F$ such that $\gr^W_w\tilde M$ is $\Z$ if $w=0$, $M$ if $w=-d$, and $0$ for other $w$. 

\end{sbpara}

\begin{sbpara}\label{(c)}

The subject (c) in \ref{S2.1} is about the conjectures 
in \cite{Be} and \cite{BK} concerning  $\Q\otimes_\Z \text{Ext}^1(\Z, M)$ or concerning  $\Q\otimes $ the $K$-group in \ref{(b)}, which we review in \ref{BBKEM} below.
We will consider how these conjectures are related to the subject (a). 
 These conjectures are not in the parts of \cite{Be} and \cite{BK} concerning zeta functions, but in  the more basic parts. The parts concerning zeta functions in the works  \cite{Be} \cite{BK} should be also related to height functions of motives as is explained in Remark \ref{BBK}.

\end{sbpara}
In the following, \ref{loc1}--\ref{loc5} are preparations on local fields, \ref{EM1}--\ref{BBKEM} are preparations on global fields, and then in \ref{Sabc}--\ref{fin}, we discuss the above main subjects. 
 
In \ref{loc1}--\ref{loc4},  we consider Galois cohomology of local fields. In \ref{loc5}, we consider the archimedean analogue. 

\begin{sbpara}\label{loc1} Let $K$ be a finite extension of $\Q_p$. 

Let $\ell$ be a prime number and let $V$ be a finite-dimensional $\Q_{\ell}$-vector space endowed with a continuous action of $\Gal(\bar K/K)$. Then we have $\Q_{\ell}$-subspaces $$H^1_e(K, V)\subset H^1_f(K, V)\subset H^1_g(K, V)$$
(\cite{BK}) of the finite-dimensional $\Q_{\ell}$-vector space $H^1(K, V)=H^1_{\text{cont}}(\Gal(\bar K/K), V)$. We review basic things on these subspaces.

If  $\ell\neq p$, $H^1_e(K, V):=0$, $H^1_g(K, V):= H^1(K, V)$, and $H^1_f(K, V)$ is defined as the kernel of $H^1(K, V)\to H^1(K_{\text{ur}}, V)$. If $\ell=p$, $H^1_e(K, V)$, $H^1_f(K, V)$, and $H^1_g(K, V)$ are defined as the kernels of the canonical maps from $H^1(K, V)$ to $H^1(K, B_{\text{crys}}^{\varphi=1}\otimes_{\Q_p} V)$, $H^1(K, B_{\crys}\otimes_{\Q_p} V$), and $H^1(K, B_{\dR}\otimes_{\Q_p} V)$, respectively. 

Let $a\in H^1(K, V)$ and let $0\to V\to E_a\to \Q_{\ell}\to 0$ be the corresponding exact sequence of continuous representations of $\Gal(\bar K/K)$ over $\Q_{\ell}$. If $\ell\neq p$ and $V$ is unramified, $a\in H^1_f(K,V)$ if and only if $E_a$ is unramified. If $\ell=p$ and $V$ is a de Rham representation of $\Gal(\bar K/K)$, $a\in H^1_g(K, V)$ if and only if $E_a$ is de Rham. If $\ell=p$ and $V$ is crystalline, $a\in H^1_f(K, V)$ if and only if $E_a$ is crystalline.

In the Tate  duality $H^1(K, V)\times H^1(K, V^*(1))\to \Q_{\ell}$, if $\ell\neq p$ (resp. $\ell=p$ and $V$ is de Rham), $H^1_e(K, V)$, $H^1_f(K, V)$, $H^1_g(K, V)$  coincide with  the annihilators of $H^1_g(K, V^*(1))$, $H^1_f(K, V^*(1))$, $H^1_e(K, V^*(1))$, respectively. 

\end{sbpara}

\begin{sbpara}\label{loc2} Let the notation be as in \ref{loc1}. In the case $\ell=p$, we assume that $V$ is a de Rham representation. For $a\in H^1(K, V)$ and for the corresponding exact sequence $0\to V\to E_a\to \Q_{\ell}\to 0$, we review the relation of the monodromy operator of $E_a$ and the image of $a$ in $H^1_g(K, V)/H^1_f(K, V)$. 

\medskip

(1) Assume first $\ell\neq p$ and that the action of the inertia subgroup $I$ of $\Gal(\bar K/K)$ on $V$ is unipotent.
Then we have  canonical isomorphisms $$H^1_g(K, V)/H^1_f(K, V)\overset{\cong}\to  H^1(K_{\text{ur}}, V)^{\sigma=1}\cong (V(-1)/NV)^{\sig=1},$$ where $\sig$ is the canonical generator of $\Gal(K_{\text{ur}}/K)$,  $N$ denotes the monodromy operator $V\to V(-1)$ of $V$, and $\sig=1$ denotes the fixed part of $\sig$. This isomorphism sends $a\in H^1(K, V)$ to the class of the monodromy operator $N_a: E_a\to E_a(-1)$. Here the class of $N_a$ is defined as the class of $N_a(e)$ where $e$ is any lifting of $1\in \Q_{\ell}$ to $E_a$. 

If $d\geq 1$ is an integer and if $V$ satisfies the weight-monodromy conjecture for the weight filtration $W$ of pure weight $-d$ (this means that $W$ is the increasing filtration defined by $W_{-d}V=V$ and $W_{-d-1}V=0$), then since $$((\gr^{\sW}_{-2}V)_{\pri}(-1))^{\sig=1}\cong (V(-1)/NV)^{\sig=1},$$ we have a canonical isomorphism
$$H^1_g(K, V)/H^1_f(K, V) \cong ((\gr^{\sW}_{-2}V)_{\pri}(-1))^{\sig=1}.$$
This isomorphism sends $a$ to $N_{a,d}$ (\S1.5). 

\medskip

(2) Assume next $\ell=p$ and that $V$ is semi-stable. Then we have a canonical isomorphism $$H^1_g(K, V)/H^1_f(K, V) \cong (D_{\st}(V)(-1)/ND_{\st}(V))^{\varphi=1}$$
which sends $a\in H^1_g(K, V)$ to the class of  the monodromy operator $N_a: D_{\st}(E_a)\to D_{\st}(E_a)(-1)$ (note that $E_a$ is semi-stable by the assumption $V$ is semi-stable and $a\in H^1_g(K, V)$). 

If $d\geq 1$ is an integer and if $D_{\st}(V)$ satisfies the weight-monodromy conjecture for the weight filtration $W$ of pure weight $-d$, since $$(\gr^{\sW}_{-2}D_{\st}V)_{\pri}(-1))^{\varphi=1}\cong (D_{\st}(V)(-1)/ND_{\st}(V))^{\varphi=1},$$ we have a canonical isomorphism
$$H^1_g(K, V)/H^1_f(K, V) \cong ((\gr^{\sW}_{-2}D_{\st}(V))_{\pri}(-1))^{\varphi=1}.$$
This isomorphism sends $a\in H^1_g(K, V)$ to $N_{a,d}$ (\S1.5). 

\medskip

(3) In general, if $\ell\neq p$ (resp. $\ell=p$), the action of $I$ on $V$ is quasi-unipotent  (resp. $V$ is potentially semi-stable). Hence by the above (1) (resp.  (2)) for a finite extension $K'$ of $K$ over which we have a unipotent action of the inertia group (resp. semi-stable representation), if $d\geq 1$ and if the weight-monodromy conjecture in (1) (resp. (2)) over $K'$ is true, 
we have a canonical injective homomorphism
$$\iota: H^1_g(K, V)/H_f(K, V)\to \gr^{\sW}_{-2}(V)_{\pri}(-1)$$ $$(\text{resp.}\; \iota: H^1_g(K,V)/H^1_f(K, V)\to \gr^{\sW}_{-2}D_{\pst}(V)_{\pri}(-1))$$
which sends $a\in H^1_g(K, V)$ to $N_{a,d}$.

\end{sbpara}

\begin{sbpara}\label{loc3} Let the notation be as in \ref{loc1} and let $T$ be a $\Gal(\bar K/K)$-stable $\Z_{\ell}$-lattice in $V$. Then $H^1(K, T)= H^1_{\text{cont}}(\Gal(\bar K/K), T)$ is a finitely generated $\Z_{\ell}$-module. By $H^1_e(K, T)$, $H^1_f(K, T)$, $H^1_g(K, T)$, we denote the inverse image of $H^1_e(K, V)$, $H^1_f(K, V)$, $H^1_g(K, V)$ in $H^1(K, T)$, respectively. 

\end{sbpara}

\begin{sblem}\label{loc4} Let the notation be as in \ref{loc1}. Let $T$ be a $\Gal(\bar K/K)$-stable $\Z_{\ell}$-lattice in $V$. 

(1) Assume $\ell\neq p$. Assume $V$ is unramified. Then the image of $H^1_g(K, T)/H^1_f(K, T) \to V(-1)^{\sig=1}$ is contained in $T(-1)$.

(2) Assume $\ell=p$. Assume $V$ is crystalline, $p$ is a prime element in $K$, and there is $r\in \Z$ such that $D^r_{\dR}(V)=D_{\dR}(V)$ and $D^{r+p-1}_{\dR}(V)=0$. Let $\sD_{\text{FL}}(T)\subset D_{\crys}(V)$ be the Fontaine-Laffaille module \cite{FL} corresponding to $T$. Then the image of $H^1_g(K, T)/H^1_f(K, T) \to D_{\crys}(V)(-1)^{\varphi=1}$ is contained in $\sD_{\text{FL}}(T)(-1)$.

\end{sblem}

\begin{pf} (1) is clear. We prove (2). The isomorphism $H^1_g(K, V)/H^1_f(K, V)\overset{\cong}\to D_{\crys}(V)(-1)^{\varphi=1}$ is the dual of $H^1_f(K, V^*(1))/H^1_e(K, V^*(1))\cong D_{\crys}(V^*(1))/(1-\varphi)D_{\crys}(V^*(1))$ 
which comes from a canonical isomorphism 
$H^1_f(K, T^*(1))\cong \sD_{\text{FL}}(T^*(1))/(1-\varphi)\sD^0_{\text{FL}}(T^*(1))$ (\cite{BK} Lemma 4.5). Hence by duality, the homomorphism $H^1_g(K, T) \to D_{\crys}(V)(-1)$ comes from  the composition 
$H^1_g(K, T) \to \Hom_{\Z_p}(H^1(K, T^*(1)), \Z_p)\to \newline \Hom_{\Z_p}(\sD_{\text{FL}}(T^*(1)), \Z_p)\cong \sD_{\text{FL}}(T)(-1).$
\end{pf}

\begin{sbpara}\label{loc5} We consider the corresponding archimedean theory. Let $d\geq 1$ and let $H$ be a pure $\R$-Hodge structure of weight $-d$. Then the group $\Ext^1_{\R \text{MHS}}(\R, H)$ of extensions in the category of mixed $\R$-Hodge structures is isomorphic to $H_\R\bs H_\C/F^0H_\C$. 

Let $H_\R^{\leq 1, \leq 1}$ be the intersection $H_\R\cap \oplus_{p\leq -1, q\leq -1,p+q=-d} \; H_\C^{p,q}$ where $H_\C^{p,q}$ is the component of $H_\C$ of Hodge type $(p,q)$. Then we have an isomorphism
$$H_{\R}^{\leq -1, \leq -1}\overset{\cong}\to H_\R\bs H_\C/F^0H_\C\;;\; \delta\mapsto i\delta.$$
Hence we have a canonical isomorphism
$$\delta: \Ext^1_{\R \text{MHS}}(\R, H)\overset{\cong}\to H_\R^{\leq -1, \leq-1}.$$
If $a\in \Ext^1_{\R \text{MHS}}(\R, H)$ and if $0\to H\to E_a\to \R \to 0$ is the corresponding exact sequence 
of mixed $\R$-Hodge structures, the last isomorphism sends $a$ to the class of $\delta$ of $E_a$. 

\end{sbpara}

\begin{sbpara}\label{EM1}  

In \ref{EM1}--\ref{2.1thm2}, let $F$ be a number field, let $d\geq 2$, and let $M$ be a pure motive with $\Q$-coefficients over $F$ of weight $-d$.

Let $S$ be a finite set of places of $F$ containing all archimedean places.  Let $H^1_{g,S}(F,M_{et})$ be the subgroup of $H^1_{\text{cont}}(\Gal(\bar F/F), M_{et})$ consisting of all elements $a$ 
whose images $a_{v,\ell}$ in $H^1(F_v, M_{et,\Q_{\ell}})$ for non-archimedean places $v$ of $F$ and  prime numbers $\ell$ satisfy the following (i) and (ii).

\medskip

(i) $a_{v,\ell}\in H^1_g(F_v, M_{\Q_{\ell}})$ for all $v$ and $\ell$.

(ii) $a_{v,\ell}\in H^1_f(F_v, M_{et,\Q_{\ell}})$ for all $\ell$ if $v\notin S$. 

\medskip
Let $$H^1_g(F, M_{et})= \cup_S \; H^1_{g,S}(F, M_{et})\subset H^1_{\text{cont}}(\Gal(\bar F/F), M_{et}).$$

We also introduce notation for archimedean theory. For an archimdean place $v$ of $F$, let $E_v^{\R}$ be the following $\R$-vector space. Let $a:F\to \C$ be a homomorphism which induces $v$. If $v$ is a complex place, let $E^{\R}_v=M_{a,B, \R}\cap M_{a,H,\C}^{\leq -1,\leq -1}$ where $M_{a,H,\C}^{\leq -1,\leq -1}$ is the part of $\C\otimes_\Q M_{a,B}$ of Hodge type $(\leq -1, \leq -1)$.
In the case $v$ is a real place, let $E^{\R}_v$ be $M_{a,B, \R}^{\sigma=-1}\cap M_{a,H,\C}^{\leq -1,\leq -1}$ where $M_{a,B, \R}^{\sigma=-1}$ denotes the part of $M_{a,B,\R}=\R\otimes_\Q M_{a,B}$ on which the complex conjugation $\sigma$ on $M_{a,B}$ acts by $-1$. 
\end{sbpara}

\begin{sbprop}\label{EM6}

 If $d\neq 2$, we have $H^1_g(F, M_{et})=H^1_{g,S}(F, M_{et})$ if  $S$ is large enough. 

\end{sbprop}
\begin{pf} If $S$ is large enough, by Weil conjecture proved by Deligne, \ref{loc2} shows that $H^1_g(F_v, M_{et,\Q_{\ell}})/H^1_f(F_v, M_{et,\Q_{\ell}})=0$ for any $v\notin S$ and any $\ell$. 
\end{pf}

\begin{sbpara}\label{Egiven}  Let the notation be as in \ref{EM1}. 

We assume that we are given a $\Q$-vector space $E$ endowed with a $\Q$-linear map
$$E\to (\oplus_{v|\infty} E_v^{\R}) \oplus H^1_g(F, M_{et}).$$

For a finite set $S$ of places of $F$ containing all archimedean places, let $E_S\subset E$ be the inverse image of $H^1_{g,S}(F, M_{et})$ under the map $E\to H^1_g(F, M_{et})$. By \ref{EM6}, if $d\neq 2$, we have $E=E_S$ for a sufficiently large $S$. 

We are interested in the following two situations.  

\medskip

(Situation I). $E:= \text{Ext}^1(\Q, M)$. 

\medskip

(Situation II). Let $X$ be a projective smooth scheme over $F$, fix $m,r\in \Z$ such that $d=2r-m$, let $M$ be the pure motive $H^m(X)(r)$ with $\Q$-coefficients over $F$ of weight $-d$, and let $E:= \gr^r_{\gamma} (\Q\otimes K_{d-1}(X))$. 

\medskip

We expect that we have $E= \text{Ext}^1(\Q, M)$ also in the situation II, but we like to discuss finiteness concerning the $K$-group without assuming it in the situation II.  

In the situation I (resp. the situation II), the map $E\to E_v^{\R}$  is obtained by taking the associated Hodge structures and using \ref{loc5} (resp. is the Beilinson regulator map). 
In the situation I 
(resp. situation II), the map $E\to H^1_g(F, M_{et})$ is obtained by the \'etale realization of extensions (resp. by the Chern class map). In the situation II, these maps should coincide with the maps in the situation I via the identification $E=\text{Ext}^1(\Q, M)$.  

\end{sbpara}

\begin{sbpara}

Assume now that $M$ is endowed with a structure of $\Z$-coefficients. Let  $E^{\Z}$ be the fiber product of
$$E\to H^1_{\text{cont}}(\Gal(\bar F/F), M_{et}) \leftarrow  H^1_{\text{cont}}(\Gal(\bar F/F), M_{et, \hat \Z}).$$

Since $G:=\Gal(\bar F/F)$ is compact, we have
$H^1(G, M_{et})= \varinjlim_n  H^1(G, \frac{1}{n}M_{et,\hat \Z})= \Q \otimes H^1(G, M_{et,\hat \Z})$. Hence we have 
$$\Q\otimes E^{\Z}= E.$$

In the situation I, we have $E^{\Z}=\text{Ext}^1(\Z, M)$ (\ref{(b)}).

For a finite set $S$ of places of $F$ which contains all archimedean places, let $E^{\Z}_S\subset E^{\Z}$ be the inverse image of $E_S$ under the map $E^{\Z}\to E$. 
By \ref{EM6}, if $d\neq 2$, we have $E^{\Z}_S=E^{\Z}$ for a sufficiently large $S$.

\end{sbpara}

\begin{sbprop}\label{WMW}

 For any $n\geq 1$, $E^{\Z}_S/nE^{\Z}_S$ is a finite group.  

\end{sbprop}

\begin{pf} This quotient group is embedded in  $\prod_{\ell |n} B_{\ell}/nB_{\ell}$ where $\ell$ ranges over all prime divisors of $n$, and 
$B_{\ell}:= H^1_{et}(\Spec(O_{F,S'}), M_{et,\Z_{\ell}})= H^1_{\text{cont}}(\Gal(F_{S'}/F), M_{et,\Z_{\ell}})$ for a finite set $S'$ of places of $F$ which contains $S$, all prime divisors of $n$ in $F$, and all non-archimedean places at which $M$ is of bad reduction. Here $O_{F,S'}$ denotes the ring of $S'$-integers in $F$, $H^1_{et}$ denotes the \'etale cohomology, and $F_{S'}$ denotes the largest Galois extension of $F$ in $\bar F$ in which all places of $F$ outside $S'$ are unramified. As is well known, $B_{\ell}$ is a finitely generated $\Z_{\ell}$-module. 
\end{pf}

\begin{sbpara}\label{BBKEM} Now assume that $M$ is endowed with a polarization. 
For each non-archimedean place $v$ of $F$, let $E_v$ be the image of $$E\to \prod_{\ell} \; H^1_g(F_v, M_{et, \Q_{\ell}})/H^1_f(F_v, M_{et,\Q_{\ell}})$$
where $\ell$ ranges over all prime numbers. 
Consider the following condition $C(v)$.

\medskip

$C(v)$: $M$ satisfies the weight-monodromy conjecture at $v$. Furthermore, for any $a\in E_v$, the element $\langle \Ad(N_{v,0}^{d-2})(\iota_v(a)), \iota_v(a)\rangle$ of $\Q_{\ell}$ ($\iota_v$ is the $\iota$ in \ref{loc2} for the local field $F_v$), where $\langle\;,\;\rangle$ is the paring defined by the polarization,  is contained in  $\Q_{\geq 0}$ for any prime number $\ell$ and it is  independent of $\ell$, and it is $0$ if and only if $a=0$. 

\medskip

In the situation I, this condition is the same as the following: For any $a\in \text{Ext}^1(\Q, M)$, the corresponding mixed motive  $M_a$ with the exact sequence $0\to M \to M_a \to \Q\to 0$ satisfies (i) and (ii) in \ref{Mdgeq2} at $v$.

Under the assumption that $C(v)$ is satisfied for all non-archimedean places $v$ of $F$, we define an element $H(a)\in \R_{>0}$ as follows. For a non-archimedean place $v$ of $F$, let $H_v(a):=\sharp(\F_v)^{l(v)}$  where $l(v)= \langle \Ad(N_{v,0}^{d-2})(\iota_v(a)), \iota_v(a)\rangle^{1/d}$. If $v$ is archimedean, let $H_v(a):= \exp(2 (\delta_v(a), \delta_v(a))^{1/d})$ if $v$ is a complex place, and let $H_v(a):= \exp((\delta_v(a), \delta_v(a))^{1/d})$ if $v$ is a real place. Here $(\;,\;)$ is the Hodge metric.
Define $H(a):= \prod_v H_v(a)$ where $v$ ranges over all places of $F$.
In the situation I,  we have $H_v(a)=H_{\diamondsuit,0,d, v}(M_a)$ and $H(a)= H_{\diamondsuit,0,d}(M_a)$ where $M_a$ is the extension $0\to M\to M_a\to \Z\to 0$ corresponding to $a$. 

\end{sbpara}

\begin{sbpara}\label{Sabc}
We consider relations between the following statements $S(a)$, $S(b)$, $S(c)$ concerning the subjects (a), (b), (c) in 2.1.1, respectively, and a stronger version $S(b)^*$ of $S(b)$. In these statements, $S$ denotes a finite set of places of $F$ containing all archimedean places. 

\medskip

$S(a)$: (Here we assume that the condition $C(v)$ is satisfied for any non-archimedean place $v$ of $F$.) For any $C\in \R_{> 0}$, there are only finitely many $a\in E^{\Z}$ such that $H(a)\leq C$.

$S(b)$: $E^{\Z}_S$ is finitely generated as an abelian group for any $S$.  

$S(b)^*$: $E^{\Z}_S$ is finitely generated as an abelian group for any $S$. Furthermore, the map $$\R\otimes_\Q E_{S_{\infty}}\to \oplus_{v|\infty} E^{\R}_v$$ is injective. Here $S_{\infty}$ denotes the set of all archimedean places of $F$.

$S(c)$: The map $\R\otimes_\Q E_{S_{\infty}}\to \oplus_{v|\infty} E^{\R}_v$ is injective. Furthermore, we have an isomorphism $${\bf A}^f_\Q\otimes_{\Q} E_S \overset{\cong}\to H^1_{g,S}(F, M_{et})$$
for any $S$. 

\medskip
In the situation II, the injectivity of $\R\otimes_\Q E_{S_{\infty}}\to \oplus_{v|\infty} E^{\R}_v$ is conjectured in \cite{Be}, the bijectivity of it in the case $d\geq 3$ is conjectured in \cite{Be}, and the bijectivity of 
${\bf A}^f_\Q\otimes_{\Q} E_S \overset{\cong}\to H^1_{g,S}(F, M_{et})$ is conjectured in \cite{BK}. 
\end{sbpara}

\begin{sbprop}\label{2.1thm} Let $F$ be a number field, let $M$ be a polarized pure motive with $\Z$-coefficients over $F$ of weight $-d$ with $d\geq 2$, and let $E$ be as in \ref{Egiven}. 

(1) Assume that $C(v)$ is satisfied at any non-archimedean place $v$ of $F$. Then  
we have the equivelence
$$S(a) \Leftrightarrow S(b)^*.$$ 

(2) We have the implication $$S(c) \Rightarrow S(b)^*.$$ 

\end{sbprop}

\begin{sbpara}\label{MWmethod} This is a preparation for the proof of $S(a) \Rightarrow S(b)^*$  in (1) of  Prop. \ref{2.1thm}.

Recall that Mordell-Weil theorem for abelian varieties is proved by using the following fact.

Let $\Gamma$ be an abelian group and assume the following (i) and (ii) are satisfied. Then $\Gamma$ is finitely generated. 

\medskip

(i) For some integer $n\geq 2$, the quotient group $\Gamma/n\Gamma$ is finite.

(ii) There is a symmetric $\Z$-bilinear form $\langle\;,\;\rangle: \Gamma \times \Gamma \to \R$ such that $\langle \gamma, \gamma\rangle \geq 0$ for any $\gamma\in \Gamma$ and such that for any $C\geq 0$, the set $\{\gamma\in \Gamma\;|\; \langle \gamma, \gamma \rangle\leq C\}$ is finite. 

\medskip

(In the proof, $\Gamma$ is the Mordell-Weil group and $\langle\;,\;\rangle$ is the height pairing.)
\end{sbpara}

\begin{sbpara}

We prove $S(a) \Rightarrow S(b)^*$ in (1) in Proposition \ref{2.1thm}. We assume $S(a)$ holds.

For $\Gamma =E^{\Z}_S$, (i) in \ref{MWmethod} is satisfied by \ref{WMW}, and 
the pairing $\langle \;,\;\rangle$ in (ii) in \ref{MWmethod}  is given by the symmetric bilinear form on $(\oplus_{v|\infty} E^\R_v) \oplus (\oplus_{v\in S\smallsetminus S_{\infty}} E_v)$. 
 The condition (ii) in \ref{MWmethod} is exactly the finiteness of the number of elements of $E_S^{\Z}$ of bounded height.
 Hence we have the implication $S(a) \Rightarrow S(b)$.

We prove that the map $f: V_1:= \R \otimes_\Q E_{S_{\infty}} \to V_2:=\oplus_{v|\infty} E^{\R}_v$ is injective. Let $L=
E^{\Z}_{S_{\infty}}/(\text{torsion})$. $S(a)$ shows that (*) the restriction of  $f$ to $L$ is injective and (**) $f(L)$ is discrete in $V_2$. Since $f(L)$ is a discrete subgroup of $f(V_1)$ by (**), we have $\text{rank}_\Z f(L) \leq \dim_\R f(V_1)$. By (*), $\text{rank}_\Z f(L)= \text{rank}_\Z L =\dim_\R V_1$. Hence  $\dim_\R f(V_1)=\dim_\R V_1$ and hence $f$ is injective.

\end{sbpara}

\begin{sbpara}\label{EMZ2} 
This is a preparation for the proof of $S(b)^*\Rightarrow S(a)$. 

Assume $d=2$. We show that for $C\in \R_{>0}$, 
$\{a\in E^{\Z}\;|\; H(a) \leq C\}\subset E^{\Z}_S$ for some finite set $S$ of non-archimedean places of $F$. 

Let $c$ be the unique element of $\Q_{>0}$ such that the image of the homomorphism $M_{et,\hat \Z}(-1) \otimes M_{et,\hat \Z}(-1) \to {\bf A}^f_{\Q}$ induced by the polarization of $M$ coincides with  $c \cdot \hat \Z\subset {\bf A}^f_\Q$. 
We show that if $S$ satisfies the following (i)--(iii), then $\{a\in E^{\Z}\;|\; H(a) \leq C\}\subset E^{\Z}_S$. 

(i) $M$ is of good reduction outside $S$. (ii) If $v$ is a place of $F$ outside $S$, then $\sharp(\F_v) > C^{1/\sqrt{c}}$. (iii) If $v$ is a place of $F$ outside $S$ and if $p=\text{char}(\F_v)$, there is $r\in \Z$ such that 
$F^rM_{\dR}=M_{\dR}$ and $F^{r+p-1}M_{\dR}=0$.

Let $a\in E^{\Z}$ and assume $H(a) \leq C$. Let $v$ be a place of $F$ outside $S$. We prove that for any prime number $\ell$,  the image of $a$ in $H^1_g(F_v, M_{et,\Q_{\ell}})$ belongs to $H^1_f(F_v, M_{et,\Q_{\ell}})$. It is sufficient to prove that $\iota_v(a)=0$. 
If $\ell\neq p$, $\langle \iota_v(a), \iota_v(a)\rangle$ defined by using  $\ell$ belongs to $ c \cdot \Z_{\ell}$. Assume $\ell=p$. By the conditions (i) and (iii), we have the Fontaine-Laffaille module $\sD_{\text{FL}}(M_{\et, \Z_p})$ and the polarization of $M$ induces 
$\sD_{\text{FL}}(M_{et,\Z_p})(-1)\times \sD_{\text{FL}}(M_{et,\Z_p})(-1)\to c  \sD_{\text{FL}}(\Z_p)=c\cdot W(\F_v)$, and hence $\langle \iota_v(a), \iota_v(a)\rangle$ defined by using $\ell=p$ belongs to $ c \cdot W(\F_v)$. 
Since we assume that $\langle \iota_v(a), \iota_v(a)\rangle$ does not depend on $\ell$ and belongs to $\Q_{\geq 0}$, we have 
$\langle \iota_v(a), \iota_v(a)\rangle\in c \cdot\Z_{\geq 0}$. Hence if $\iota_v(a)\neq 0$, we have $H_v(a) \geq \sharp(\F_v)^{c^{1/2}}$. But $C\geq H(a)\geq H_v(a)$, and this would imply $C^{1/\sqrt{c}}\geq \sharp(\F_v)$, a contradiction.

\end{sbpara}

\begin{sbpara}\label{EMZ3} We prove the implication $S(b)^*\Rightarrow S(a)$ in \ref{2.1thm} (1). Assume $S(b)^*$.  We show that for $C\in \R_{>0}$, the set $\{a\in E^{\Z}\;|\; H(a) \leq C\}$ is finite. By \ref{EM6} and \ref{EMZ2}, it is sufficient to prove that the set $\{a\in E^{\Z}_S\;|\; H(a) \leq C\}$ is finite.
  
  Let $v$ be a non-archimedean place $v$ of $F$. By the assumption, $E_v$ is of finite dimension over $\Q$. The positive definite symmetric bilinear form on $E_v$ induces
     a positive definite 
     symmetric  bilinear form $E^{\R}_v \times E^{\R}_v\to \R$ over $\R$ where $E^{\R}_v=\R\otimes_\Q E_v$. 

 $E^{\Z}_S$ is a finitely generate $\Z$-module, we have an embedding of finite-dimensional $\R$-vector spaces 
 $\R\otimes_\Z E^{\Z}_S\to \oplus_{v\in S} E^{\R}_v$, and we have a positive definite symmetric bilinear form $\langle\;,\rangle$ on $\oplus_{v\in S} E^{\R}_v$. If $C\in \R_{>0}$, there is $C'\in \R_{\geq 0}$ such that $H(a) \leq C$ for $a\in E^{\Z}_S$ implies $\langle a, a\rangle\leq C'$, and hence
  there are only finitely many such $a$.

\end{sbpara}

\begin{sbpara}\label{2.1thm2}  We prove (2) of \ref{2.1thm}, that is, the implication $S(c) \Rightarrow S(b)^*$. 

Assume that ${\bf A}^f_\Q\otimes_\Q E_S\to H^1_{g,S}(F, M_{et})$ is an isomorphism. We show that $E_S^{\Z}$ is finitely generated as an abelian group.  Since the $\Q_{\ell}$-component of $H^1_{g,S}(F, M_{et})$ is contained in the finite-dimensional $\Q_{\ell}$-vector space $H^1_{et}(\Spec(O_{F, S'}), M_{et, \Q_{\ell}})$ for a finite set $S'$ of places of $F$ which contains $S$ and all prime divisors of $\ell$ in $F$ such that $M$ is of good reduction outside $S'$, we have that $E_S$ is of finite dimension over $\Q$. Let $\hat L$ be the fiber product of 
$$H^1_{g, S}(F, M_{et})\to H^1_{\text{cont}}(\Gal(\bar F/F), M_{et})\leftarrow H^1_{\text{cont}}(\Gal(\bar F/F), M_{et, \hat \Z}).$$
Then $\Q\otimes \hat L= H^1_{g,S}(F, M_{et})$. The torsion part of $\hat L$ is isomorphic to $H^0(F, M_{et}/M_{et, \hat \Z})$ which is finite by Weil conjecture proved by Deligne. Each $\ell$-adic component of $\hat L$ is a finitely generated $\Z_{\ell}$-module because it is contained in the finitely generated $\Z_{\ell}$-module $H^1_{et}(\Spec(O_{F,S'}), M_{et,\Z_{\ell}})$ for $S'$ as above. These prove that $\hat L/(\text{torsion})$ is identified with a $\hat \Z$-lattice in ${\bf A}^f_{\Q} \otimes_\Q E_S$ and these prove that $E^{\Z}_S$, which is the fiber product of 
$E_S\to {\bf A}^f_{\Q} \otimes_{\Q} E_S= H^1_{g,S}(F, M_{et}) \leftarrow \hat L$, is finitely generated as an abelian group.

\end{sbpara}

\begin{sbpara}\label{MWd=1}  
Let $F$ be a number field and let $M$ be a polarized pure motive with $\Q$-coefficients over $F$ of weight $-1$. 

Let $E$ be a $\Q$-vector space and assume that we are given a $\Q$-linear map $E\to H^1_g(F, M_{et})$.
 We are interested in the following two situations I, II. In the situation I, $E=\text{Ext}^1(\Q, M)$ and the above map is given by \'etale realizations. In the situation II, 
$E=(\Q\otimes CH^r(X))_0$ assuming $M$ is $H^{2r-1}(X)(r)$ for a projective smooth scheme $X$ over $F$ and for $r\in \Z$ and the above map is given by the Chern class map. We expect that $E$ is $\text{Ext}^1(\Q, M)$ also in the situation II and the map coincides with the map of the situation I, but we like to discuss the situation II not assuming these are true.  Assume that we are given a  positive definite symmetric $\Q$-bilinear form $E\times E\to \R$ and  a structure of $\Z$-coefficients of $M$. Let $E^{\Z}$ be the fiber product 
of $E\to H^1_{\text{cont}}(F, M_{et})\leftarrow H^1_{\text{cont}}(\Gal(\bar F/F), M_{et,\hat \Z})$. Then by analogous arguments as in the proof of \ref{2.1thm}, we have the implications 
$$S(a) \Leftrightarrow S(b) \Leftarrow S(c)$$
between the following statements $S(a)$, $S(b)$, $S(c)$. 
(In the situation II, this S(c) is conjectured in \cite{BK} like S(c) in \ref{Sabc}.)

 \medskip

S(a): The set $\{a\in E^{\Z}\;|\; \langle a, a\rangle\leq C\}$ is finite for any $C\in \R_{\geq 0}$.

S(b): $E^{\Z}$ is finitely generated as an abelian group.

S(c): ${\bf A}^f_\Q\otimes_\Q E \overset{\cong}\to H^1_g(F, M_{et})$.

\end{sbpara}

\begin{sbpara}\label{fin} Assume that the conjectures (i) and (ii) in \ref{Mdgeq2} and the conjectures (i)--(iv) in \ref{BeBl} (which appeared when we discussed the definition of the height function $H_{\diamondsuit}$) are true. Assume that $S(c)$ in \ref{Sabc} and \ref{MWd=1} are true in the situations I. 

Then we can prove (1) in \ref{(a)}, that is, the finiteness of the number of isomorphism classes of mixed motives with $\Z$-coefficients over $F$ with fixed pure graded quotients and of bounded height,  as follows. 

By induction on the length of the weight filtration of $M$, we may assume that
$W_wM=M$, $W_{w-d-1}M=0$ ($w\in \Z$ and $d\geq 1$ are fixed), and the subquotients $W_wM/W_{w-d}M$ and $W_{w-1}M$ of $M$ are fixed. Then if such $M$ exists, isomorphism classes of such $M$ form a torsor over $\text{Ext}^1(\gr^W_wM, \gr^W_{w-d}M)=\text{Ext}^1(\Z, M')$ where $M'= (\gr^W_wM)^*\otimes\gr^W_{w-d}M$. This can be seen from the exact sequence
$$0\to \text{Ext}^1(\gr^W_wM, \gr^W_{w-d}M) \to \text{Ext}^1(\gr^W_wM, W_{w-1}M)\to \text{Ext}^1(\gr^W_w, M/W_{w-d}M)$$ in which the last arrow sends the class of various $M$ to the  class of the fixed $M/W_{w-d}M$. Write the action of $a\in \text{Ext}^1(\Z, M')$ on the set of isomorphism classes of such $M$ as $M\mapsto M+a$. 
Then there  are constants $c,c'\in \R_{>0}$ such that
$H_{\diamondsuit,0,d}(a)\leq c'H_{\diamondsuit}(M)^cH_{\diamondsuit}(M+a)^c$ for a fixed $M$ and for any $a\in \text{Ext}^1(\Z, M')$,  where $H_{\diamondsuit,0,d}(a)$ denotes $H_{\diamondsuit}(M'_a)$ for  the extension  $0\to M'\to M'_a\to \Z\to 0$ given by $a$. This shows that if $H_{\diamondsuit}(M)$ for such $M$ are bounded, then $H_{\diamondsuit,0,d}(a)$ is bounded, and hence by \ref{2.1thm} and \ref{MWd=1}, we have the finiteness of the classes of such $M$.

\end{sbpara}

\subsection{Period domains}\label{s:2.2}

In this \S\ref{s:2.2}, we consider a period domain $X(\C)$ which classifies Hodge structures, a  toroidal partial compactification $\bar X(\C)$ of $X(\C)$, and a set $X(F)$ of motives over a number field $F\subset \C$ with a map $X(F)\to X(\C)$, as in 
$$X(F) \to  X(\C) \subset \bar X(\C).$$

\begin{sbpara}\label{Phi}
Fix
$$\Phi=((h(w,r))_{w,r\in\Z}, H_{0,\Q}, W, (\langle\;,\;\rangle_{0,w})_{w\in \Z})$$
where  $h(w,r)$ are integers $\geq 0$ given for $w,r\in \Z$ such that $h(w,r)=0$ for almost all $(w,r)$ and such that $h(w,w-r)=h(w,r)$ for all $(w,r)$, $H_{0,\Q}$  is a $\Q$-vector space of dimension $\sum_{w,r} h(w,r)$, $W$ is an increasing filtration  on $H_{0,\Q}$ such that $\dim_\Q(\gr^W_w)= \sum_r h(w,r)$ for any $w$,  and $\langle\;,\;\rangle_{0,w}: \gr^W_w\times \gr^W_w\to \Q(-w)=\Q\cdot (2\pi i)^{-w}$ is a non-degenerate bilinear form for each $w\in \Z$ which is symmetric if $w$ is even and anti-symmetric if $w$ is odd.

\end{sbpara}

\begin{sbpara}\label{Dmu10a}
Define algebraic groups $G_1$ and $G$ over $\Q$ as 
$$G_1=\{g\in \text{Aut}(H_{0,\Q})\;|\; gW_{\bullet}=W_{\bullet},\; \langle gx,gy\rangle_{0,w}=\langle x,y\rangle_{0,w}\;(w\in \Z, x,y\in \gr^W_w)\},$$
$$G= \{(g,t)\in \text{Aut}(H_{0,\Q})\times {\bf G}_m\;|\; gW_{\bullet}=W_{\bullet},\; \langle gx, gy\rangle_{0,w}= t^w\langle x,y\rangle_{0,w}\;(w\in \Z, x,y\in \gr^W_w)\}.$$
We identify $\Lie(G_1)$ with $\text{End}_{\Q,W, \langle\;,\;\rangle}(H_{0,\Q})$.

Let $G_u= \{g \in G_1\; |\; (g-1)W_w\subset W_{w-1}\;\forall \;w\in \Z\}$ be the unipotent radical of $G_1$ which is also identified with the unipotent radical of $G$.

\end{sbpara}

\begin{sbpara}\label{Dmu1}

Let $\check{D}$ be the set of all decreasing filtrations $F$ on $H_{0,\C}=\C\otimes_\Q H_{0,\Q}$ satisfying the following (i).

\medskip

(i) $\dim_{\C}\gr^r\gr^W_wH_{0,\C}= h(w,r)$ for any $w,r\in \Z$ and the annihilator of $F^r\gr^W_wH_{0,\C}$ with respect to $\langle\;,\;\rangle_w$ coincides with $F^{w+1-r}\gr^W_wH_{0,\C}$ for any $w,r\in \Z$. 

\medskip

Let $D\subset \check{D}$ (resp. $D^{\pm}\subset \check{D}$)  be the subset consisting of all $F\in \check{D}$ satisfying the following (ii).

\medskip

(ii) The triples $(\gr^W_wH_{0,\Q}, \gr^W_w(F), \langle\;,\;\rangle_{0,w})$ (resp. There is  $\epsilon \in \{\pm 1\}$ such that the triples $(\gr^W_wH_{0,\Q}, \gr^W_w(F), \epsilon^w\langle \;, \;\rangle_{0,w})$) are polarized Hodge structures of weight $w$ for all $w\in \Z$.

\medskip

Then $D$ (resp. $D^{\pm}$) is an open set of the complex analytic manifold $\check{D}$ and hence is also a complex analytic manifold. $D^{\pm}=D$  if $\gr^W_wH_{0,\Q}=0$ for all odd integers $w$, and $D^{\pm} = D \coprod gD$ for some $g\in G(\Q)$ otherwise.

It is known that the natural action of $G_1(\R)G_u(\C)$ (resp. $G(\R)G_u(\C)$) on $D$ (resp. $D^{\pm}$) is transitive.

If the weight filtration $W$ is pure, $D$ is the classical period domain of Griffiths \cite{G1}. In general, $D$ is the generalization of Griffiths period domain to mixed Hodge structure in Usui \cite{U}.

\end{sbpara}

\begin{sbpara}\label{neat2} Let $K$ be an open compact subgroup of $G({\bf A}^f_\Q)$ satisfying the following condition (*). 
Let $G_{(1)}= \{(g, t)\in G\;|\; g\in G_1\}=\{(g,t)\in G\;|\; t^w=1\;\text{if}\; \gr^W_wH_{0,\Q}\neq 0\}$. 

\medskip

(*) If $\nu\in G({\bf A}^f_\Q)$ and  if $(\gamma_1,\gamma_2)\in G_{(1)}(\Q)\cap \nu K \nu^{-1}\subset G({\bf A}^f_\Q)$, 
the subgroup of $\C^\times$ generated by eigenvalues of $\gamma_1: H_{0,\C}\to H_{0,\C}$ is torsion free.

\medskip

For example, if $H_{0,\Z}$ is a $\Z$-lattice in $H_{0,\Q}$ and $n\geq 3$ is an integer, and if any element $(g,t)$ of $K$ satisfies the following (i) and (ii), then $K$ satisfies the neat condition (*).

(i) $g H_{0, {\hat \Z}}= H_{0,\hat \Z}$ where $H_{0,\hat \Z}=\hat \Z\otimes_{\Z} H_{0,\Z}$. (ii) The automorphism of $H_{0, \Z}/nH_{0,\Z}$ induced by $g$ is the identity map.

\end{sbpara}

\begin{sbpara}\label{XK}
Let
$$X(\C)=X_K(\C):=G(\Q) \bs (D^{\pm} \times (G({\bf A}^f_\Q)/K))=G(\Q)^+ \bs (D \times (G({\bf A}^f_\Q)/K))$$
where $G(\Q)^+=\{(g,t)\in G(\Q)\;|\; t^w>0 \;\text{if}\;\gr^W_wH_{0,\Q}\neq 0\}$.

For a representative $(\nu_i)_{1\leq i\leq n}$ ($\nu_i\in G({\bf A}^f_\Q)$) of the finite set $G(\Q)^+\bs G({\bf A}^f_\Q)/K$, the maps $D \to X(\C)\;;\;F \mapsto \text{class}(F,\nu_i)$ for $1\leq i\leq n$ induce a bijection

(*) $\;\;\coprod_{i=1}^n \Gamma_i\bs D  \overset{\cong}\to X(\C)$

\noindent
where $\Gamma_i$ is the image of $G(\Q)^+\cap \nu_iK\nu_i^{-1}= G_{(1)}(\Q) \cap G({\bf A}^f_\Q)$ in $G_1(\Q)$. 

Since $D\to \Gamma_i \bs D$ are homeomorphisms by the neat condition on $K$, we have a unique structure of a complex analytic manifold on $X(\C)$ for which 
the map $D^{\pm}\times G({\bf A}^f_\Q)/K\to X(\C)$ is locally an isomorphism and for which the above (*) is an isomorphism of complex analytic manifolds.

\end{sbpara}

\begin{sbpara}\label{X(F)} 

For a number field $F\subset \C$, let $X(F)$ be (resp. Let $X'(\C)$ be) the set
 of isomorphism classes of mixed motives with $\Q$-coefficients over $F$
 (resp. of mixed $\Q$-Hodge structures) such that 
$$\dim_F \gr^r \gr^W_w M_{dR} = h(w,r)\quad (\text{resp.} \;\; \dim_\C\gr^r\gr^W_wH_\C=h(w,r))\;\;\text{for any}\;w,r$$ 

\medskip
\noindent
endowed with a {\it polarized $K$-level structure} in the following sense. 

For a mixed motive  $M$  with $\Q$-coefficients over $F$ (resp. For a mixed $\Q$-Hodge structure $H$), a polarized $K$-level structure on $M$ (resp. $H$) means the equivalence class of a
quadruple   $(\la_1, \theta_1, \la_2, \theta_2)$ where 
\medskip

$\la_1$ is an isomorphism of ${\bf A}^f_\Q$-modules ${\bf A}^f_\Q\otimes_\Q H_{0, \Q} \overset{\cong}\to   M_{et}$ such that there is a homomorphism $k=(k_1, k_2): \Gal(\bar F/F)\to K$ satisfying $\sigma \la_1(x)= \la_1(k_1(\sigma)x)$  for any $\sigma\in \Gal(\bar F/F)$ and any $x\in M_{et}$ and $k_2: \Gal(\bar F/F)\to ({\bf A}^f_\Q)^\times$ is the inverse of the cyclotomic character (resp. $\la_1$ is an isomorphism of ${\bf A}^f_\Q$-modules 
${\bf A}^f_\Q\otimes_\Q H_{0,\Q}\overset{\cong}\to {\bf A}^f_\Q\otimes_\Q H_\Q$), 
which preserves the weight filtrations, 

\medskip

$\theta_1$ is an isomorphism of $\Q$-vector spaces $M_B\overset{\cong}\to H_{0,\Q}$ 
where $M_B$ is defined with respect to the embedding $F\to \C$ which we fixed (resp. an isomorphism of $\Q$-vector spaces $H_\Q\overset{\cong}\to H_{0,\Q}$),
which preserves weight filtrations, and 

\medskip

 $\la_2\in ({\bf A}^f_\Q)^\times$, $\theta_2\in \Q^\times$,
 
 \medskip
 \noindent
  satisfying the following conditions (*) and (**).

\medskip

(*) There is $\epsilon\in \{\pm 1\}$ such that we have  a polarization on $\gr^W_wM$ (resp. $\gr^W_wH$) for each $w$ which induces (resp. which is given by)
$$\gr^W_w M_B \times \gr^W_w M_B \;(\text{resp.}\;\;\gr^W_wH_\Q\times \gr^W_w H_\Q) \to \Q\cdot (2\pi i)^{-w}\;;\;(x,y)\mapsto \epsilon^w\langle \theta_1(x), \theta_1(y)\rangle_{0,w}.$$

\medskip

(**) $({\bf A}^f_\Q \otimes \theta_1)\circ \la_1\in G({\bf A}^f_\Q)$. 

\medskip

By our definition, such $(\la_1, \theta_1, \la_2, \theta_2)$ and $(\la'_1,\theta'_1, \la'_2, \theta'_2)$ are equivalent if and only if there are $(k_1, k_2)\in K$ and $(\gamma_1, \gamma_2)\in G(\Q)$ satisfying the following (i)--(iii).

\medskip
(i) $\la'_1=  \la_1\circ k_1$. (ii) $\theta'_1= \gamma_1 \circ \theta_1$. 
(iii)  $\la'_2 \theta'_2= \la_2 \theta_2k_2\gamma_2$. 

\medskip

By our definition, $M$ (resp. $H$) with $\text{class}(\la_1, \theta_1, \la_2, \theta_2)$ and $M'$ (resp. $H'$) with $\text{class}(\la'_1, \theta'_1, \la'_2, \theta'_2)$ are isomorphic if and only if there is an isomorphism $\alpha: M\overset{\cong}\to  M'$ (resp. $H\overset{\cong}\to H'$) such that $\text{class}(\la'_1, \theta'_1, \la'_2, \theta'_2)$ coincides with $\text{class}(\alpha \circ \la_1, \theta_1\circ \alpha^{-1}, \la_2, \theta_2)$. 

We have a canonical bijection 
 $X(\C) \to X'(\C)$ which sends the class of $(F, (g,t))$ ($F\in D^{\pm}, (g,t)\in G({\bf A}^f_\Q)$) to the class of 
$(H, \la_1, \theta_1, \la_2, \theta_2)$ where $H$ is $H_{0,\Q}$ endowed with the Hodge filtration $F$, $\la_1=g$, $\theta_1$ is the identify map, and $\la_2=t$, $\theta_2=1$. The converse map $X'(\C)\to X(\C)$ 
sends the class of $(H, \la_1,\theta_1, \la_2, \theta_2)$ to the class of $(F, (g,t))$ where 
$F$ is the image of the Hodge filtration of $H$ under $\theta_1$, $g= ({\bf A}^f_\Q\otimes \theta_1) \circ \la_1$, and $t=\la_2\theta_2$.

We have a canonical map
$X(F)\to X'(\C)$ by taking the mixed Hodge structure associated  to a mixed motive. Hence we have a canonical map $X(F)\to X(\C)$. 

If $F\subset F'\subset \C$ and $F'$ is a number field, we have a canonical map $X(F) \to X(F')$. 
\end{sbpara}

\begin{sbprop}\label{Aut} (1) Let $M$ be a mixed motive with $\Q$-coefficients over a number field $F\subset \C$ endowed with a polarized $K$-level structure $\la$. Then $\text{Aut}(M, \la)=\{1\}$.

(2) Let $H$ be a mixed $\Q$-Hodge structure endowed with a polarized $K$-level structure $\la$. Then $\text{Aut}(H, \la)=\{1\}$. 
\end{sbprop}

\begin{pf} We write the proof for (2). The proof for (1) can be given in the same way (actually (1) follows from (2) by applying (2) to $M_H$). Let $\alpha$ be an automorphism of $(H,\la)$. Let $(\la_1,\theta_1,\la_2, \theta_2)$ be a representative of $\la$. Then since $(\alpha\circ\la_1, \theta_1\circ \alpha^{-1}, \la_2, \theta_2)$ is also a representative of $\la$, there are $(\gamma_1, \gamma_2)\in G(\Q)$ and $(k_1,k_2)\in K$ satisfying 

(i) $\alpha\circ \la_1= \la_1\circ k_1$,\;\; (ii) $\theta_1\circ \alpha^{-1}= \gamma_1 \circ \theta_1$,\;\; (iii) $k_2\gamma_2=1$.

From (i)--(iii), we have 

(iv) $(\gamma_1,\gamma_2)= (g,t)(k_1, k_2)^{-1}(g,t)^{-1}$ where $(g,t)\in G({\bf A}^f_\Q)$ is defined by $g= ({\bf A}^f_\Q\otimes \theta_1)\circ \lambda_1$ and $t= \theta_2\la_2$.

 For some $\epsilon \in \{\pm 1\}$, $(\gr^W H_{0,\Q}, \gr^W_w(\theta_1F), \epsilon^w\langle\;,\;\rangle_{0,w})$ is a polarized Hodge structure of weight $w$ for any $w\in \Z$. Since $\alpha$ fixes the Hodge filtration $F$ of $H$, 
 
 (v) $\gamma_1$  preserves the Hodge filtration $\theta_1F$. 
 
 Hence $(\gr^WH_{0,\Q}, \gr^W_w(\theta_1F), \epsilon^w\langle \gamma_1(\;), \gamma_1(\;)\rangle_{0,w})$ is also a  polarized Hodge structure of weight $w$ for any $w\in \Z$.
 Since $\langle\gamma_1(\;), \gamma_1(\;)\rangle_{0,w}= \gamma_2^w\langle \;,\;\rangle_{0,w}$, we have 
 
 (vi) $\gamma_2^w>0$ for any $w\in \Z$ such that $\gr^W_wH_{0,\Q}\neq 0$. That is, $(\gamma_1,\gamma_2)\in G(\Q)^+$.

 By (iv) and (vi), we have

 (vii) $\gamma_1\in G_1(\Q)$. 
 
 By (v), $\gamma_1$ belongs to the isotropy group of the associated Hodge metric of the polarization which is a compact subgroup of $G_1(\R)$. 
Since the image of $G_{(1)}(\Q)\cap (g,t)K(g,t)^{-1}$ in $G_1(\R)$ is discrete, we see that 

(viii) $\gamma_1$ is of finite order.

By (iv) (vii) and (viii), and by the neat condition on $K$, 
we have $\gamma_1=1$. Hence $\alpha=1$ by (ii). 
\end{pf}

\begin{sbrem}\label{fdef} (1) By using \ref{Aut}, we can show that for a finite Galois extension $F'\subset \C$ of $F$, the map from $X(F)$ to the $\Gal(F'/F)$-fixed part of $X(F')$ is bijective though we do not give the proof in this paper.

(2) We expect that the map $X(F)\to X(\C)$ is injective. 

(3) It is probable that any $M\in X(F)$ is semi-stable (\ref{semi-st}) as a mixed motive with $\Q$-coefficients over $F$. 

(4) Proposition \ref{Aut} and the statement in the above (1) remain true even if the neat condition on $K$ is replaced by a weaker condition that the image of $G_{(1)}(\Q)\cap \nu K\nu^{-1}$ in $G_1(\Q)$ is torsion free for any $\nu\in G({\bf A}^f_\Q)$. But the neat condition is important in the theory of toroidal partial compactifications which appear below. 
\end{sbrem}

\begin{sbpara}\label{canpol}

For a mixed motive $M$  with $\Q$-coefficients over a number field $F\subset \C$ (resp. For a mixed $\Q$-Hodge structure $H$) endowed with a polarized $K$-level structure,  we have  a unique family of polarizations $(\langle\;,\;\rangle_w)_{w\in \Z}$ on $\gr^W_wM$ (resp. $\gr^W_wH$), which we call the canonical polarizations, satisfying the following condition. For any 
representative $(\la_1, \theta_1, \la_2, \theta_2)$ of this polarized $K$-level structure, there is $\epsilon\in \{\pm 1\}$ such that for each $w\in \Z$, $\langle\;,\;\rangle_w$ induces (resp.  coincides with) the  
map $\gr^W_wM_B \times \gr^W_wM_B$ (resp. $\gr^W_w H_\Q \times \gr^W_wH_\Q$) $\to \Q\cdot (2\pi i)^{-w}$  
given by
 $(x,y)\mapsto \epsilon^w|\theta_2\la_2|^w_f\langle \theta_1(x), \theta_1(y)\rangle_{0,w}$. Here $|\;|_f$ denotes the absolute value $(a_v)_{v\neq \infty}\mapsto \prod_{v \neq \infty} |a_v|_v$ on $({\bf A}^f_\Q)^\times$. 
 In fact, if $(\la_1,\theta_1, \la_2, \theta_2)$ and $(\la'_1, \theta'_1,\la'_2, \theta'_2)$ are two representatives and if $(k_1,k_2)$ and $(\gamma_1, \gamma_2)$ are as in \ref{X(F)},  
 by  (iii) in \ref{X(F)} and $|k_2|_f=1$, we have
 $$|\theta'_2\la'_2|^w_f \langle \theta'_1 x, \theta' _1 y\rangle_{0,w}
 = |\theta'_2\la'_2|^w_f \langle \gamma_1 \theta_1 x, \gamma_1\theta_1  y\rangle_{0,w}$$
 $$= |\theta'_2\la'_2|^w_f\gamma_2^w\langle \theta_1 x, \theta_1 y\rangle_{0,w}=
|\theta'_2\la'_2\gamma_2^{-1}|_f^w \delta^w \langle \theta_1x, \theta_1y\rangle_{0,w}= |\theta_2\la_2|^w_f\delta^w \langle \theta_1 x, \theta_1 y\rangle_{0,w}$$
 where $\delta$ is the sign of $\gamma_2$. 
 
 Canonical polarizations are preserved by isomorphisms between mixed motives with $\Q$-coefficients over $F$ (resp. mixed $\Q$-Hodge structures) endowed with polarized $K$-level structures. 

\end{sbpara}

\begin{sbpara}\label{Siegel} {\bf Example.} Let $g\geq 1$ and let $\Phi$ be as follows. $h(w,r)=g$ if $(w,r)=(1,0), (1,1)$, and $h(w,r)=0$ otherwise. $H_{0,\Q}= \Q^{2g}$ with the standard basis $(e_i)_{1\leq i\leq 2g}$. $W_1=H_{0, \Q}$ and $W_0=0$. $\langle\;,\;\rangle_{0,1}$ is such that $\langle e_i, e_j\rangle_{0,1}$ is $(2\pi i)^{-1}$ if $i=g+j$, $-(2\pi i)^{-1}$ if $j=i+g$, and is $0$ otherwise. 

Then $D$ is the upper half space of Siegel of degree $g$ and $X(\C)$ is the space of $\C$-points of the Siegel modular variety of degree $g$ and of level $K$.  

For a number filed $F\subset \C$, let $X'(F)$ be the set of $F$-points of this modular variety, which is a fine moduli space of principally polarized abelian varieties $A$ of dimension $g$ with $K$-level structure. We have a canonical map $X'(F)\to X(F)\;;A\mapsto H^1(A)$. We expect that the last map  is a bijection.

\end{sbpara}

\begin{sbpara}\label{Dmu10}\label{frak S}

We define a toroidal partial compactification $\bar X(\C)$ of $X(\C)$ using our previous  works  \cite{KU} and  \cite{KNU} with S. Usui and C. Nakayama. 

Let $\Sig$ be the set of all cones  in $\Lie(G_1)_\R:=\R \otimes_\Q \Lie(G_1)$ of the form $\R_{\geq 0}N$ where  $N$ ranges over all elements of $\Lie(G_1)$ satisfying the following condition: $N$ is nilpotent as a linear operator on $H_{0,\Q}$ and the relative monodromy filtration of $N$ with respect to $W$ (\ref{2.4.b}) exists. 

We have the set $D_{\Sig}\supset D$ (resp. $D^{\pm}_{\Sig}\supset D^{\pm}$) defined as follows. $D_{\Sig}$ (resp. $D^{\pm}_{\Sig}$) is the set of all pairs $(\sig, Z)$ where $\sig\in \Sig$ and $Z$ is a non-empty subset of $\check{D}$ satisfying the following conditions (i)--(iii). Write $\sig=\R_{\geq 0}N$ with $N\in \Lie(G_1)$, and let $F\in Z$.

\medskip

(i) $Z=\exp(\C N)F$. 

(ii) $N F^p\subset F^{p-1}$ for any $p\in \Z$.

(iii) If $z\in \C$ and $\text{Im}(z)$ is sufficiently large, $\exp(zN)F\in D$ (resp. $\exp(zN)F\in D^{\pm}$). 

\medskip
A point $F$ of $D^{\pm}$ is identified with $(\{0\}, \{F\})\in D^{\pm}_{\Sig}$.

Define
$${\bar X}(\C):= G(\Q)\bs (D^{\pm}_{\Sig} \times (G({\bf A}_\Q^f)/K)) = G(\Q)^+ \bs (D_{\Sig} \times G({\bf A}^f_\Q)/K))\supset X(\C).$$

For a representative $(\nu_i)_{1\leq i\leq n}$ of $G(\Q)^+\bs G({\bf A}^f_\Q)/K)$, we have a bijection

(*) $\;\;\coprod_{i =1}^n \Gamma_i \bs D_{\Sig} \overset{\cong}\to \bar X(\C)$

\noindent
extending the bijection (*) in \ref{XK}. By 
 \cite{KU}, \cite{KNU} Part III, $\Gamma_i\bs D_{\Sig}$ has a structure of a {\it  log manifold} in the sense of \cite{KU} 3.5.7. From this, we have a structure of log manifold on $\bar X(\C)$  for which (*) is an isomorphism. This structure is independent of the choice of the representative $(\nu_i)_i$. The log manifold $\bar X(\C)$ contains $X(\C)$ as a dense open subspace, it has the sheaf $\sO_{\bar X(\C)}$ of holomorphic functions which extends the sheaf $\sO_{X(\C)}$ of holomorphic functions of $X(\C)$, and has a log structure whose restriction to $X(\C)$ is trivial. With these structures, $\bar X(\C)$ is a fine moduli space of log mixed Hodge structures.
Concerning this, we will describe in \ref{curve0} more precise things about curves over $\C$.

The universal mixed Hodge structure $\sH_{X(\C)}$ over $X(\C)$  extends to the universal log mixed Hodge structure  $\sH_{\bar X(\C)}$  over $\bar X(\C)$. 
The universal Hodge bundle  $\sH_{X(\C),\sO}:=\sO_{X(\C)} \otimes_\Q \sH_{X(\C),\Q}$ with Hodge filtration  on $X(\C)$ extends canonically to a vector bundle $\sH_{\bar X(\C), \sO}$ on $\bar X(\C)$ with Hodge filtration. 

For each $w\in \Z$, we  have a canonical polarization $\langle \;,\;\rangle_w: \gr^W_w\sH\times \gr^W_w \sH\to \Q(-w)$ on $X(\C)$ which induces the canonical polarization (\ref{canpol}) at each point of $X(\C)$, and this extends to a perfect duality $\langle\;,\;\rangle_w: \gr^W_w \sH_{\sO}\times \gr^W_w\sH_{\sO}\to \sO_{\bar X(\C)}$ on $\bar X(\C)$.

\end{sbpara}

\begin{sbrem} The works \cite{KU} and \cite{KNU} Part III give more general toroidal partial compactifications of $X(\C)$ for which we can use higher dimensional cones of mutually commuting nilpotent operators in $\Lie(G_1)_\R$. The above $\bar X(\C)$ is a special case of this general theory, where we use only cones of dimension $\leq 1$.

  \end{sbrem}

\begin{sbpara}\label{Siegle3} In \ref{Siegel}, if $g=1$, ${\bar X}(\C)$ is the compactification of the modular curve $X(\C)$. If $g\geq 2$, $\bar X(\C)$ is not compact, and it is a toroidal partial compactification of $X(\C)$.

\end{sbpara}

\begin{sbpara}\label{dual}  

We have the sheaf $\Omega_{\bar X(\C)}^1(\log)$ of holomorphic differential forms on $\bar X(\C)$ with (at worst) log poles at the boundary  $\bar X(\C) \smallsetminus X(\C)$. It is a vector bundle on $\bar X(\C)$.

Let 
$$\sE:= \sE nd_{\sO_{\bar X(\C)}, W, \langle\;,\;\rangle}(\sH_{\bar X(\C), \sO}).$$
We have

\medskip

(1)  The vector bundle $T_{\bar X(\C)}:= \sE /F^0\sE$ (here $F^0$ is the Hodge filtration) and
  $\Omega^1_{\bar X(\C)}(\log)$ are canonically the dual vector bundles of each other. 
  
  \medskip
  
  In other words, $T_{\bar X(\C)}$ is  the log tangent bundle of $\bar X(\C)$.

In the pure case, this (1) follows from  \cite{KU} Proposition 4.4.3. The mixed case can be  proved in a similar way.

The restriction of $T_{\bar X(\C)}$ to $X(\C)$ is identified with the tangent bundle of $X(\C)$. 

The vector bundle $T_{\bar X(\C),\hor}: = F^{-1}\sE/F^0\sE\subset T_{\bar X(\C)}$ on $\bar X(\C)$ is called the {\it horizontal log tangent bundle} of $\bar X(\C)$, and its restriction $T_{X(\C), \hor}$ to $X(\C)$ is called the {\it horizontal tangent bundle} of $X(\C)$.
\end{sbpara}

\begin{sbpara}\label{idealI}

Let $I_{\bar X(\C)}$ be the ideal of $\sO_{\bar X(\C)}$ generated by the images in $\sO_{\bar X(\C)}$ of local generators of the log structure. Then $I_{\bar X(\C)}$ is an invertible ideal.  

\end{sbpara}

\begin{sbpara}\label{dif}  Let $\Omega^1_{\bar X(\C)}$ be the $\sO_{\bar X(\C)}$-submodule of $\Omega^1_{\bar X(\C)}(\log)$ generated by $df$ for local sections $f$ of $\sO_{\bar X(\C)}$. Then $\Omega^1_{\bar X(\C)}$ is a vector bundle. 

We have an exact sequence of $\sO_{\bar X(\C)}$-modules 
$$0\to \Omega^1_{\bar X(\C)}\to \Omega^1_{\bar X(\C)}(\log) \to \sO_{\bar X(\C)}/I_{\bar X(\C)}\to 0$$
where the third arrow sends $d\log(q)$ to $1$ for a local generator $q$ of the log structure.

\end{sbpara}

\begin{sbpara}\label{K'}
If $K'$ is an open subgroup of $K$, we have a morphism $\bar X_{K'}(\C)\to \bar X_K(\C)$ of log manifolds. This morphism is log \'etale, and this implies that the pullback of $\Omega^1_{\bar X_K(\C)}(\log)$ on  $\bar X_{K'}(\C)$ coincides with $\Omega^1_{\bar X_{K'}(\C)}(\log)$.

The pullback of the ideal $I_{\bar X_K(\C)}$ on $\bar X_{K'}(\C)$ need not coincide with $I_{\bar X_{K'}(\C)}$. 
The pullback of $\Omega^1_{\bar X_K(\C)}$ (without log) on $\bar X_{K'}(\C)$ need not coincide with $\Omega^1_{\bar X_{K'}(\C)}$. 
\end{sbpara}

 \begin{sbpara}\label{curve0} For a projective smooth curve $C$ over $\C$ and for a finite subset $R$ of $C$, let $C_R$ be $C$ endowed with the log structure associated to the divisor $R$, and let $\Mh(C_R, \bar X(\C))$ be the set of morphisms $C_R\to \bar X(\C)$ of log manifolds which is horizontal. Here ''horizontal'' means that the induced morphism of logarithmic tangent bundles $T_{C_R}\to T_{\bar X(\C)}$ factors through the horizontal log tangent bundle of $\bar X(\C)$ as $T_{C_R} \to T_{\bar X(\C),\hor}\subset T_{\bar X(\C)}$. This set $\Mh(C_R, \bar X(\C))$ is identified with the set of morphisms $C\to \bar X(\C)$ of locally ringed spaces over $\C$ such that the image of $U=C\smallsetminus R$ is contained in $X(\C)$ and such that the morphism $U\to X(\C)$ is horizontal (that is, the morphism of tangent bundles $T_U\to T_{X(\C)}$ factors through $T_{X(\C),\hor}$).

 By
  \cite{KU}, \cite{KNU} Part III, this set $\Mh(C_R, \bar X(\C))$ is identified with the set of isomorphism classes of variations of mixed $\Q$-Hodge structure $\sH$ on $C\smallsetminus R$ 
  satisfying the conditions (i)--(iii) of \ref{Hodgeside} 
  such that 
  
  \medskip
  
   $\text{rank}_{\sO_C}(\gr^r \gr^W_w\sH_{\sO})= h(w,r)$ for all $w,r$
  
  \medskip
 \noindent
  endowed with a polarized $K$-level structure in the following sense.
 A polarized $K$-level structure on $\sH$ is a global section on $C\smallsetminus R$ of the sheaf of equivalence classes of $(\la_1, \theta_1, \la_2, \theta_2)$  
 where $\la_1$ is an isomorphism of sheaves of ${\bf A}^f_\Q$-modules ${\bf A}^f_\Q\otimes_\Q H_{0,\Q}\overset{\cong}\to  {\bf A}^f_\Q\otimes_\Q \sH_\Q$ preserving weight filtrations, $\theta_1$ is an isomorphism of sheaves of $\Q$-vector spaces $\sH_\Q\overset{\cong}\to H_{0,\Q}$ preserving weight filtrations, and $\la_2\in ({\bf A}^f_\Q)^\times$, $\theta_2\in \Q^\times$, satisfying the following condition (*).

\medskip

(*) There exists $\epsilon\in \{\pm 1\}$ such that for any $w\in \Z$, the pairing $$\gr^W_w\sH_\Q\times \gr^W_w\sH_\Q\to \Q \cdot (2\pi i)^{-w} \;;\; (x,y)\mapsto \epsilon^w\langle \theta_1(x), \theta_1(y)\rangle_{0,w}$$ 
is a polarization of $\gr^W_w\sH$. 

\medskip

 The equivalence of  $(\la_1, \theta_1, \la_2, \theta_2)$ and $(\la'_1,\theta'_1, \la'_2, \theta'_2)$ is defined in the same way as in \ref{X(F)}.     
 
   Note that for $\sH\in \Mh(C, \bar X(\C))$ and for $w\in \Z$, we have a canonical polarization  $\langle\;,\;\rangle_w$ on $\gr^W_w\sH$  which is the pull back of the canonical polarization of the universal object on $X(\C)$ (\ref{frak S}) and which 
 induces the canonical polarization \ref{canpol} on its fiber at any point of $C\smallsetminus R$.

  Let 
  $$\Mh(C, \bar X(\C)):= \bigcup_R \; \Mh(C_R, \bar X(\C)).$$
  This set is identified with the set of morphisms $C\to \bar X(\C)$ of locally ringed spaces over $\C$ such that for some finite subset $R$ of $C$, the image of $C\smallsetminus R$ in $\bar X(\C)$ is contained in $X(\C)$ and the morphism $C\smallsetminus R\to X(\C)$ is horizontal.

\end{sbpara}

\subsection{Height functions and period domains}\label{s:2.2b}

\begin{sbpara}\label{2.3.1}
Let $\Phi= ((h(w,r))_{w,r\in\Z}, H_{0,\Q}, W, (\langle\;,\;\rangle_{0,w})_{w\in \Z})$ be as in \ref{Phi} and fix an open compact subgroup  $K$ of $G({\bf A}^f_\Q)$ which satisfies the neat condition \ref{neat2}. 

In the rest of this paper, for a number field $F\subset \C$, we will consider the height functions  
$$H_{\La}, H_{\spadesuit}, H_{\heartsuit,S}: X(F) \to \R_{>0}$$
which are defined fixing somethings and assuming somethings. We will also consider the Hodge theoretic analogues
$$h_{\La}, h_{\spadesuit}, h_{\heartsuit} : \Mh(C, \bar X(\C)) \to \R$$
for a projective smooth curve $C$ over $\C$. 

Here 
$$\La= ((c(w,r))_{w,r\in \Z}, (t(w,d))_{w\in \Z, d\geq 1})\quad \text{with}\;\; c(w,r), t(w,d)\in \R$$
and $S$ is a finite set of places of $F$ which contains all archimedean places of $F$. 

We will consider the relations of these height functions to the extended period domain $\bar X(\C)$, and  describe some philosophies. 
\end{sbpara}

\begin{sbpara}\label{defh}
We first define Hodge theoretic height functions $h_{\La}, h_{\spadesuit}, h_{\heartsuit}: \Mh(C, \bar X(\C))\to \R$. Let $\sH\in \Mh(C, \bar X(\C))$. 

By using the height functions in \S1, define
$$h_{\La}(\sH)= \sum_{w,r} c(w,r)\text{deg}(\gr^r\gr^W_w\sH_{\sO}) + \sum_{w,d} t(w,d)h_{\diamondsuit,w,d}(\sH),$$
 where $h_{\diamondsuit, w,d}$ are defined as in \S1.6 by using 
 the canonical polarizations on $\gr^W_w\sH$ (\ref{canpol}, \ref{curve}). 
 
 $h_{\spadesuit}(\sH)$ is defined as in \S1.6 and is a special case $c(w,r)=b(h,w,r)$ and $t(w,d)=0$ of $h_{\La}(\sH)$ (\ref{ab} (2), \ref{Hanalo}). 
 
Define 
$$h_{\heartsuit}(\sH):=\sum_{x\in C}  h_{\heartsuit, x}(\sH), \quad 
h_{\heartsuit, x}(\sH):=e(x)\in \Z_{\geq 0}$$ where $e(x)$ is as follows. Let $N'_x: \sH_{\Q,x}\to \sH_{Q,x}$ be the local monodromy operator (\ref{Hodgeside}) at $x$. Via the polarized $K$-level structure $\text{class}(\la_1,\theta_1, \la_2, \theta_2)$ of $\sH$, $N'_x$ corresponds to  a nilpotent operator $N:=\la^{-1}_1\circ N'_x\circ\la_1$ on ${\bf A}^f_\Q\otimes_\Q H_{0,\Q}$, which is determined mod the adjoint action of $K$.
If $N'_x=0$, $e(x)=0$.  if $N'_x\neq 0$,  $e(x)\in \Z_{>0}$ is defined by 
$$\{b\in \Q\;|\; (\exp(bN),1)\in K\}=e(x)^{-1}\Z.$$

\end{sbpara}

\begin{sbpara}\label{HX(F)} To define 
$H_{\La}, H_{\spadesuit} : X(F) \to\R_{>0}$, we fix 
 a $\Z$-lattice $H_{0,\Z}$ in $H_{0,\Q}$ such that $H_{0,\hat \Z}:= \hat \Z\otimes_{\Z} H_{0,\Z}$ is stable under the action of $K$ (such $H_{0, \Z}$ exists). 
 
 Let $M\in X(F)$. Concerning the polarized $K$-level structure $\la=\text{class}(\la_1,\theta_1,\la_2,\theta_2)$ of $M$ (\ref{X(F)}), 
 the $\hat \Z$-lattice $\la_1(H_{0,\hat \Z})\subset M_{et}$   is independent of the choice of a representative 
 $(\la_1, \theta_1, \la_2, \theta_2)$ of $\la$ and is 
 stable under the action of $\Gal(\bar K/K)$. 
 Hence $M$ becomes a mixed motive over $F$ with $\hat \Z$-coefficients. By using the height functions in \S1, define
 $$H_{\La}(M):=(\prod_{w,r} H_r(\gr^W_wM)^{c(w,r)})\cdot (\prod_{w,d} H_{\diamondsuit,w, d}(M)^{t(w,d)})$$
where $H_{\diamondsuit,w,d}(M)$ is defined with respect to the canonical polarizations on $\gr^W_wM$ (\ref{canpol}).

Here and in the rest of this paper, when we treat the height function $H_{\La}$, for $(w,d)$ such that $t(w,d)\neq 0$, $\gr^W_wH_{0,\Q}\neq 0$ and $\gr^W_{w-d}H_{0,\Q}\neq 0$, we assume that the conjectures in \S1.7 which we need in the definition of $H_{\diamondsuit, w,d}$ are true.

 $H_{\spadesuit}(M)$ is defined as in \S1.4 and is a special case $c(w,r)=b(h,w,r)$ and $t(w,d)=0$ of $H_{\La}(M)$ (\ref{ab} (2), \ref{rsc} (2)).

\end{sbpara}

\begin{sbpara}\label{Hheart} To define $H_{\heartsuit, S}:X(F)\to \R_{>0}$, we fix 
$$m\geq 1, \;\; (L_i)_{1\leq i\leq m},\;\; (L'_i)_{1\leq i\leq m}, \;\;\text{and}\;\;  J$$
as in (i) below, and we assume that $K$ satisfies (ii) below.

\medskip

(i) $L_i$ and $L'_i$ are $\Z$-lattices in $H_{0,\Q}$ such that $L_i\supset L'_i$ and  $J$ is a subgroup of $\oplus_{i=1}^m (\Z/n_i\Z)^\times$ where $n_i$ is the smallest integer $>0$ which kills $L_i/L'_i$. 
  
  \medskip
  
(ii) For an element $g$ of $\text{Aut}_{{\bf A}^f_\Q}({\bf A}^f_\Q\otimes_\Q H_{0,\Q})$ which belongs to the image of $G({\bf A}^f_\Q)$, $g$ belongs to the image of $K$ if and only if $g\hat L_i=\hat L_i$, $g\hat L'_i=\hat L'_i$ for $1\leq i\leq m$, and for some $h\in J$, the action of $g$ on $L_i/L_i'=\hat L_i/\hat L'_i$ coincides with the scaler action of $h$ via $J\to (\Z/n_i\Z)^\times$ for any $i$. 
Here $\hat L_i=\hat \Z\otimes_\Z L_i$ and $\hat L'_i$ is defined similarly.

   \medskip
   For $M\in X(F)$, define
   $$H_{\heartsuit,S}(M):= \prod_{v\notin S} H_{\heartsuit,v}(M), \quad H_{\heartsuit,v}(M)=\sharp(\F_v)^{e(v)}$$
   where $v$ ranges over all places of $F$ outside $S$ and $e(v)\in \Q_{\geq 0}$ are defined as follows.

    Let $\la=\text{class}(\la_1, \theta_1, \la_2, \theta_2)$ be the polarized $K$-level structure of $M$. Let $T_i=\la_1(\hat L_i)$, $T'_i= \la_1(\hat L'_i)$. Then $T_i$ and $T'_i$ are independent of the choice of the representative $(\la_1,\theta_1,\la_2,\theta_2)$ of $\la$, and are stable under the action of $\Gal(\bar F/F)$. For a prime number $\ell$, let $T_{i,\Z_{\ell}}$ be the $\ell$-adic component of $T_i$ and define $T'_{i,\Z_{\ell}}$ similarly. 
 Let 
 $p=\text{char}(\F_v)$  and let $D=D_{\pst, W(\bar \F_v)}$ (\ref{st}). Let $N'_v$ on $M_{et, \Q_{\ell}}$ for $\ell\neq p$ and $N'_v$ on $D_{\text{pst}}(F_v, M_{et,\Q_p})$ be the local monodromy operators at $v$.

   Let $I$ be the subgroup of $\Q$ consisting of all $b\in \Q$ satisfying the following condition:

For any $\ell\neq p$, the action of $\exp(bN_v')$ on $M_{et, \Q_{\ell}}$ satisfies
$\exp(bN'_v)T_{i,\Z_{\ell}}=T_{i,\Z_{\ell}}$ and $\exp(bN_v')T'_{i,\Z_{\ell}}= T'_{i, \Z_{\ell}}$ for $1\leq i\leq m$, 
the action of $\exp(bN_v')$  on $\hat F_{v,0,\text{ur}}\otimes_{F_{v,0, \text{ur}}} D_{\text{pst}}(F_v,M_{et, \Q_p})$ satisfies $\exp(bN_v')D(T_{i, \Z_p})= D(T_{i, \Z_p})$  and $\exp(bN'_v)D(T'_{i,\Z_p})= D(T'_{i, \Z_p})$ for $1\leq i\leq m$, and there is $h\in J$ such that for $1\leq i\leq m$,  the action of $\exp(bN_v')$ on $T_{i, \Z_{\ell}}/T'_{i, \Z_{\ell}}$ 
coincides with the scaler action of $h$ via $J\to (\Z/n_i\Z)^\times$ for any $\ell\neq p$, 
 and the action of $\exp(bN_v')$ on 
 $D(T_{i,\Z_p})/D(T'_{i,\Z_p})$ coincides with the scaler action of $h$ via $J\to (\Z/n_i\Z)^\times$.

Then $I$ is non-zero. $I=\Q$ if $M$ is of good reduction at $v$.

We define $e(v)=0$ if $I$ is not isomorphic to $\Z$, and we define $e(v)\in \Q_{>0}$ by $I=e(v)^{-1}\Z$ if $I\cong \Z$. 

We expect that the following (*) is true. 

\medskip

(*)  $e(v)\in \Z$. $e(v)=0$ if and only if $M$ is of good reduction at $v$.

 \end{sbpara}

The proof of the following Proposition is easy.

\begin{sbprop} (1) Let $\sH\in \Mh(C, \bar X(\C))$. Then $h_{\heartsuit}(\sH^*)= h_{\heartsuit}(\sH)$ for the dual, $h_{\heartsuit}(\sH(r))=h_{\heartsuit}(\sH)$ for $r\in \Z$ for Tate twists, and $h_{\heartsuit}(\sH')=[C':C]h_{\heartsuit}(\sH)$ if $C'$ is the integral closure of $C$ in a finite extension of the function field of $C$ and $\sH'\in \Mh(C', \bar X(\C))$ is induced by $\sH$. 

(2) Let $M\in X(F)$. Then $H_{\heartsuit,S}(M^*)= H_{\heartsuit,S}(M)$ for the dual, $H_{\heartsuit,S}(M(r))=H_{\heartsuit,S}(M)$ for $r\in \Z$ for Tate twists, and $H_{\heartsuit,S'}(M')=[C':C]H_{\heartsuit,S}(M)$ if $F'$ is  a finite extension of $F$ and $M'\in X(F')$ is the image of  $M\in X(F)$ where $S'$ denotes the set of all places of $F'$ lying over $S$.

 \end{sbprop}
 
 \begin{sbrem}\label{WDconj} If we assume that the following conjecture (i) on $M\in X(F)$ is true, we can define $H_{\heartsuit,S,v}(M)$ without using (i), (ii) in \ref{Hheart}. Let $W'_{F_v}$ be the Weil-Deligne group defined as a group scheme over $\Q$ (\cite{WD} 8.3.6). 
 We expect that the following (i) is true.
 
 \medskip
 
 (i)  Take an embedding $\bar F\to \bar F_v$ over $F$. Then there are a homomorphism $\rho: W'_{F_v}\to G$ and a representative $(\la_1,\theta_1, \la_2, \theta_2)$ 
 of the polarized $K$-level structure of $M$ satisfying the following condition. For any prime number $\ell\neq p$, the homomorphism 
 $$\Gal(\bar F_v/F_v)\to G(\Q_{\ell})\;;\;\sigma\mapsto (\theta_{1,\Q_{\ell}}\circ \sigma_{M_{et, \Q_{\ell}}}
  \circ \theta_{1,\Q_{\ell}}^{-1}, \chi_{\text{cyclo}}(\sigma)^{-1})$$
 corresponds to $\rho$. 
 Here $\sig_{M_{et,\Q_{\ell}}}$ denotes the action of $\sig$ on $M_{et,\Q_{\ell}}$.

 \medskip
 
 Assuming (i), let $N$ be the image of the canonical base of $\Lie(W'_{F_v})$ in $\Lie(G_1)$ under $\rho$. We define
 $H_{\heartsuit, S, v}(M):=e(v)$ where $e(v)=0$ if $N\neq 0$ and $e(v)$ is determined by
 $$\{b\in \Q\;|\; \exp(b N_{\la_1, \theta_1}) \in K\}= e(v)^{-1}\Z$$
 if $N\neq 0$, where $N_{\la_1, \theta_1}=\text{Ad}(g)(N) \in {\bf A}^f_\Q\otimes_\Q \text{Lie}(G_1)$  with $g=({\bf A}^f_\Q\otimes \theta_1)\circ \la_1$. 
 
 \end{sbrem}

We next consider the relations of $h_{\spadesuit}$ and $h_{\heartsuit}$ to the extended period domain $\bar X(\C)$. 

\begin{sbprop}\label{handper} Let $C$ be a projective smooth curve over $\C$, let $\sH\in \Mh(C, \bar X(\C))$, and denote this morphism $C\to \bar X(\C)$ of $\sH$ as $f$. 

(1) $h_{\spadesuit}(\sH)= \text{deg}(f^*\Omega^1_{\bar X(\C)}(\log)).$

(2) $h_{\heartsuit}(\sH)=\sum_{x\in C} h_{\heartsuit,x}(\sH)$ where $h_{\heartsuit,x}(\sH)$ is the integer $e(x)\geq 0$ determined by
$(f^*I_{\bar X(\C)})_x =m_x^{e(x)}$. Here $I_{\bar X(\C)}$ is as in \ref{idealI} and $m_x$ denotes the maximal ideal of $\sO_{C,x}$.

(3)  $h_{\spadesuit}(\sH)-h_{\heartsuit}(\sH)=\text{deg}(f^*\Omega^1_{\bar X(\C)})$. 

\end{sbprop}

\begin{pf}\label{spade0} (1) follows from 
 $$\text{det}_{\sO_C}(\sE nd_{\sO_C, W, \langle\;,\;\rangle}(\sH_{\sO})/F^0)^{-1}\cong f^*\text{det}_{\sO_{\bar X(\C)}}(\sE/F^0\sE)^{-1}$$ 
 $$\cong  f^*\text{det}_{\sO_{\bar X(\C)}}(\Omega^1_{\bar X(\C)}(\log)).\quad(\ref{dual})$$

We prove (2). Let $R$ be a finite subset of $C$ such that $\sH$ is non-degenerate outside $R$, 
let $\sM$ be the log structure of $C$ associated to the divisor $R$, and let $\sM_{\bar X(\C)}$ be the log structure of $\bar X(\C)$. Let $x\in R$, let $y$ be the image of $x$ in $\bar X(\C)$. Then by \cite{KU} and \cite{KNU} Part III, $\exp(e(x)^{-1}\Z N'_x)$ is identified with the local monodromy group of $\bar X(\C)$ at $y$ and the inclusion map $\exp(\Z N'_x) \to \exp(e(x)^{-1} \Z N'_x)$ is identified with the canonical map from the local monodromy group of $C$ at $x$ to the local monodromy group of $\bar X(\C)$ at $y$. 
By \cite{KU} and \cite{KNU} Part III, we have a commutative diagram 
$$\begin{matrix}  
\exp(\Z N'_x) &\overset{\cong}\to  & \Hom(\sM^{\gp}_x/\sO^\times_{C,x}, \Z(1))\\
\downarrow &&\downarrow\\
\exp(e(x)^{-1}\Z N'_x) & \overset{\cong}\to & \Hom(\sM^{\gp}_{\bar X(\C), y}/\sO^\times_{\bar X(\C),y}, \Z(1))\\
\end{matrix}$$
where $\sS^{\gp}$ for a commutative monoid $\sS$ means the associated group $\{ab^{-1}\;|\; a,b \in \sS\}$, the horizontal isomorphism are as in \cite{KU} 2.2.9, and the right vertical arrow is given by the pullback map $\sM_{\bar X(\C),y}\to \sM_x$. This shows that for a generator $q$ of the log structure of $\bar X(\C)$ at $y$, the pull back of $q$ at $x$ generates $m_x^{e(x)}$. 

(3) follows from (1) and (2) and the exact sequence in \ref{dif}. 
\end{pf}

\begin{sbpara}\label{Cinfmet}  
 In \ref{Cinfmet}--\ref{spadeper}, we consider a canonical metric on $\text{det}_{\sO_{X(\C)}}(\Omega^1_{X(\C)})$, a canonical measure on $X(\C)$, and their relations to the height function $H_{\spadesuit}$.

 Let $T_{X(\C)}= \sE nd_{\sO_{X(\C)}, W, \langle\;,\;\rangle}(\sH_{X(\C), \sO})/F^0$ be the tangent bundle of $X(\C)$ (\ref{dual}).  The canonical Hodge metrics of $\gr^W_w\sH_{X(\C)}$ ($w\in \Z$) 
 give a Hodge metric on $\text{det}_{\sO_{X(\C)}}(T_{X(\C)})$. By the duality \ref{dual}, this gives a canonical Hermitian $C^{\infty}$ metric
 on  $\text{det}_{\sO_{X(\C)}}(\Omega^1_{X(\C)})$.
  (Here a metric $|\;\;|$ of a line bundle on a complex analytic manifold is said to be $C^{\infty}$ if for a local base $e$ of the line bundle, $|e|$ is a $C^{\infty}$-function.) It gives a canonical measure on $X(\C)$. 

This canonical measure on $X(\C)$   is actually a classical object. The pull back of this measure on $Y:= D^{\pm} \times G({\bf A}^f_\Q)/K$ via the local homeomorphism $Y\to X(\C)$  is a $G(\R)G_u(\C)\times G({\bf A}^f_\Q)$-invariant measure on $Y$ (which is unique upto $\R_{>0}$-multiple). Hence it is a classical object.

By using the theory of degenerations of Hodge metrics as in \cite{CKS}, \cite{KNU} Part II, etc., we can show that this canonical $C^{\infty}$ metric on $\text{det}_{\sO_{X(\C)}}(\Omega^1_{X(\C)})$ extends to a metic of $\text{det}_{\sO_{\bar X(\C)}}(\Omega^1_{\bar X(\C)}(\log))$ with (at worst) log singularities at the boundary. 
\end{sbpara}

\begin{sbpara}\label{poincare}

{\bf Example.} Consider the case $g=1$ of \ref{Siegel}. In this case,  $D$ is canonically identified with 
the upper half plane $\frak H$: $z\in \frak H$ corresponds to $F\in D$ defined by $$F^0=H_{0, \C}\supset F^1=\C\cdot (ze_1+e_2)\supset F^2=0.$$
On $D$, the canonical metric $|\;\;|$ on $\text{det}_{\sO_{X(\C)}}(\Omega^1_{X(\C)})$ is characterized by the property that $|(2y)^{-1} dz|=1$ where $y= \text{Im}(z)$. The corresponding measure of $D$ is the Poincar\'e measure $y^{-2}dxdy$. 
In the  toroidal compactification $\bar X(\C)$ of $X(\C)$, at the standard cusp (the limit of the image of $x+iy\in \frak H$ in $\bar X(\C)$ with $y\to \infty$), if $q$ denotes  $\exp(2\pi i z)$ (so $q\to 0$ at this boundary),  
$dq/q$ is a basis of $\Omega^1_{\bar X(\C)}(\log)=\text{det}_{\bar X(\C)}(\Omega^1_{\bar X(\C)}(\log))$ at this cusp, and $|(4 \pi y)^{-1}dq/q|=1$. Thus the metric on $\Omega^1_{\bar X(\C)}(\log)$ has 
a log singularity $4\pi y$ at this cusp.

\end{sbpara}

\begin{sbpara}\label{spadeper} The height function 
$H_{\spadesuit}$ is related to  period domain $X(\C)$ as follows. For a number field $F\subset \C$ and for 
$M\in X(F)$, and for the archimedean place $v$ of $F$ corresponding to this embedding $F\to \C$, 
the canonical metric $|\;\;|_v$ on $\text{det}_F(\text{End}_{F, W,\langle\;,\;\rangle}M_{\dR})/F^0)^{-1}$ in \S1.4 used in the definition of $H_{\spadesuit}$  is induced by the canonical metric on $\text{det}_{\sO_{X(\C)}}(\Omega^1_{X(\C)})$ which corresponds to the canonical measure on $X(\C)$. 

\end{sbpara}

\begin{sbrem} The author plans to discuss  non-archimedean versions of these stories elsewhere by using $p$-adic period domains (\cite{Ra}, \cite{Ka1}) and using toroidal partial compactifications of $p$-adic period domains in \cite{Ka1}. For $F$ and $M$ as in \ref{spadeper} and for a non-archimedean place $v$ of $F$, the canonical metric $|\;\;|_v$ on 
$\text{det}_F(\text{End}_{F, W,\langle\;,\;\rangle}M_{\dR})/F^0)^{-1}$ in \S1.4 also comes from a metric on 
$\text{det}_{\sO_{X(F_v)}}(\Omega^1_{X(F_v)})$ for a $v$-adic period domain $X(F_v)$ which extends to a metric of 
$\text{det}_{\sO_{\bar X(F_v)}}(\Omega^1_{\bar X(F_v)}(\log))$ on a toroidal partial compactification $\bar X(F_v)$ of $X(F_v)$. 

\end{sbrem}

\begin{sbpara}\label{X?} We think that $X(F)\to X(\C)\to \bar X(\C)$ is similar to $V(F)\subset V(\C)\subset \bar V(\C)$ in (2) of 0.1.

From the above observations on period domains and their toroidal partial compactifications, we think that in the analogy between (1), (2) and (3) in \S0, the height functions $H_{\spadesuit}$ in (1) and $h_{\spadesuit}$ in (3) are  like the height function $H_{K+D}$ in (2), and the height function $h_{\heartsuit}$ in (3) is like the height function $h_D$ in (2). 
 $H_{\heartsuit,S}(M)$ is like the height function, which lacks the part of places in $S$,  associated to the boundary of the moduli space of motives. Roughly speaking, 
$H_{\heartsuit,S}(M)$  is like the intersection number of $M$ and the boundary of the moduli space of motives (outside $S$).

\end{sbpara}

\begin{sbpara}\label{Dmu20} Like $H_{\spadesuit}$, the height function  $H_{\clubsuit}$ in \S1 is also related to differential forms and measures on some period domain.

Fix an integer $w$ and integers $h(r)\geq 0$ ($r\in \Z$) such that $h(r)=0$ for almost all $r$ and $h(w-r)=h(r)$ for all $r$. We fix an $\R$-vector space $H_{0,\R}$ of dimension $\sum_r h(r)$. 
Let $H_{0,\C}:=\C\otimes_\R H_{0,\R}$.
Let $\check{D_{\clubsuit}}$ be the set of all decreasing filtrations on $H_{0,\C}$ satisfying 
 $\dim_{\C}\gr^r H_{0,\C}= h(r)$ for any $r\in \Z$.
Let $D_{\clubsuit}\subset \check{D_{\clubsuit}}$ be the subset consisting of all $F\in \check{D_{\clubsuit}}$ such that $(H_{0,\R}, F)$ is a pure $\R$-Hodge structure of weight $w$. 
Then $D_{\clubsuit}$ is an open set of the complex analytic manifold $\check{D_{\clubsuit}}$ and hence is also a complex analytic manifold. 
If $\sH_{D_{\clubsuit},\sO}$ denotes the universal Hodge bundle on $D_{\clubsuit}$, $T_{D_{\clubsuit}}:=\text{End}_{\sO_{D_{\clubsuit}}}(\sH_{D_{\clubsuit},\sO})/F^0$ and $\Omega_{D_{\clubsuit}}^1$ are canonically dual vector bundles of each other.  By the method of \S1.4, we have a canonical $C^{\infty}$ metric on $\text{det}_{\sO_{D_{\clubsuit}}}(T_{D_{\clubsuit}})^{-1}\cong \text{det}_{\sO_{D_{\clubsuit}}}(\Omega^1_{D_{\clubsuit}})$, whose fibers give the metrics at archimedean places used in the definition of $H_{\clubsuit}$ in \S1.4.  It corresponds to a canonical measure on $D_{\clubsuit}$. 

 (But we do not use $H_{\clubsuit}$ in the rest of this paper.)

\end{sbpara}

\subsection{Speculations on curvatures and motives}\label{cumo}\label{2.3b}

In \ref{curve}--\ref{pos}, we  review a work of Griffiths (\ref{Gcurva}) on curvature forms on Griffiths domains, works of Griffiths and Peters on variations of Hodge structure on curves, and give complements.

In \ref{Q1Q2}--\ref{mix7}, 
we present questions about the relation of the curvature forms  on the Griffiths domains and motives over number fields.

\begin{sbpara} 

Let $\Phi= ((h(w,r))_{w,r\in\Z}, H_{0,\Q}, W, (\langle\;,\;\rangle_{0,w})_{w\in \Z})$ be as in \ref{Phi} and fix an open compact subgroup  $K$ of $G({\bf A}^f_\Q)$ which satisfies the neat condition \ref{neat2}.

Except in \ref{mix100}--\ref{GA1}, we assume that $\Phi$ is pure of weight $w$ (that is, $W_w=H_{0,\Q}$ and $W_{w-1}=0$). So, $h(w',r)=0$ if $w'\neq w$. Let $h(r):=h(w,r)$.

\end{sbpara}

\begin{sbpara}\label{curve} Let $Z$ be a complex analytic manifold. 
 For a holomorphic line bundle $L$ on $Z$ with a $C^{\infty}$ metric $|\;\;|$ (\ref{Cinfmet}),  the associated curvature form $\kappa$ (called also the Chern form), which is a $(1,1)$-form on $Z$,  is defined as follows. Let $e$ be a local base of $L$, and let $\kappa= \bar \partial \partial \log(|e|)$. Then this is independent of $e$, and defined globally. 
 
 It gives an Hermitian form on the tangent bundle $T_Z$ of $Z$.

\end{sbpara}

\begin{sbpara}\label{kappa}
Let $\kappa_r$ for $r\in \Z$ be the restriction to $T_{X(\C),\hor}$ of the curvature form of the canonical Hodge metric of $\text{det}_{\sO_{X(\C)}}(\gr^r\sH_{X(\C),\sO})$.

We have $$\kappa_{w-r}= -\kappa_r$$ since $\gr^r\sH_{X(\C),\sO}$ and $\gr^{w-r}\sH_{X(\C),\sO}$ are dual vector bundles with dual metrics by the polarization. Hence $2\kappa_r= \kappa_r-\kappa_{w-r}$, and this $\kappa_r-\kappa_{w-r}$ 
is the curvature form of the canonical metric on $\text{det}_{\sO_{X(\C)}}(\gr^r\sH_{X(\C),\sO})\otimes_{\sO_{X(\C)}} \text{det}_{\sO_{X(\C)}}(\gr^{w-r}\sH_{X(\C),\sO})^{-1}$ defined in 1.4.1 without using polarization. Hence $\kappa_r$ can be defined without using the polarization of $\sH_{X(\C)}$.

\end{sbpara}

\begin{sbpara}\label{geqr}

For $r\in \Z$, let
$$\kappa_{\geq r}:= \sum_{i\geq r} \kappa_i$$
This is the restriction to $T_{X(\C),\hor}$ of the curvature form of the Hodge metric of $\text{det}_{\sO_{X(\C)}}(F^r\sH_{X(\C),\sO})$.

From  $\kappa_{w-r}=-\kappa_r$ ($r\in \Z$), we have $$\kappa_{\geq r}=\kappa_{\geq w+1-r}\quad \text{for any}\; r, \quad\kappa_{\geq r}=0 \quad \text{for almost all}\;r.$$

\end{sbpara}

\begin{sbpara}\label{La5} For $\La=(c(r))_{r\in \Z}$ with $c(r)\in \R$, let
$$\kappa_{\La}:= \sum_r c(r) \kappa_r.$$

The pullback of $\kappa_{\La}$ to $D^{\pm}\times G({\bf A}^f_\Q)/K$ is $G(\R) \times G({\bf A}^f_\Q)$-invariant. 
Since this group acts transitively on $D^{\pm}\times G({\bf A}^f_\Q)/K$, 
  $\kappa_{\La}$ is determined by its fiber at any one point of $D$. In particular, $\kappa_{\La}$ is positive (resp. positive semi-definite)  if and only if its pullback to $T_{D,\hor}$ is positive definite (resp. positive semi-definite).

We consider the height functions $h_{\La}$ and $H_{\La}$ of \S2.3 in the present pure situation. 
For a projective smooth curve $C$ and for $\sH\in \Mh(C, \bar X(\C))$, let $h_{\La}(\sH)= \sum_r c(r) \text{deg}(\gr^r\sH_{\sO})$. Fixing a $\Z$-lattice $H_{0,\Z}$ in $H_{0,\Q}$ such that $H_{0, \hat \Z}= \hat \Z\otimes_\Z H_{0,\Z}$ is stable under the action of $K$, for a number field $F\subset \C$ and for $M\in X(F)$, let $H_{\La}(M) = \prod_r H_r(M)^{c(r)}$. 

\end{sbpara}

\begin{sbpara}\label{Gcurva} 
We have a canonical isomorphism $$T_{X(\C),\hor}\cong  \{(f_r)_r\in \oplus_{r\in \Z} \; \sH om_{\sO_{X(\C)}}(\gr^r\sH_{X(\C), \sO}, \gr^{r-1}\sH_{X(\C),\sO}) \;|\; f_{w+1-r}=-{}^tf_r\;(\forall r)\}.$$
Here ${}^tf_r$ denotes the transpose of $f_r$ with respect to the canonical polarization $\langle\;,\;\rangle_w$. For $r\in \Z$, we will denote by $T_{X(\C),\hor,r}$ the direct summand $\{(f_i)_i\in T_{X(\C),r}\;|\; f_i=0 \;\text{for}\; i\neq r,w+1-r\}$. Then $T_{X(\C),\hor}$ is the sum of $T_{X(\C), \hor,r}$ for $r\in \Z$. We have $T_{X(\C), \hor, r}=T_{X(\C), \hor, w+1-r}$.

The following was proved by Griffiths (\cite{Gri}).

\medskip

Let $r\in \Z$. Then
the Hermitian form $\kappa_{\geq r}$ on the horizontal tangent bundle $T_{X(\C),\hor}$ is positive semi-definite. Its restriction to  the component  $T_{X(\C),\hor,r}$ of $T_{X(\C),\hor}$  (\ref{hor})  is  positive definite.

In particular, the Hermitian form $\sum_{r\in \Z} \kappa_{\geq r} = \sum_{r\in \Z} r\kappa_r$ on $T_{X(\C),\hor}$ is positive definite. 

\end{sbpara}

\begin{sbpara}\label{Higgs} Here we review roughly how Griffiths (\cite{Gri}) and Peters (\cite{Pe}) proved 
the case $\sH$ is pure of Proposition \ref{Hprop}.

 Let $C$ be a projective smooth curve over $\C$ and let $\sH\in \Mh(C, \bar X(\C))$.  
 Let $\La=(c(r))_{r\in \Z}$, $\kappa_{\La}$ and $h_{\La}$ be as in \ref{La5}. We have:
\medskip

(1) Assume $\kappa_{\La}$ is positive semi-definite as an Hermitian form on $T_{X(\C),\hor}$. Then 
$h_{\La}(\sH)\geq 0.$

(2) Assume $\kappa_{\La}$  is positive definite as an Hermitian form on $T_{X(\C),\hor}$ and assume $h_{\La}(\sH)= 0$.  
Then $\sH$ is constant. That is, the map $C\to \bar X(\C)$ is a constant map.

\medskip

By \ref{Gcurva}, taking $c(r)=r$, this  (1) (resp. (2)) proves the case $\sH$ is pure of  \ref{Hprop} (1) (resp. (2)).

The above (1) and (2) are proved, roughly speaking, as follows. Write the morphism $C \to \bar X(\C)$ of  $\sH$ as $f$. 
Take a finite subset $R$ of $C$ such that $f$ induces $C\smallsetminus R\to X(\C)$. The curvature form of $\text{det}_{\sO_C}\gr^r H_{\sO}$ on $C\smallsetminus R$ coincides with the 
pullback of $\kappa_r$ under $f$.
On the other hand,

(3) $\text{deg}(\gr^r \sH_{\sO})$ is the integration of the curvature form of $\text{det}_{\sO_C}\gr^r_F H_{\sO}$  over the curve $C$.

(1) follows from this. If the assumptions in (2) are satisfied, the map $T_C\to T_{X(\C),\hor}$ induced by $f$ is trivial, and hence $f$ is a constant map.

\end{sbpara}

\begin{sbpara} Let
$$\kappa_{\star} :=\sum_{r\in \Z} \kappa_{\geq r}, \quad  \kappa_{\spadesuit}:=\sum_{r\in \Z} b(h, w,r)\kappa_r$$
where $b(h,w,r)$ is as in \ref{ab} (2).

\end{sbpara}

\begin{sblem}\label{HHspade2}

 Let $S= \{r\in \Z\;|\; h(r) \neq 0\}$ and let $S^*\subset S$ be the complement of the smallest element of $S$. For $s\in S^*$, let $s'$ be the largest element of $S$ such that $s>s'$.

(1)  
$$\kappa_{\star}= \sum_{s\in S^*} (s-s') \kappa_{\geq s}.$$

(2) Let $m$ be the minimal element of $S^*$ such that $m>w/2$ where $w$ is the weight. Then
$$\kappa_{\spadesuit}
=(\sum_{s\in S^*} h(s)\kappa_{\geq s})-(-1)^w \kappa_{\geq m}.$$

\end{sblem}

\begin{pf} This can be obtained by elementary computation. 
\end{pf}

Recall that $\kappa_{\star}$ is positive definite (\ref{Gcurva}). 
 
\begin{sbprop}\label{pos}  $\kappa_{\spadesuit}$ is positive definite.

\end{sbprop}

\begin{pf} If $w$ is odd, this follows from (2)  of \ref{HHspade2}. Assume $w$ is even. If $h(m-1)=0$, the component $T_{D,\hor,m}$ of $T_{D,\hor}$ (\ref{Gcurva})
is trivial and hence by (2) of \ref{HHspade2} and by the result \ref{Gcurva} of Griffiths, $\kappa_{\spadesuit}$ is positive definite. Assume $h(m-1)\neq 0$. Then by the definition of $m$, we have $m-1=w/2$ and $w/2\in S^*$. Since 
 $\kappa_{\geq m}= \kappa_{\geq w/2}$ (\ref{geqr}), $$\kappa_{\spadesuit}= (\sum_{s\in S^*\smallsetminus \{m,w/2\}} h(s)\kappa_{\geq s}) + \frac{1}{2}h(w/2) \kappa_{\geq m} + \frac{1}{2}h(w/2) \kappa_{\geq w/2}.$$
By the result of Griffiths \ref{Gcurva}, this proves that $\kappa_{\spadesuit}$ is positive definite.
 \end{pf}

As Koshikawa shows in \cite{Ko1}, \cite{Ko2}, we can often deduce properties of a height function for motives from corresponding properties of the Hodge analogues. 
The above Hodge analogues  lead to the following questions.

Let $\La=(c(r))_{r\in \Z}$ with $c(r)\in \R$. Let the notation be as in \ref{La5}.

\begin{sbpara}\label{Q1Q2} 

{\bf Question 1.} Are the following three conditions (i)--(iii) equivalent?

\medskip

(i) The Hermitian form $\kappa_{\La}$ (\ref{La5}) on the horizontal tangent bundle $T_{X(\C),\hor}$ of $X(\C)$ is 
positive semi-definite.

\medskip

(ii)  Fix a $\Z$-lattice $H_{0,\Z}$ in $H_{0,\Q}$ such that $H_{0,\hat \Z}$ is stable under the action of $K$. Then there exists $C\in \R_{>0}$ such that 
$$H_{\La}(M) \geq C^{[F:\Q]}$$
for any number field $F\subset \C$ and for any  $M\in X(F)$.

\medskip
(iii) $h_{\La}(\sH)\geq 0$ for any projective smooth curve $C$ over $\C$ and for any $\sH\in \Mh(C, \bar X(\C))$. 
\medskip

By the works of Griffiths and Peters,  (i) $\Rightarrow$ (iii) is true. 

\medskip

{\bf Question 2.} Are the following (i)--(iii) equivalent?

\medskip

(i) The Hermitian form $\kappa_{\La}$ on the horizontal tangent bundle $T_{X(\C),\hor}$ of $X(\C)$ is 
positive definite.

\medskip

(ii) If we fix a $\Z$-lattice $H_{0,\Z}$ in $H_{0,\Q}$ such that $H_{0,\hat \Z}$ is stable under the action of $K$ and fix an integer $d\geq 1$ and $C>\R_{>0}$, then  there exist finitely many pairs $(F_i, M_i)$ ($1\leq i\leq n$) of number fields $F_i\subset \C$ and $M_i\in X(F_i)$ having the following property: For any number field $F\subset \C$ and any $M\in X(F)$ such that $[F:\Q]\leq d$ and $H_{\La}(M)\leq C$, there is $i$ ($1\leq i\leq n$) such that $F_i\subset F$ and $M$ is the image of $M_i$ under $X(F_i)\to X(F)$.

\medskip

(iii) If $C$ is a projective smooth curve over $\C$ and if $\sH\in \Mh(C, \bar X(\C))$ and $h_{\La}(\sH)\leq 0$, then $\sH$ is constant.

\medskip

By the works of Griffiths and Peters, the implication (i) $\Rightarrow$ (iii) is true. 
\end{sbpara}

\begin{sbrem} The above Questions 1 and 2 are based on the analogous known results (1) and (2) on usual height functions, respectively. Let $V$ and $\bar V$ be as in (2) in 0.1. 

(1) Let $E$ be an effective divisor on $\bar V$ and let $H$ be a height function associated to $E$. Then there is a constant $C\in \R_{>0}$ such that $H(x)\geq C^{[F':F]}$ for any finite extension $F'$ of $F$ and for any $x\in V(F')$. 

(2) Let $A$ be an ample line bundle on $\bar V$ and let $H$ be an associated height function. Let $d\geq 1$ and fix $C\in \R_{>0}$. Let $\bar F$ be an algebraic closure of $F$. Then there is a finite subset $I$ of  $ V(\bar F)$ 
having the following property: For any number field $F'\subset \bar F$ and any $x\in V(F')$ such that $[F':F] \leq d$ and $H(x) \leq C$, we have $x\in I$.

\end{sbrem} 

If the answer to the part (i) $\Rightarrow$ (ii) of the above Question 2 is Yes, it gives a very strong consequence in number theory as in \ref{mform} below. 

\begin{sblem}\label{hor0} Assume $T_{X(\C),\hor}=0$. Assume that the implication (i) $\Rightarrow$ (ii) in Question 2 is true. Fix $d\geq 1$. Then
there exist finitely many pairs $(F_i, M_i)$ ($1\leq i\leq n$) of number fields $F_i\subset \C$ and $M_i\in X(F_i)$ having the following property: For any number field $F\subset \C$ and any $M\in X(F)$ such that $[F:\Q]\leq d$, there is $i$ ($1\leq i\leq n$) such that $F_i\subset F$ and $M$ is the image of $M_i$ under $X(F_i)\to X(F)$.

\end{sblem}

\begin{pf} Let $\La=(c(r))_{r\in \Z}$ with $c(r)=0$ for all $0$. Then $H_{\La}(M)=1$ for any $(F, M)$. Since $T_{X(\C),\hor}=0$, $\kappa_{\La}$ is positive definite.
\end{pf}

\begin{sblem}\label{mform} Assume that the implication (i) $\Rightarrow$ (ii) in Question 2 is true. Let $k\geq 4$ be an even integer. Then modulo quadratic twists,  there are only finitely many new forms of weight $k$ (in the usual sense) without complex multiplication whose Fourier expansions are with $\Q$-coefficients. 
\end{sblem}

(We consider here only even $k$ because any new form of odd weight whose Fourier expansion is with $\Q$-coefficients has complex multiplication by Ribet \cite{Ri}.)

\begin{pf} 
Let $\Phi$ be as follows: $h(w,r)=1$ for $(w,r)=(1,0), (1,k-1)$ and $h(w,r)=0$ for all other $(w,r)$. $H_{0,\Q}=\Q^2$ with the standard basis $(e_1, e_2)$. $W$ is pure of weight $k-1$. $\langle\;,\;\rangle_{0,1}$ is defined by $\langle e_2, e_1\rangle_{0,k-1}=(2\pi i)^{1-k}$, $\langle e_1,e_2\rangle_{0,k-1}=-(2\pi i)^{1-k}$, $\langle e_i,e_i\rangle_{0,k-1}=0$ ($i=1,2$). Fix $K$ which satisfies the neat condition. In this case, $T_{X(\C),\hor}=0$.

Let $A$ be the set of all such new forms. For $f\in A$, let $M(f)$ be the pure motive with $\Q$-coefficients of weight $k-1$ over $\Q$ associated to $f$. For $d\geq 1$, let $B(d)$ be the set of all triples   $(f, F, \la)$ where $f\in A$, $F\subset \C$ is a number field such that $[F:\Q]\leq d$, and $\la$ is a polarized $K$-level structure on $M(f)|_F$ where $M(f)|_F$ denotes the motive over $F$ induced by $M(f)$. Then if $d$ is sufficiently large, the map $B(d) \to A\;;\;(f, F, \la)\mapsto f$ is surjective. 
By \ref{hor0}, we have $(F_i, M_i)$ ($1\leq i\leq n$) having the property in \ref{hor0} for this $d$. 
Let $C(d)$ be the set of all  quadruples $(f, F, \la, i)$ where $(f, F, \la)\in B(d)$ and $i\in \{1,\dots, n\}$ such that $F_i\subset F$ and $(M(f)_F,\la)\in X(F)$ is the image of $M_i\in X(F_i)$ under the map $X(F_i)\to X(F)$. Then the map $C(d)\to B(d)\;;\;(f,F, \la, i)\mapsto (f,F, \la)$ is surjective. Let $I\subset \{1, \dots, n\}$ be the image of the map $C(d)\to \{1,\dots, n\}\;;\;(f, F, \la, i)\mapsto i$. For each $i\in I$, choose an element $(f_i, F'_i, \la'_i, i)$ of $C(d)$. Then if $f\in A$, there is an element $(f, F, \la, i)$ of $C(d)$, and $M(f)|_{F'}\cong M(f_i)|_{F'}$ for any number field $F'$ such that $F\subset F'$ and $F_i'\subset F'$. From $M(f)|_{F'}\cong M(f_i)|_{F'}$, we obtain that either $f= f_i$ or $f$ is a quadratic twist of $f_i$. 
(The last thing is obtained by looking at the representations $M(f)_{et,\Q_p}$ and $M(f_i)_{et, \Q_p}$ of $\Gal(\bar \Q/\Q)$ 
($p$ is a prime number) and by using the fact that $f\in A$ is determined by the representation $M(f)_{et,\Q_p}$ of $\Gal(\bar \Q/\Q)$.) 
\end{pf}

\begin{sbrem}\label{Robert}  
The work of Robert \cite{Ro} suggests that the surprising finiteness in \ref{mform} may be true. 

The author thanks T. Koshikawa for informing the importance of this paper of Robert for the subject of this paper. 

\end{sbrem}

\begin{sbpara}\label{mix100} In \ref{mix100}--\ref{GA1}, 
we consider mixed motives (we do not assume $\Phi$ is pure).

For each $w\in \Z$, let $\Phi_w= ((h'(w',r))_{w',r\in \Z}, \gr^W_wH_{0,\Q}, W',  \langle\;,\;\rangle_{0,w})$ where $h'(w', r)=h(w,r)$ if $w'=w$ and $h'(w',r)=0$ if $w'\neq w$ and $W'$ is pure of weight $w$. Let $D_w$ be Griffiths domain $D$ of  the pure situation $\Phi_w$. 
For $w,r\in \Z$, let  the Hermitian form $\kappa_{w,r}$ on $T_{D_w,\hor}$ be the $\kappa_r$ of \ref{kappa} of the pure situation $\Phi_w$.

Let  $\La=((c(w,r)_{w,r\in \Z}, (t(w,d))_{d\geq 1})$ ($c(w,r), t(w,d)\in \R$) as in \ref{2.3.1}.

We say $\La$ is {\it ample} if $\sum_r c(w,r)\kappa_{w,r}$ is positive definite for any $w\in \Z$ and if $t(w,d)>0$ for any $w,d$.

For example, if we take $c(w,r)=r$ for all $w,r$ and $t(w,d)=1$ for all $w,d$,  $\La$ is ample and 
$$H_{\La}= H_{\star\diamondsuit}, \quad h_{\La}=h_{\star\diamondsuit}.$$

\end{sbpara}

\begin{sbpara}\label{mix7}  We consider the mixed motive version of the above Question 2, which is stronger than the finiteness statement considered in \S2.1. (The number field $F$ was fixed in \S2.1 but $F$ moves below.) 
Fix $H_{0,\Z}$ as in \ref{HX(F)} to define $H_{\La}$. 

\medskip

{\bf Question 3.} Is the following true?

Assume $\La$ is ample (\ref{mix100}). 
Fix $d\geq 1$ and $C\in \R_{>0}$. Then
there exist finitely many pairs $(F_i, M_i)$ ($1\leq i\leq n$) of number fields $F_i\subset \C$ and $M_i\in X(F_i)$ having the following property: For any number field $F\subset \C$ and any $M\in X(F)$ such that $[F:\Q]\leq d$ and $H_{\La}(M)\leq C$, there is $i$ ($1\leq i\leq n$) such that $F_i\subset F$ and $M$ is the image of $M_i$ under $X(F_i)\to X(F)$.

\end{sbpara}

 \begin{sbrem}\label{GA1} The arithmetic height function $H_{\La}$ and the Hodge theoretic height function  $h_{\La}$ are related in the theory of asymptotic behaviors of $H_{\La}$, in the following way, supporting the speculations in this \S\ref{2.3b}.  The details will be discussed elsewhere. 
 
 Let $F\subset \C$ be a number field and let $C$ be a projective smooth curve over $F$. Let $M\in X(F(C))$ where $F(C)$ is the function field of $C$ and $X(F(C))$ is defined in the similar way as $X(F)$.

 T. Koshikawa, S. Bloch and the author proved that the following formula holds in certain cases:
 
 \medskip
 
 (*) $\;\;\log(H_{\La}(M(x)))/\log(H(x))\to h_{\La}(\sH)\;\;$ when $\;\;H(x)^{1/[F':F]} \to \infty$. 
 
 \medskip
 
Here $F'$ ranges over finite extensions of $F$, $x$ ranges over $F'$-rational points of $C$ at which $M$ does not degenerate, $M(x)\in X(F')$ is the specialization of $M$ at $x$, $\sH\in \Mh(C, \bar X(\C))$ is the Hodge realization of $M$, and $H(x)$ is the height of the $F'$-rational point $x$ of $C$ defined as $H_{L}(x)^{1//\text{deg}(L)}$, where $L$ is a fixed ample line bundle  on $C$ and $H_L$ is a height function associated to $L$. In the case $\Phi$ is pure, this (*) follows from Koshikawa  \cite{Ko2} \S4.3. 

In  (*),  the height function $H_{\La}$ for (1) in \S0, the height function $h_{\La}$ for (3) in \S0, and the height function $H$ for (2) in \S0 are connected actually, not only philosophically. 

For example, let $E$ be an elliptic curve over $F(C)$, let $a\in E(F(C))$, and let $M$ be the mixed motive over $F(C)$ corresponding to $a$ such that $W_0M=M$, 
$\gr^W_0=\Q$, $\gr^W_{-1}$ is the $H_1$ of $E$, and $W_{-2}M=0$. Then $H_{\diamondsuit, 0,1}(M(x))$ coincides with the N\'eron-Tate height of the specialization $a(x)$ of $a$ at $x$ and $h_{\diamondsuit, 0,1}(\sH)$ coincides with the geometric height of $a$ defined by using the intersection theory. The above formula for $H_{\diamondsuit, 0,1}$ in this case was proved by Tate in  \cite{Ta}. A generalization to abelian varieties was obtained in Green \cite{Gre}.

\end{sbrem}

\subsection{Speculations on Vojta conjectures}\label{2.4b}

In this \S\ref{2.4b}, we consider how to formulate Vojta conjectures for motives. In \ref{Sha}, we give a related comment on  Shafarevich conjecture for motives.

Let $\Phi= ((h(w,r))_{w,r\in\Z}, H_{0,\Q}, W, (\langle\;,\;\rangle_{0,w})_{w\in \Z})$ be as in \ref{Phi} and fix an open compact subgroup  $K$ of $G({\bf A}^f_\Q)$ which satisfies the neat condition \ref{neat2}.

\begin{sbpara}\label{generic} Let $M\in X(F)$.
We say $M$ is generic if
the following (*) is satisfied for some prime number $p$. Let $\text{class}(\la_1,\theta_1,\la_2, \theta_2)$ be the polarized $K$-level structure of $M$. 
Let $\chi_{\text{cyclo}}: \Gal(\bar F/F)\to \Q_p^\times$ be the cyclotomic character. 

\medskip

(*) There is no algebraic subgroup $H$ of $G$ over $\Q$ such that $\dim(H)<\dim(G)$ and such that the image of
$$\Gal(\bar F/F)\to G(\Q_p)\;;\; \sigma \mapsto (\theta_{1,\Q_p}\circ \sigma \circ\theta^{-1}_{1,\Q_p}, \;\chi_{\text{cyclo}}(\sig)^{-1})$$
is contained in $H(\Q_p)$. 

\medskip

The philosophy of Mumford-Tate groups  predicts that the condition (*) is independent of $p$. It also predicts that (*) is equivalent to

\medskip

(**) The image of 
$$\Gal(\bar F/F) \to G(\Q_p)\;;\; \sigma\mapsto (\la_{1,\Q_p}^{-1}\circ \sigma \circ\la_{1,\Q_p}, \;\chi_{\text{cyclo}}(\sigma)^{-1})$$
is open in $G(\Q_p)$. Here $\la_{1,\Q_p}$ denotes the $\Q_p$-component of $\la_1$.

\end{sbpara}
In the rest of \S2.5, 
we fix $H_{0,\Z}$ as in \ref{HX(F)}. When we talk about $H_{\heartsuit,S}$, we fix $(L_i)_i$, $(L'_i)_i$, and $J$ as in \ref{Hheart}, and we assume that $K$ satisfies the condition (ii) in \ref{Hheart}.

 The part of the following Question 4 concerning (1) (resp. (2)) asks whether a motive version of Conjecture 24.1 (b) (resp. 24.3 (b)) of Vojta in \cite{Vo} is true.

\begin{sbpara}\label{Qvo}  {\bf Question 4.} Is the following statement true?

 Assume $\La$ is ample (\ref{mix100}). Fix  $d\geq 1$,  a finite set of places $S_0$ of $\Q$ which contains the archimedean place of $\Q$, and  fix $\epsilon\in \R_{>0}$. 
Then there exists $C\in \R_{>0}$ such that

\medskip

(1) \;\; $H_{\heartsuit, S}(M)\cdot |D_F|\geq C\cdot  H_{\spadesuit}(M)\cdot H_{\La}(M)^{-\epsilon}$ 

\medskip

(resp. (2) \;\; $(\prod_{v\in\Sig_S(M)} \sharp(\F_v))\cdot |D_F|\geq C\cdot  H_{\spadesuit}(M)\cdot H_{\La}(M)^{-\epsilon})$ 

\medskip
\noindent
 for any number field $F\subset \C$ such that $[F:\Q]\leq d$ and for any $M\in X_{\text{gen}}(F)$.

Here $S$ denotes the set of places of $F$ lying over $S_0$, $\Sig_S(M)$ denotes the set of all places of $F$ which do not lie over $S$ and at which $M$ has bad reduction,  and $D_F$ denotes the discriminant of $F$.

\end{sbpara}

\begin{sbrem} (1)
In the original conjecture of Vojta, an exceptional set is removed from the algebraic variety. The author hopes that putting the condition $M$ is generic  corresponds to it. A remark on this point is put in \ref{strem} (2).

\medskip

(2) Vojta proved that his conjectures 24.1 (b) and 24.3 (b) in \cite{Vo} are equivalent (\cite{Vo0}) though the latter is seemingly  stronger than the former.

In our motivic version, the answer Yes to the part (2) of Question 4 implies the answer Yes to the part (1) if (*) at the end of \ref{Hheart} is true.

\end{sbrem}

In the pure case, basing on \ref{pos}, we present the following question.

\begin{sbpara}\label{QVo3}

{\bf Question 5.} Is the following statement true?

Fix  $d$, $S_0$ and $\epsilon$ as in Question 4, and assume $\Phi$ is pure. 
Then there exists  $C\in \R_{>0}$ such that

\medskip

(1)\;\;  $H_{\heartsuit,S}(M)\cdot |D_F| \geq C\cdot  H_{\spadesuit}(M)^{1-\epsilon}
$

\medskip

(resp. (2) 
\;\; $(\prod_{v\in\Sig_S(M)} \sharp(\F_v))\cdot |D_F| \geq C\cdot  H_{\spadesuit}(M)^{1-\epsilon}
$)

\medskip
\noindent
for any number field $F\subset \C$ such that $[F:\Q]\leq d$ and for any $M\in X_{\text{gen}}(F)$.

\end{sbpara}

\begin{sbrem} In the case $g=1$ of \ref{Siegel}, 
 the answer Yes to Question 4 is essentially  the case of  Vojta conjectures 24.1 (b), 24.3 (b) of \cite{Vo} for modular curves.  The preprint \cite{Mo2} of Mochizuki proves this case of Vojta conjectures. As far as the author knows, the long referee work on this preprint is not finished yet. In \cite{Mo2} and in the related work \cite{Mo1}, Mochizuki studies heights of elliptic curves. Are his methods useful in the study of height functions for motives and Questions in this paper?

\end{sbrem}

\begin{sbpara}\label{Qweak} We do not put the condition generic in the following question.

{\bf Question 6.}  Assume $\Phi$ is pure and $\La$ is ample (\ref{mix100}).
 Is there $b\in \R_{>0}$ satisfying the following?

 If we fix $d$, $S_0$ as in Question 4,  there is $C\in \R_{>0}$ such that 
 
 \medskip
 
 (1) \;\; $H_{\heartsuit,S}(M)\cdot |D_F| \geq C\cdot  H_{\La}(M)^b$
 
 \medskip
 
(resp. (2)\;\;  $(\prod_{v\in\Sig_S(M)} \sharp(\F_v))\cdot |D_F|\geq C\cdot  H_{\La}(M)^b$)

\medskip
\noindent
for any number field $F\subset \C$ such that $[F:\Q]\leq d$ and for any $M\in X(F)$.

\end{sbpara}

\begin{sbpara} \label{Sha} If the answer to the part (i) $\Rightarrow$ (ii) in  Question 2 (\ref{Q1Q2}) is Yes and if the answer to either the part (1) or (2) to Question 6 is Yes, the following version of Shafarevich conjecture for pure motives is true. 

Recall that the original Shafarevich conjecture (the finiteness of abelian varieties over $F$ which are of good reduction outside $S$) was proved by Faltings \cite{Fa0}. 

\medskip

{\bf Conjecture.} Assume $\Phi$ is pure. Let $F\subset \C$ be a number field and let $S$ be a finite set of places of $F$ which contains all archimedean places of $F$. Then there are only finitely many $M\in X(F)$ which are of good reduction outside $S$.

\medskip
Another form of  Shafarevich conjecture for pure motives was  formulated in \cite{Ko1} 10.12.

\end{sbpara}

\begin{sbrem} In \cite{Vo}, Vojta compares height functions for varieties over  number fields with height functions in Nevanlinna theory. Imitating this, in a sequel (Part II) of this paper, we will compare the height functions of this paper with height functions $T_{f,\La}(r)$, $T_{f,\spadesuit}(r)$, $N_{f, \heartsuit}(r)$ of Nevanlinna type. Here $f$ is a morphism $\C\to \bar X(\C)$ from the complex plain $\C$ (not from a projective smooth curve) to the extended period domain $\bar X(\C)$ which induces a horizontal morphism $\C\smallsetminus R\to X(\C)$ for some  discrete subset $R$ of $\C$,  $\La$ is as in \ref{2.3.1}, and $r\in \R_{>0}$.

\end{sbrem}

\subsection{Speculations on the number of motives of bounded height}\label{2.5b}

\begin{sbpara}\label{BM} Let $V$, $\bar V$ and $D=\bar V\smallsetminus D$ be as in (2) in 0.1 and let $A$ be an ample line bundle on $\bar V$ and let $H=H_A$ be a height function on $V(F)$ associated to $A$.
 For $B\in\R_{>0}$, let 
$$N(H,B)= \sharp\{x\in V(F)\;|\; H(x) \leq B\}<\infty.$$
Then it is known that in many cases, we have
$$N(H,B) = c B^a \log(B)^b(1+o(1))$$
for some constants $a\in \Q_{\geq 0}$, $b\in \frac{1}{2}\Z$, $c\in \R_{\geq 0}$. There are many studies on these $a,b,c$. (Manin, Batyrev,Tschhinkel, Peyre, and other people). 

We wonder whether there are similar theories for  motives.

\end{sbpara}

\begin{sbpara}\label{NHB1}

Let $\Phi= ((h(w,r))_{w,r\in\Z}, H_{0,\Q}, W, (\langle\;,\;\rangle_{0,w})_{w\in \Z})$ be as in \ref{Phi} and fix an open compact subgroup  $K$ of $G({\bf A}^f_\Q)$ which satisfies the neat condition \ref{neat2}. 

Fix an ample $\La$ (\ref{mix100}) and let $$H:=H_{\La},\quad h:=h_{\La}.$$ 
Here and in the rest of \S2.6, we fix 
 $H_{0,\Z}$ as in \ref{HX(F)}.

\end{sbpara}

\begin{sbpara}\label{NHB2}  Fix a number field $F$.
For $B\in \R_{>0}$, 
let $$N(H,B)=\{M\in X(F)\;|\; H(M)\leq B\}, \quad 
N_{\text{gen}}(H, B)=\{M \in X_{\text{gen}}(F)\;|\; H(M)\leq B\}.$$
See \ref{generic} for the definition of $X_{\text{gen}}(F)$. 
We expect that these numbers are finite. 

\end{sbpara}

\begin{sbpara}\label{Q7}  {\bf Question 7.} Do we have 
$$N(H,B) = c B^a \log(B)^b(1+o(1)), \quad N_{\text{gen}}(H, B)= c' B^{a'}\log(B)^{b'}(1+o(1))$$
for some constants $a, a'\in \R_{\geq 0}$, $b, b'\in \frac{1}{2}\Z$, $c, c'\in \R_{\geq 0}$? Do we have $a,a'\in \Q$ if $c(w,r), t(w,d)\in \Q$ in the definition of $\La$?

\end{sbpara}
\begin{sbpara} Note that if $X(F)\neq \emptyset$ (resp. $X_{\text{gen}}(F)\neq \emptyset$) and the answer to Question 7 for $N(H,B)$ (resp. $N_{\text{gen}}(H, B)$) is Yes, $c>0$ (resp. $c'>0$) and the constants $c,a,b$ (resp. $c', a', b'$) are uniquely determined. 

\end{sbpara}

\begin{sbpara}\label{rEx1}  {\bf Example.} Fix an integer $r\geq 1$ and let $\Phi$ be as follows: $h(w,r)=1$ if $(w,r)$ is $(0,0)$ or $(-2r,-r)$, and  $h(w,r)=0$ for other $(w,r)$. $H_{0,\Q}=\Q\cdot (2\pi i)^r \oplus \Q$. Write $((2\pi i)^r, 0), (0,1)\in H_{0,\Q}$ as $e_1, e_2$, respectively. $W_{-2r-1}=0\subset \Q\cdot (2\pi i)\oplus 0=W_{-2r}=W_{-1}\subset W_0=H_{0,\Q}$. $\langle e_1, e_1\rangle_{0,-2r}=(2\pi i)^{2r}$, $\langle e_2, e_2\rangle_{0,0}= 1$. 

Then (for the matrix presentation with respect to the basis $(e_1,e_2)$)
$$G= \{(g,t)\;|\; g\in GL_2, \; t\in {\bf G}_m, \; g_{11}^2=t^{-2r}, g_{21}=0, g_{22}^2=1\}.$$
Fix $n\geq 3$ which is a power of a prime number,  and let 
$K\subset G({\bf A}^f_\Q)$ be the subgroup consisting of all elements $(g,t)$ such that 
$$g_{11}, g_{22}\in \text{Ker}(\hat \Z^\times \to (\Z/n\Z)^\times), \; g_{12}\in \hat \Z, \; t\in \hat \Z^\times.$$
Then $K$ satisfies the neat condition \ref{neat2}. 
We have
$$X(\C) = \C/(\Z\cdot (2\pi i)^r)$$
Here $z\in \C$ corresponds to $\text{class}(f(z), 1, 1)\in X(\C)$ where $f(z)\in D= D^{\pm}$ is the decreasing filtration on $H_{0,\C}$ defined by: $f(z)^0=H_{0,\C}$, $f(z)^1=f(z)^r$ is generated by $(z,1)$, and $f(z)^{r+1}=0$. 

Let $F\subset \C$ be a number field and assume that $F$ contains a primitive $n$-th root $\zeta_n$ of $1$. Then $X(F)$ 
is identified with the group of extensions $\text{Ext}^1(\Z, \Z(r))$  in the category of mixed motives over $F$ with $\Z$-coefficients. 
Let $X'(F)$ be the fiber product of $\Q\otimes K_{2r-1}(F) \to H^1_{\text{cont}}(\Gal(\bar F/F), {\bf A}^f_\Q(r)) \leftarrow H^1_{\text{cont}}(\Gal(\bar F/F), \hat \Z(r))$. Then we have a canonical homomorphism $X'(F) \to X(F)$ and 
the usual philosophy of motives suggests that this homomorphism is an isomorphism. 

We define the height functions $H$ and $h$ (\ref{NHB1}) by taking $t(0,2r)=1$. That is, $H=H_{\star\diamondsuit}=H_{\diamondsuit,0,2r}$, $h=h_{\star\diamondsuit}=h_{\diamondsuit,0,2r}$.  We will consider $N(H, B)$ and $N_{\text{gen}}(H, B)$ replacing $X(F)$ by $X'(F)$ and replacing $X_{\text{gen}}(F)$ by the inverse image of $X'_{\text{gen}}(F)$ in $X'(F)$. In the case $r\geq 2$, we will assume  $F=\Q(\zeta_n)$ with $n=3$ or $n=4$ for simplicity. 

Assume $r=1$. Then $X'(F)= F^\times$. An element of $F^\times$ belongs to $X'_{\text{gen}}(F)$ if and only it is not a root of $1$. By \ref{1.7ex1}, we have
$$N(H, B)= c \cdot B(1+o(1)),     \quad N_{\text{gen}}(H, B)= c\cdot B(1+o(1))$$
for the set $X'(F)$ for some  $c\in \R_{>0}$.

Assume $r\geq 2$ and $F=\Q(\zeta_n)$ with $n=3$ or $n=4$. We show that
$$N(H, B)= c \cdot \zeta_F(r)^{-1}\log(B)^{r/2}(1+o(1)), \quad N_{\text{gen}}(H, B)= c \cdot \zeta_F(r)^{-1}\log(B)^{r/2}(1+o(1))$$
for $X'(F)$ for some  $c\in \Q_{>0}$ where $\zeta_F(s)$ is the Dedekind zeta function of $F$. 

It is known that $X'(F) \cong \Z \oplus J$ for some finite abelian group $J$ and that the regulator of the element of $X'(F)$ corresponding to $(1,0)\in \Z\oplus J$ in
$\text{Ext}^1_{\R\text{MHS}}(\R, \R(r))\cong \R$ is  $u\cdot \zeta_F(r)$ for some $u\in \Q_{>0}$. Hence for the element $M\in X'(F)$ corresponding to $(m, x)\in \Z\oplus J$, 
we have  $H(M) =\exp((2u|m|\zeta_F(r))^{2/r})$. $M$ belongs to $X'_{\text{gen}}(F)$ if 
and only if $m\neq 0$. Hence $N(H, B)$ (resp. $N_{\text{gen}}(H,B)$) is $\sharp(J)$ times the number of 
integers (resp. non-zero integers) $m$ such that $|m|\leq (2u)^{-1}\zeta_F(r)^{-1}\log(B)^{r/2}$.

\end{sbpara}

\begin{sbpara}\label{beta} Define $\beta\in \R_{\geq 0}\cup\{\infty\}$ as 
$$\beta:= \lim\sup\; \frac{\log(N_{\text{gen}}(H,B))}{\log(B)}\quad (B\to \infty).$$

If $X_{\text{gen}}(F)\neq \emptyset$ and the answer to Question 7 is Yes,  $\beta$ is $a'$ in \ref{Q7}. 

In conjectures in \cite{BM}, the constant $a$ in \ref{BM} (defined  using $\lim\sup$ as above, and denoted by $\beta$ in \cite{BM}) 
 is related to a geometric invariant $\alpha$ which is defined by using a relation between the ample divisor $A$ and the canonical divisor $K$ of $\bar V$. Following the analogy, we will define a geometric invariant $\alpha$ using the height functions $h_{\spadesuit}$ and $h_{\heartsuit}$ for variations of Hodge structure, and compare it in Question 8 with the above $\beta$. In the analogy, our height function $h_{\spadesuit}$ plays the role of $H_{K+D}$, $h_{\heartsuit}$ plays the role of  $H_D$, and $h_{\spadesuit}- h_{\heartsuit}$ plays the role of $H_K=H_{K+D}H_D^{-1}$ (\ref{handper}).

\end{sbpara}

\begin{sbpara}\label{alpha}
Define 
$\alpha\in \R_{\geq 0}\cup \{\infty\}$ to be the inf of all $s\in \R_{\geq 0}$ satisfying the following condition (*).

(*) For any projective smooth curve $C$ over $\C$ and for any $\sH\in \Mh(C, \bar X(\C))$, we have $$sh(\sH) + h_{\spadesuit}(\sH) - h_{\heartsuit}(\sH)\geq 0$$
\medskip
($\alpha$ is defined to be $\infty$ if such $s$ does not exist).

\medskip

{\bf Question 8.}  Assume $X_{\text{gen}}(F)\neq \emptyset$. Do we have $\alpha,\beta<\infty$ always?
Do we have 
$$\alpha=\beta?$$ 
\end{sbpara}

For the reason why we consider only generic motives in Question 8, see \ref{strem} (1). A remark on this point is put in \ref{strem} (2). 

\begin{sbpara}\label{rEx2} {\bf Example.} Let the notation be as in \ref{rEx1}. 

Assume $r=1$. Then $\beta=1$ by \ref{rEx1}. In this case, $X'(F) \to X(\C) \subset \bar X(\C)$ is identified with $F^\times \to \C^\times \subset {\bf P}^1(\C)$, $h(\sH)$ and $h_{\heartsuit}(\sH)$ coincide with the degree of the pull back of the divisor $\{0,\infty\}$  on ${\bf P}^1(\C)$ under the  morphism $C\to {\bf P}^1(\C)$ of $\sH$, and $h_{\spadesuit}=0$. Hence $\alpha=1$. Thus $\alpha=\beta$.

Assume $r\geq 2$ and $F=\Q(\zeta_n)$ with $n=3$ or $n=4$. Then $\beta=0$ by \ref{rEx1}. We have $h_{\spadesuit}=0$ and $h_{\heartsuit}=0$. Hence $\alpha=0$. Thus $\alpha=\beta$. 

\end{sbpara}

\begin{sbpara}\label{Q10.0} The following Question 9 is a refined version of Question 8. In Question 9, we consider the type of monodromy. 

Consider the quotient set $\Sig/\sim$ of $\Sig$ where $\sim$ is the following equivalence relation. For nilpotent operators $N, N'\in \Lie(G_1)$, $\R_{\geq 0}N\sim \R_{\geq 0}N'$ if and only if there are $(g,t)\in G(\C)$ and $c\in \C^\times$ such that $N'=c\text{Ad}(g)(N)$ in $\Lie(G_1)_\C$. Let $\Sig'$ be a subset of $\Sig/\sim$ which contains the class of the cone $\{0\}$.

\medskip
We define the height function $h_{\heartsuit, \Sig'}(\sH)$ and the counting function $N_{\text{gen},S, \Sig'}(H, B)$ for a finite set $S$ of places of $F$ containing all archimedean places of $F$, which take the shapes of the monodromy operators into account.

For  a projective smooth curve $C$ and for $\sH\in \Mh(C, \bar X(\C))$, 
define $h_{\heartsuit, \Sig'}(\sH)= \sum_x  h_{\heartsuit, x}(\sH)\in \Z_{\geq 0}$ where $x$ ranges over all points of $C$ such that $N'_x$ satisfies the following condition (i).

\medskip

(i) The image of $x$ under the composite map $$C\overset{\sH}\to \bar X(\C) \to G(\Q)\bs (\Sig \times (G({\bf A}^f_\Q)/K))\to G(\Q)\bs \Sig \to \Sig/\sim$$ belongs to $\Sig'$. Here $G(\Q)$ acts on $\Sig$ by conjugation.

\medskip

Let $N_{\text{gen}, S, \Sig'}(H, B)$ be the number of $M\in X_{\text{gen}}(F)$ which appear in the definition of $N_{\text{gen}}(H,B)$ and which satisfy  the following condition (ii) at each place $v\notin S$ of $F$.

(ii) There is an element $\text{class}(\R_{\geq 0}N)$ (with $N\in \Lie(G_1)$ nilpotent) 
of $\Sig'$ which satisfies the following:
Let $p=\text{char}(\F_v)$. For any prime number $\ell\neq p$, there are $t,c\in \bar \Q_{\ell}^\times$ and  an  isomorphism $\bar \Q_{\ell}\otimes_{\Q_{\ell}} M_{et,\Q_{\ell}}\cong H_{0,\bar \Q_{\ell}}$ of $\bar \Q_{\ell}$-vector spaces which preserves $W$ and via which $\langle\;,\;\rangle_w$ on $\bar \Q_{\ell} \otimes_{\Q_{\ell}} \gr^W_wM_{et,\Q_{\ell}}$ corresponds to $t^w\langle \; , \;\rangle_{0,w}$ on $\gr^W_wH_{0,\bar \Q_{\ell}}$ for any $w\in \Z$ and  
the local monodromy operator $N'_v$ corresponds to $cN$.   Furthermore, there are $t,c\in \bar F_v^\times$ and  an  isomorphism $\bar F_v \otimes_{F_{v,0,\text{ur}}} D_{\pst}(F_v, M_{et,\Q_p})\cong H_{0, \bar F_v}$ of $\bar F_v$-vector spaces 
which preserves $W$ and via which  $\langle\;,\;\rangle_w$ on $\bar F_v\otimes_{F_{v,0,\text{ur}}} \gr^W_wD_{\pst}(F_v, M_{et,\Q_p})$ corresponds to $t^w\langle \; , \;\rangle_{0,w}$ on $\gr^W_wH_{0,\bar F_v}$  for any $w\in \Z$ and 
the local monodromy operator $N'_v$ corresponds to $cN$.

\end{sbpara}

\begin{sbpara}　The above $h_{\heartsuit,\Sig'}$ can be expressed 
in the same form 
$$h_{\heartsuit,\Sig'}(\sH)= \sum_{x\in C} e(x), \quad (f^*I_{\Sig'})_x=m_x^{e(x)}$$
as in \ref{handper} (2), where $I_{\Sig'}$ is the invertible ideal of $\sO_{\bar X(\C)}$ which coincides with $I_{\bar X(\C)}$ in \ref{idealI} at points of  $\bar X(\C)$ whose classes in $\Sig/\sim$ belong to $\Sig'$ and coincides with $\sO_{\bar X(\C)}$ at  other points of $\bar X(\C)$. 
 
 \end{sbpara}
 
\begin{sbpara}\label{Q10}  Let the notation be as in \ref{Q10.0}. Let $F$ be a number field. Let $S$ be a finite set of places of $F$ containing all archimedean places. 
Define the modified versions $$\beta_{S,\Sig'}, \;\; \alpha_{\Sig'} \;\;   \quad \text{of}\quad  \beta, \;\;  \alpha$$ by replacing  $$N_{\text{gen}}(H, B)\quad \text{by}\quad N_{\text{gen}, S, \Sig'}(H,B)$$  and replacing $$sh(\sH)+ h_{\spadesuit}(\sH) -h_{\heartsuit}(\sH)\geq 0\quad \text{by} \quad sh(\sH)+ h_{\spadesuit}(\sH) - h_{\heartsuit,\Sig'}(\sH)\geq 0,$$ respectively. 

\medskip
{\bf Question 9}. Assume $X_{\text{gen}}(F)\neq \emptyset$. Do we have $\alpha_{\Sig'}= \beta_{S,\Sig'}$?

\end{sbpara}

\begin{sbpara}\label{rEx3}  {\bf Example.} Let the notation be as in \ref{rEx1}. Consider the case $\Sig'$ consists just of the class of the cone $\{0\}$ and consider the case $r=1$. We have $\alpha_{\Sig'}=0$ because $h_{\heartsuit, \Sig'}=0$. 
On the other hand, it is easy to prove that $\beta_{S,\Sig'}=0$ for the set $X'(F)$. Thus $\alpha_{\Sig'}=\beta_{S,\Sig'}$. 
\end{sbpara}

\begin{sbrem} Assume that the conjecture (i) in \ref{WDconj} is true and that the weight monodromy conjecture (\ref{Mdgeq2}) on $M$ at any $v\notin S$ is true. Then for the element $N$ in \ref{WDconj}, $\R_{\geq 0}N\in \Sig$ and the class of $(\R_{\geq 0}N, ({\bf A}^f_\Q\otimes \theta_1)\circ \la_1, \theta_2\la_2)$ in $G(\Q) \bs (\Sig\times (G({\bf A}^f_\Q)/K))$ is well defined. Hence we can formulate a generalization (a finer version) of Question 9 for a subset $\Sig'$ of $G(\Q)\bs (\Sig \times (G({\bf A}^f/\Q)/K))$ containing $G(\Q) \bs (\{0\} \times (G({\bf A}^f_\Q)/K))$. 

\end{sbrem}

\begin{sbrem}\label{strem}

(1) In  the motivic version of Batyrev-Manin conjecture, we consider only generic motives. In the original Batyrev-Manin conjecture for (2) in 0.1, some closed subset of the algebraic variety $V$, which has too many $F$-rational points, is removed if necessary. The author hopes that putting the condition generic corresponds to it. 

\medskip

(2) In \S\ref{2.4b} and \S\ref{2.5b}, we considered mainly generic motives (in the sense of this paper), that is, motives whose 
 Mumford-Tate groups are of finite index in  $G$. In a sequel (Part II) of this paper, we will treat more general mixed motives. To treat a motive whose Mumford-Tate group is any linear algebraic group $G$ over $\Q$, we will consider the height functions for $G$-motives over a number field $F$. Here a $G$-motive over $F$ means an exact tensor functor from the category of finite-dimensional representations of $G$ over $\Q$ to the category of mixed motives with $\Q$-coefficients over $F$. Just as we removed motives whose Mumford-Tate group $H$ are of smaller dimensions in \S\ref{2.4b} and \S\ref{2.5b}, we will remove $G$-motives which are indued from $H$-motives for algebraic subgroups $H$ of $G$ over $\Q$ such that $\dim(H) < \dim(G)$. 
\end{sbrem} 

\begin{sbrem}\label{BBK} We can also formulate the problem to count the numbers of mixed motives whose pure graded quotients are fixed. This is closely related to Beilinson conjectures in \cite{Be} and Tamagawa number conjectures in \cite{BK} on zeta values. Let $F$ be a number field. Assume that for each $w\in \Z$, we are given a pure motive $M_w$ with $\Z$-coefficients over $F$ of weight $w$. Assume $M_w=0$ for almost all $w$. Let $\sX(F)$ be the set of all isomorphism classes of mixed motives $M$ with $\Z$-coefficients over $F$ such that $\gr^W_wM=M_w$ for all $w$.
For $B\geq 0$, let $N(B)=\sharp\{M\in \sX(F)\;|\; H_{\diamondsuit}(M)\leq B\}$.

(1) The author imagines that it may be possible to deduce $N(B)=O(\log(B)^n)$ for some $n\geq 0$ from the conjectures in \cite{Be} and \cite{BK} which appear in \S2.1, but it may be harder to show that $N(B)= c\log(B)^m(1+o(1))$ for some $m\geq \Z_{\geq 0}$ and $c\in \R_{>0}$, generalizing 2.6.5. Assuming  such $c,m$ exist, 
we ask how this $c$ is expressed by zeta values.

(2) Assume that $M_0=\Z$ and that there is  $w<0$ such that  $M_{w'}=0$ for $w'\neq 0,w$. Then $\sX(F)= \text{Ext}^1(\Z, M_w)$ is a group. In \cite{BK}, the adelic period domain $\MyProd_v \sX(F_v)$ was defined by using the Ext groups $\sX(F_v)$ in the category of mixed Hodge structures for $v$ archimedean and in the category of $p$-adic Hodge structures for $v$ non-archimdean, and the Tamagawa measure $\mu_{\diamondsuit}$ on $\MyProd_v \sX(F_v)$ was defined. For simplicity, assume $w<-2$. Then Tamagawa number conjecture in \cite{BK} predicts that the volume $\mu_{\diamondsuit}(\sX(F)\bs \MyProd_v \sX(F_v))$ of the quotient group $\sX(F) \bs \MyProd_v \sX(F_v)$ 
 is a rational number (this rationality is a reformulation of Beilinson conjecture in \cite{Be} on zeta values) which is related to Tate-Shafarevich group of $M_w$. 

(3) It is interesting how this (2) can be extended  to the situation where the order of the set $\{w\;|\; M_w\neq 0\}$ can be $\geq 3$. 
The author plans to define the adelic period domain $\MyProd_v \sX(F_v)$ and the Tamagawa measure $\mu_{\diamondsuit}$ on $\MyProd_v \sX(F_v)$ in another paper. Since $\sX(F)$ and $\sX(F_v)$ are not groups any more, we can not consider the quotient space $\sX(F)\bs \MyProd_v \sX(F_v)$ but we can consider how many points of $\sX(F)$ are distributed inside $\MyProd_v \sX(F_v)$ with respect to the Tamagawa measure. 
For simplicity, assume that if $a>b$ and if $M_a\neq 0$ and $M_b\neq 0$, then $a-b>2$. The height function $H_{\diamondsuit}:\sX(F) \to \R_{>0}$ is induced from a continuous map $H_{\diamondsuit}: \MyProd_v \sX(F_v)\to \R_{>0}$ which is the product of continuous maps $H_{\diamondsuit,v}:\sX(F_v)\to \R_{>0}$. The set $C(B):= \{x\in \MyProd_v \sX(F_v)\;|\;H_{\diamondsuit}(x) \leq B\}$ is compact and of finite volume for $\mu_{\diamondsuit}$. The problem in (1) should be extended to the problem of the study of the limit of the ratio $(N(B): \mu_{\diamondsuit}(C(B)): \log(B))$ when $B\to \infty$. N. Tung showed that in (2),  $\mu_{\diamondsuit}(\sX(F)\bs \MyProd_v \sX(F_v))$  is equal to  the limit $\lim_{B\to \infty} \mu_{\diamondsuit}(C(B))/N(B)$.

Thus Beilinson  conjectures in \cite{Be} and Tamagawa number conjectures in \cite{BK} on zeta values may be generalized to the problem to count the number of mixed motives with fixed pure graded quotients with bounded height, which may be understood as a motive version of the work of Peyre \cite{P} who studied $c$ in \ref{BM} for Fano varieties $V$ by using Tamagawa measures on the adelic spaces $\MyProd_v V(F_v)$. 
\end{sbrem}

\end{document}